\documentclass[reqno,12pt]{amsart}
\usepackage{amssymb,amsmath,amsthm,amstext,amsfonts,euscript}
\usepackage{graphicx}
\usepackage{a4wide}
\usepackage{color}
\usepackage{psfrag}
\usepackage[shortlabels]{enumitem} 

\usepackage{hyperref}

\makeatletter \@addtoreset{equation}{section} \makeatother

\theoremstyle{plain}
\newtheorem{maintheorem}{Theorem}
\newtheorem{maincorollary}[maintheorem]{Corollary}

\newtheorem{theorem}{Theorem}[section]
\newtheorem{proposition}[theorem]{Proposition}
\newtheorem{lemma}[theorem]{Lemma}

\newtheorem{corollary}[theorem]{Corollary}
\theoremstyle{definition} \theoremstyle{remark}
\newtheorem{remark}[theorem]{Remark}
\newtheorem{definition}[theorem]{Definition}
\newtheorem{conjecture}{Conjecture}
\newtheorem{example}{Example}

\newcommand{\diam}{\operatorname{diam}}
\newcommand{\dist}{\operatorname{dist}}
\newcommand{\supp}{\operatorname{supp}}

\newcommand{\EE}{{\mathbb E}}
\newcommand{\RR}{{\mathbb R}}
\newcommand{\NN}{{\mathbb N}}
\newcommand{\CC}{{\mathbb C}}

\newcommand{\PP}{{\mathbb{P}}}

\newcommand{\vfi}{\varphi}

\newcommand{\cB}{{\mathcal B}}
\newcommand{\cC}{{\mathcal C}}
\newcommand{\cD}{{\mathcal D}}
\newcommand{\cF}{{\mathcal F}}
\newcommand{\cI}{{\mathcal I}}
\newcommand{\cL}{{\mathcal L}}
\newcommand{\cM}{{\mathcal M}}
\newcommand{\cN}{{\mathcal N}}
\newcommand{\cQ}{{\mathcal Q}}
\newcommand{\cR}{{\mathcal R}}
\newcommand{\cS}{{\mathcal S}}
\newcommand{\cW}{\mathcal{W}}
\newcommand{\cP}{\mathcal{P}}
\newcommand{\cO}{\mathcal{O}}
\newcommand{\cU}{\mathcal{U}}
\newcommand{\cV}{\mathcal{V}}

\newcommand{\cH}{\mathcal{H}}

 \newcommand{\inter}{\operatorname{int}}

\newcommand{\Lip}{\operatorname{Lip}}
\newcommand{\Leb}{\operatorname{Leb}}
\newcommand{\leb}{\operatorname{Leb}}

\DeclareMathOperator{\Div}{div}
\DeclareMathOperator{\sgn}{sgn}

\newcommand{\loc}{\text{loc}}

 \def \fX {{\mathfrak X}}

\begin{document}

\title{Robust Exponential Mixing and Convergence to Equilibrium for Singular-Hyperbolic Attracting Sets}

\author{Vitor Ara\'ujo}

\address[V.A.]{
 Departamento de Matem\'atica, Universidade Federal da Bahia\\
 Av. Ademar de Barros s/n, 40170-110 Salvador, Brazil.}

\email{vitor.d.araujo@ufba.br, vitor.araujo.im.ufba@gmail.com}
\urladdr{https://sites.google.com/view/vitor-araujo-ime-ufba}

\keywords{singular-hyperbolic attracting set, physical/SRB measures,
  robust exponential mixing, exponential convergence to equilibrium}

\subjclass[2010]{Primary: 37D25. Secondary: 37D30, 37D20, 37D45.}

\author{Edvan Trindade}

\address[E.T.]{Instituto Federal de Educa\c{c}\~ao Ci\^encia e Tecnologia da
  Bahia (IFBA)\\Campus Porto Seguro\\Rod. Br 367 Km 57,5 - Fontana I,
   45810-000, Porto Seguro - BA, Brazil}

\email{edvan.trindade@ifba.edu.br, trindade.matematica@gmail.com}

\thanks{V.A. was partially supported by CNPq-Brazil
and E.T. was partially supported by CAPES-Brazil.}

  \begin{abstract}
We extend results on robust exponential mixing for geometric Lorenz
attractors, with a dense orbit and a unique singularity, to
singular-hyperbolic attracting sets with any number of (either Lorenz- or
non-Lorenz-like) singularities and finitely many ergodic physical/SRB
invariant probability measures, whose basins cover a full Lebesgue
measure subset of the trapping region of the attracting set.

We obtain exponential mixing for any physical probability measure
supported in the trapping region and also exponential convergence to
equilibrium, for a $C^2$ open subset of vector fields in any
$d$-dimensional compact manifold ($d\ge3$).

\end{abstract}

\date{\today}

\maketitle
\tableofcontents

 \section{Introduction}  \label{sec:intro}

The expression ``statistical properties'' of a Dynamical System refers
to the statistical behavior of \emph{typical} trajectories of the
system. It is well-known that this behavior is related to properties
of the evolution of measures under the dynamics. Statistical
properties are frequently a simpler object to study than pointwise
behaviour of trajectories, which is most of the time
unpredictable. However, statistical properties are regular for most
known systems and mostly admit a simple description.

The statistical tools provided by Differentiable Ergodic Theory are
among the most powerful techniques available to study the global
asymptotic behavior of Dynamical Sytems. A central concept is that of
\emph{physical} measure (or \emph{Sinai-Ruelle-Bowen} measure) for a
flow or a transformation. Such measure for a flow $X_t$ on a compact
manifold is an invariant probability measure $\mu$ for which the
family of points $z$ satisfying
\begin{align*}
  \lim_{t\to+\infty}\frac1t\int_0^t \psi\big(X_s(z)\big) \,
  ds = \int\psi\,d\mu,
\end{align*}
for all continuous observables (functions) $\psi$. That is, the time
averages of a continuous observable along the trajectory of $z$
converge to the space average of the same observable with respect to
$\mu$; and the set of all these points is a positive Lebesgue
(volume) measure subset $B(\mu)$ (the \emph{ergodic basin}) of the
ambient space.

These time averages are considered  \emph{a priori} physically
observable when dealing with a mathematical model of some real
phenomenon whose properties can be measurable.

This kind of measures was first rigorously obtained for (uniformly)
hyperbolic diffeomorphisms by Sinai, Ruelle and Bowen
\cite{Si72,Ru76,BR75}. For non uniformly hyperbolic transformations
and flows these measures were studied more recently: we mention only
the results closer to the present text in \cite{APPV, LeplYa17,
  ArSzTr}, on the existence of physical measures for
singular-hyperbolic attractors.
Statistical properties of such measures are an active field of study:
among the article used in this work we stress
\cite{alves-luzzatto-pinheiro2005,gouezel,ArVar,ArBuVa,ArMel18,AMV15}.

The general motivation is that the family $\{ \psi\circ X_t\}_{t\ge0}$
should behave asymptotically as a family of independent and
identically distributed random variables.

An important property is the speed of convergence of the time average
to the space average among many others.
Considering $\vfi$ and $\psi\circ X_t:M\to\RR$ as random variables
with law $\mu$, \emph{mixing} means that  \emph{the random variables
    $\vfi$ and $\psi\circ X_t$ are asymptotically 
  independent}: $\EE\big(\vfi\cdot(\psi\circ X_t)\big)$ converges
to $\EE(\vfi)\cdot\EE(\psi)$ when $t$ grown without bound.
Writting the correlation function
\begin{align*}
C_t(\vfi,\psi)
&=
\EE\big(\psi\cdot(\vfi\circ X_t)\big)
- \EE(\vfi)\cdot\EE(\psi)
=
\int\psi\cdot\big(\vfi\circ X_t\big) \, d\mu
-\int\vfi\,d\mu\int \psi\, d\mu
\end{align*}
we get
 $|C_t(\vfi,\psi) | \xrightarrow[t\to\infty]{}0$ for all integrable
 observables in case of mixing.
 \emph{Exponential mixing} means that there exist
 $ C, \gamma>0$ so that
$$
|C_t(\psi,\varphi)| \le C e^{-\gamma t} \|\psi\| \|\varphi\|, \quad
t>0;
$$
while \emph{superpolynomial mixing} holds if for all $\beta>0$ we can
find $C_\beta>0$ for which
    $$
    |C_t(\psi,\varphi)| \le C_\beta t^{-\beta} \|\psi\|
        \|\varphi\|, \quad t>0;
    $$
on a Banach space of usually more regular observables than just
integrable ones (mostly H\"older continuous, some times differentiable).

To ascertain the speed of mixing is a subtle issue for flows. In spite
of exponential mixing having been prived for hyperbolic
diffeomorphisms for Sinai, Ruelle and Bowen \cite{Si72,Ru76,BR75} in
the 70's, only in the final years of the XXth century a significant
breakthrough was obtained in the fundamental work of Dolgopyat
\cite{Do98}. Here the author obtained for the first time exponential
mixing for Anosov flows with respect to physical measures under rather
strong assumptions (global smoothness of stable and unstable
foliations and their uniform non-integrability). These assumptions are
not robust, i.e., the family of systems which satisfy these
assumptions loose these properties by small perturbations.

Later superpolynomial mixing was obtained for  open and dense families
(hence robust) of hyperobolic flows by Field, Melbourne and
Torok \cite{FMT} refining Dolgopyat techniques, but only achieving a
slower mixing speed.

Singular-hyperbolicity is a non-trivial recent extension of the
notion of uniform hyperbolicity that encompassses systems like the
Lorenz attractor in a unified theory, founded on the work of Morales,
Pacifico and Pujals~\cite{MPP04}. This allows to rigorously frame
Lorenz-like attractors after the the work of Tucker \cite{Tu99}.

For singular-hyperbolic attracting sets the existence of physical
measures and some of their properties were obtained for the first time in
\cite{APPV}. Surprinsingly it was easier to obtain robust exponential mixing
for physical measures among Lorenz-like attractors -- this was first
proved by Araujo and Varandas in \cite{ArVar} for an open subset of
vector fields with a geometric Lorenz attractor -- than among
hyperoblic attractors or even Anosov flows.

For the original Lorenz attractor exponentially mixing was proved by
the works of Araujo, Melbourne and Varandas \cite{AMV15,ArMel16} and
recently Araujo and Melbourne \cite{ArMel18} proved superpolynomial
mixing for an open an dense subset of singular-hyperbolic attracting
sets.

The same techniques allow us to obtain robust exponential mixing for
Axiom A attractors \cite{ArBuVa} and have been recently extended
to achieve robust exponential mixing for Anosov flows
\cite{BuWar20}. Still more recently \cite{tsujii2020zhang} explores
the same technique to get exponential mixing for all equilibrium
states (of which physical measure are but an example) with respect to
H\"older continuous potentials for an open and dense subset of
topologically mixing $C^\infty$ Anosov flows on $3$-manifolds.

An interesting variation of the theme is the \emph{convergence to
  equilibrium}: replacing $\mu$ by the Lebesgue (volume) measure we
consider the following function
  \begin{align*}
    E_t(\psi,\varphi) = \int (\varphi\circ
      X_t)\,\psi\; d\Leb - \int \psi\; d\Leb \int
      \varphi\; d\mu
  \end{align*}
  and if we have convergence $|E_t(\psi,\varphi)|\to0$ this means, in
  particular (letting $\psi\equiv1$), that for a certain (usually fairly
  regular) class of observables
  \begin{align*}
    \lim_{t\to\infty}\int \vfi\circ X_t \,d\Leb = \int\vfi\,d\mu
  \end{align*}
  which, \emph{a priori}, allows us to use a ``natural'' measure to
  estimate $\mu$ through experimental observations of the system.

  In this work we extend the result of robust exponential mixing for
  the $3$-dimensional geometric Lorenz attractor, with a unique
  singularity and a dense orbit, to singular-hyperbolic attracting
  set, with any number of singularities (Lorenz-like or not), finite
  number of invariant ergodic physical probability measures and higher
  dimensional stable bundle. We obtain exponential mixing for all
  physical measures supported on the trapping region of the attracting
  set and also exponential convergence to equilibrium, for a
  $C^2$-open subset of vector fields on compact $d$-manifold
  ($d\ge3$).


\subsection{Preliminary definitions}
\label{sec:main-results}

Let $M$ be a compact boundaryless $d$-dimensional manifold. Given an
integer $k \geq 1$, we denote by $\fX^k(M)$ the set of $C^k$ vector
fields on $M$ endowed with the $C^k$ topology. We fix some smooth
Riemannian structure on $M$ and we denote the distance induced by this
structure by $\dist$ and the volume measure by $\leb$. We may assume
that both $\dist$ and $\leb$ are normalized, that is, the diameter of
$M$, denoted here by $\diam(M)$, and $\leb(M)$ are equal to $1$.

Given $X\in{\fX}^k(M)$ we denote by $X_t:M \to M$, $t \in \RR$, the flow
induced by $X$. For each $x \in M$ and each interval $I \subset \RR$ we set
$X_I(x) := \{X_t(x):\ t \in I\}$. In general, given a point $x \in M$ we
denote the \emph{orbit} of $x$ by the flow of $X$ by the set $\cO_X(x) =
X_{\RR}(x)$.

We say that $x \in M$ is \emph{regular} for the vector field $X$ if
$X(x) \neq 0$. Otherwise we say that $x$ is an \emph{equilibrium} or
\emph{singularity} of $X$. We also say that the corresponding orbit is
\emph{regular} or \emph{singular}, respectively. If $\sigma\in M$ is a
singularity for $X$ then $\sigma$ is a fixed point for the flow of
$X$, that is, $X_t(\sigma) = \sigma$ for all $t \in \RR$. We say that
$p\in M$ is a \emph{periodic point} (or the orbit of $p$ is
\emph{periodic}) for $X$, if the set $\{t \in \RR^+:\ X_t(p) = p\}$ is
nonempty and the number $T := \inf\{t \in \RR^+:\ X_t(p) = p\}$ is
positive. In this case we call $T$ the \emph{period} of $p$.

We say that a set $\Lambda \subset M$ is \emph{invariant} by $X$ if
$X_t(\Lambda) = \Lambda$ for all $t \in \RR$. A compact invariant set
$\Lambda$ for $X$ is said to be \emph{isolated} if we can find an open
neighborhood $U\supset\Lambda$ so that
$\Lambda =\bigcap_{t\in\RR}X_t(U)$. If $U$ also satisfies
$\overline{X_t(U)}\subset U$\footnote{We write $\overline{A}$ to
  denote the topological closure of a set $A$.} for all $t>0$ then we
say that $\Lambda$ is an \emph{attracting set} and that $U$ is a
\emph{trapping region} for $\Lambda$. In this case we have that
$\Lambda = \cap_{t > 0}X_t(U)$.  The \emph{topological basin} of an
attracting set $\Lambda$ is given by
$$W^s(\Lambda)= \left\{ x\in M : \lim_{t\to+\infty}\dist\big(
  X_t(x) , \Lambda\big) = 0 \right\}.$$

Given $x \in M$ the \emph{$\omega$-limit set} of $x$ by the flow $X_t$
is given by the set
$$\omega(x)=\omega_X(x) = \left\{y \in M:\ \exists t_k \nearrow +\infty\text{
    such that } \lim_{k \to +\infty}\dist(X_{t_k}(x), y) =
  0\right\}.$$

An invariant set $\Lambda$ is \emph{transitive} for $X$ if there
exists a regular point $x \in M$ such that $\Lambda = \omega_X(x)$. We
say that $\Lambda$ is \emph{non-trivial} if it is neither a finite set
of periodic orbits nor a finite set of equilibria. Otherwise we say
that $\Lambda$ is \emph{trivial}.

A compact invariant set $\Lambda \subset M$ is an \emph{attractor} for
a vector field $X$ if it is a transitive attracting set for $X$. We
say that the attractor is \emph{proper} if it is not the whole ambient
manifold $M$.

\subsubsection{Singular-hyperbolic attracting sets}
\label{sec:PH}

Let $\Lambda$ be a compact invariant set for $X\in \fX^r(M)$ for some
$r\ge1$.  We say that $\Lambda$ is {\em partially hyperbolic} if the
tangent bundle over $\Lambda$ can be written as a continuous
$DX_t$-invariant sum
$$
T_\Lambda M=E^s\oplus E^{cu},
$$
where $d_s=\dim E^s_x\ge1$ and $d_{cu}=\dim E^{cu}_x\ge2$ for $x\in\Lambda$,
and there exist constants $C>0$, $\lambda\in(0,1)$ such that
for all $x \in \Lambda$, $t\ge0$, we have
\begin{itemize}
\item uniform contraction along $E^s$:
$
\|DX_t | E^s_x\| \le C \lambda^t;
$
\item domination of the splitting:
$
\|DX_t | E^s_x\| \cdot \|DX_{-t} | E^{cu}_{X_tx}\| \le C \lambda^t.
$
\end{itemize}
We refer to $E^s$ as the stable bundle and to $E^{cu}$ as the
center-unstable bundle.  A {\em partially hyperbolic attracting set}
is a partially hyperbolic set that is also an attracting set.

The center-unstable bundle $E^{cu}$ is \emph{volume expanding} if
there exists $K,\theta>0$ such that
$|\det(DX_t| E^{cu}_x)|\geq K e^{\theta t}$ for all $x\in \Lambda$,
$t\geq 0$.

\begin{definition} \label{def:singhypatt} Let $\Lambda$ be a compact
  invariant set for $X \in {\fX}^r(M)$.  We say that $\Lambda$ is a
  \emph{singular-hyperbolic set} if all equilibria in $\Lambda$ are
  hyperbolic, and $\Lambda$ is partially hyperbolic with volume
  expanding two-dimensional center-unstable bundle ($d_{cu}=2$).  A
  singular-hyperbolic set which is also an attracting set is called a
  {\em singular-hyperbolic attracting set}.
\end{definition}

\begin{remark} \label{rmk:per} A singular-hyperbolic attracting set
  contains no isolated periodic orbits.  For such a periodic orbit
  would have to be a periodic sink, violating volume expansion.
\end{remark}

\begin{theorem}{\cite[Lemma 3]{MPP99}}
  \label{thm:hyplemma}
  Every compact invariant set without singularities of a
  singular-hyperbolic set is hyperbolic.
\end{theorem}

A subset $\Lambda \subset M$ is \emph{transitive} if it has
a full dense orbit, that is, there exists $x\in \Lambda$ such that
$\overline{\{X_tx:t\ge0\}}=\Lambda= \overline{\{X_tx:t\le0\}}$.

\begin{definition}\label{def:attractor}
  A \emph{singular-hyperbolic attractor} is a transitive
  singular-hyperbolic attracting set.
\end{definition}

\begin{proposition}{\cite[Proposition
    2.6]{ArMel18}} \label{pr:Lorenz-like} Suppose that $\Lambda$ is a
  singular-hyperbolic attractor and let $\sigma\in\Lambda$ be an
  equilibrium.  Then $\sigma$ is \emph{Lorenz-like}.  That is,
  $DG(\sigma)|E^{cu}_\sigma$ has real eigenvalues $\lambda^s$,
  $\lambda^u$ satisfying $-\lambda^u<\lambda^s<0<\lambda^u$.
\end{proposition}

\begin{remark}
  \label{rmk:no-recur-non-Lorenz}
  Some consequences of singular-hyperbolicity follow.
  \begin{enumerate}
  \item Partial hyperbolicity of $\Lambda$ implies that the direction
    $X(x)$ of the flow is contained in the center-unstable bundle
    $E^{cu}_x$ at every point $x$ of $\Lambda$ (see \cite[Lemma
    5.1]{ArArbSal}).

  \item The index of a singularity $\sigma$ in a singular-hyperbolic
    set $\Lambda$ equals either $\dim E^s$ or $1+\dim E^s$.  That is,
    $\sigma$ is either a hyperbolic saddle with $\dim M - \dim E^s_\sigma = 2$ (that
    is, the \emph{codimension} of $E^s_\sigma$ equals $2$) or a
    Lorenz-like singularity.

  \item If a singularity $\sigma$ in a singular-hyperbolic set
      $\Lambda$ is not Lorenz-like, then there is no regular orbit of
      $\Lambda$ that accumulates $\sigma$ in the positive time
      direction. In other words, there is no $x \in \Lambda$ regular such that $\sigma
      \in \omega(x)$ (see \cite[Remark 1.5]{ArSzTr})
  \end{enumerate}
\end{remark}

\begin{definition}\label{def:nontrivialsinghyp}
  A singular-hyperbolic invariant set is \emph{nontrivial} if it is a
  non-trivial compact invariant subset which contains some Lorenz-like
  equilibrium.
\end{definition}

\subsubsection{Physical measures}
\label{sec:physical-measures}

The existence of a unique invariant and ergodic physical measure for
singular-hyperbolic attractors was first proved for $3$-dimensional
manifolds in \cite{APPV} and extended to singular-hyperbolic
attracting sets in e.g.~\cite{ArSzTr}. For sectional-hyperbolic
attractors\footnote{That is the same as singular-hyperbolicity, but
  allowing $\dim E^{cu}>2$ and demanding that volume expansion holds
  along every two-dimensional subspace of $E^{cu}$.}, existence and
  uniqueness of physical measure was obtained in~\cite{LeplYa17} and
  recently extended to attracting sets in~\cite{Araujo19}. In fact,
  sectional-hyperbolic attracting sets have finitely many ergodic
  physical measures which are equilibrium states for the
  central-unstable Jacobian, just like Axiom A attracting sets.

\begin{theorem}
  {\cite[Theorem 1.7]{ArSzTr}}
  \label{theorem:existence-physical-measures}
  Let $\Lambda$ be a singular-hyperbolic attracting set for
  a $C^2$ vector field $X$ with the open subset $U$ as
  trapping region. Then
  \begin{enumerate}
  \item there are finitely many ergodic physical/SRB
    measures $\mu_1,\dots,\mu_k$ supported in $\Lambda$ such
    that the union of their ergodic basins covers $U$
    Lebesgue almost everywhere:
    $
      \leb\left(U\setminus\big(\bigcup_{i=1}^k
      B(\mu_i)\big)\right)=0.
    $
  \item Moreover, for each $X$-invariant ergodic probability
    measure $\mu$ supported in $\Lambda$ the following are
    equivalent
    \begin{enumerate}
    \item
      $h_\mu(X_1)=\int\log|\det DX_1\mid_{E^{cu}}|\,d\mu>0$;
    \item $\mu$ is a $SRB$ measure, that is, admits an
      absolutely continuous disintegration along unstable
      manifolds;
    \item $\mu$ is a physical measure, i.e., its basin
      $B(\mu)$ has positive Lebesgue measure.
    \end{enumerate}
  \item The family $\EE$ of all $X$-invariant
    probability measures which satisfy item (2a) above is the
    convex hull
    $\EE=\left\{\sum_{i=1}^k t_i \mu_i : \sum_i t_i=1; 0\le
    t_i\le1, i=1,\dots,k\right\}. $
  \end{enumerate}
\end{theorem}

We note that there are many examples of singular-hyperbolic attracting
sets, non-transitive and containing non-Lorenz-like singularities; see
Subsection~\ref{eq:lorenz-equations}.

\subsection{Statement of results}
\label{sec:statement-results}

We can now state our main results.  In what follows, we write
$C^{k+\eta}(M)$, where $\eta\in(0,1]$ is a real number and $k\ge0$ is
a non-negative integer, for the set of functions $\vfi:M\to\RR$ which
are of class $C^k$ and the $k$th derivative $D^k\vfi$ is
$\eta$-H\"older. This is a Banach space with norm given by
\begin{align*}
  \|\vfi\|_{k+\eta}:=\sum_{i=0}^k|D^i\vfi|_\infty+|D^k\vfi|_\eta,
\end{align*}
where for any function $\psi:M\to\RR$ we set
$|\psi|_\infty := \sup_{x\in M}|\psi(x)|$ and
$|\psi|_\eta := \sup_{x \neq y}|\psi(x) - \psi(y)| / \dist(x,y)^\eta$.

\begin{maintheorem}[Exponential mixing]
  \label{mthm:decay-of-correlations-singular-hyperbolic}
  There exists an open subset $\cU \subset \fX^2(M)$ such that each
  vector field $X\in\cU$ admits a non-trivial connected singular-hyperbolic attracting set $\Lambda$ such that, given $X \in \cU$ and
  $\mu$ a physical measure supported in $\Lambda$, there exist
  constants $c, C > 0$ such that for any $\eta\in(0,1]$ we have
  $|C_t(\vfi,\psi)|\leq C e^{-c t}\|\varphi\|_{\eta}\|\psi\|_{\eta},$
  for all $\varphi,\psi \in C^\eta(M)$ and $t > 0$.
\end{maintheorem}

We can also present this result with a different appearence.  If
$\mu_1,\dots,\mu_k$ are the ergodic invariant physical probability
measures of $X$ supported in $\Lambda$ as given by
Theorem~\ref{theorem:existence-physical-measures} and
$\vartheta_i=\Leb(B(\mu_i))$ is the volume of each of their basins,
then the normalized Lebesgue measure, on a trapping region $U$ for
$\Lambda$, can be written as a linear convex combination
$\Leb=\sum_i\vartheta_i\Leb_i$ where
$\Leb_i=\vartheta_i^{-1}\Leb\mid_{B(\mu_i)}$.

\begin{maincorollary}[Exponential convergence to equilibrium]
  \label{mcor:convergence-equilibrium}
  In the same setting of
  Theorem~\ref{mthm:decay-of-correlations-singular-hyperbolic},
  for each $0<\eta\le 1$ 
  there exist constants $c, C > 0$ such that
  \begin{align*}
  \left| \int (\varphi \circ X_t) \psi\, d\Leb - \int \varphi\,
    d\widetilde{\mu}\int \psi\, d\Leb\right| \leq Ce^{-c
    t}\|\varphi\|_\eta\|\psi\|_\eta,
  \end{align*}
  for all $\vfi,\psi \in C^\eta(M)$ and $t > 0$, where
  $\widetilde\mu=\sum_i\vartheta_i\mu_i$.
\end{maincorollary}

If we are dealing with an attractor, that is, if $\Lambda$ is
transitive, then there is a unique physical measure and putting
$\psi \equiv 1$ in the statement of
Corollary~\ref{mcor:convergence-equilibrium} we get
$ \big| \int \varphi \circ X_t\, d\Leb - \int \varphi\, d\mu \big|
\leq Ce^{-c t} |\varphi|_{C^1} $ for all $\varphi \in C^1(M)$ and all
$t > 0$, that is, $(X_t)_*\Leb$ converges exponentially fast to the
physical (also known as ``natural'') measure $\mu$ in the weak-*
topology when $t$ goes to infinity.


\subsection{Some consequences of fast mixing}
\label{sec:some-conseq-fast}

It is known that fast decay of correlations for a dynamical system
implies many other statistical properties.

The base map of a hyperbolic skew-product semiflow is known to satisfy
exponential mixing for H\"older observables with respect to its
physical measures; see e.g. \cite{ArGalPac}. This in turn
automatically implies certain statistical properties for the induced
measure on the suspension flow: the Central Limit Theorem, the Law of
the Iterated Logarithm and the Almost Sure Invariance Principle; see
e.g. \cite{MeTo04}.

For mixing and the speed of mixing, the properties of the base map do
not extend to the suspension flow in general: the suspension flow does
not even have to be mixing. More precisely, see \cite{Pol85}, rates of
mixing of suspension flows can be arbitrarily slow even if the base
map is exponentially mixing.

Having a flow which mixes exponentially fast should imply more subtle
statistical properties. In fact, some statistical properties of the
time-$1$ map of singular-hyperbolic flows near attracting sets can be
obtained in this way. 

\begin{maincorollary}[Consequences of exponential mixing for the
  time-$1$ map]
  Let $\cU$ be as in Theorem
  \ref{theorem:conjugacy-skew-product-hyperbolic}.  Given $X \in \cU$,
  let $\mu$ be an ergodic physical measure for the singular-hyperbolic
  attracting set $\Lambda$. For all $\varphi \in C^1(M)$ it holds:
\begin{enumerate}[(1)]
\item \emph{(Central Limit Theorem (CLT) for the time-$1$ map).}  There
  exists $\sigma \geq 0$ such that we have the following convergence
  in distribution
  \begin{align*}
    \dfrac{1}{\sqrt{n}}\left[ \sum^{n - 1}_{j = 0} \varphi \circ X_j - n \int \varphi\,
    d\mu \right] \xrightarrow[n \to +\infty]{\cD} \cN(0, \sigma^2)
  \end{align*}
  Moreover, if $\sigma = 0$, then for every periodic point
  $q \in \Lambda$, there exists $T > 0$ (independent of $\varphi$)
  such that $\int^T_0 \varphi(X_t(q))\, dt = 0$.
\item \emph{(Almost Sure Invariant Principle (ASIP) for the time-$1$ map).}
  Passing to an enriched probability space, there exists a sequence
  $Y_0, Y_1, \dots$ of iid normal random variables with mean zero and
  variance $\sigma^2$ such that
  $$\sum^{n - 1}_{j = 0}\varphi \circ X_j = n \int \varphi\, d\mu + \sum^{n - 1}_{j =
    0} Y_j + \cO(\sqrt[4]{n \log \log n}\sqrt{\log n}),\ a.e.$$
\end{enumerate}
\end{maincorollary}

This corollary follows from the proof of exponential mixing just as in
\cite{AMV15}, where the same was deduced from superpolynomial decay of
correlations.

The ASIP implies the CLT and also the functional CLT (weak invariance
principle), and the law of the iterated logarithm together with its
functional version, as well as numerous other results.  The reader
should consult~\cite{PhilippStout75} for a comprehensive list.

\subsubsection{Organization of the text}
\label{sec:organization-text}

The remainder of the paper is organized as follows.  

In Section~\ref{sec:strategy-comments-co} we present the overall
organization of the proof, open classes of examples in the setting of
our main results and some conjectures to extend the results presented
in the text.

In Section~\ref{chapter:decay-correlations-suspension}, we present
general properties of partially hyperbolic attracting sets and
singular-hyperbolic attracting sets which will enable us to find a
global Poincar\'e section for the flow in a neighborhood of the
attracting set. The corresponding global Poincar\'e return map is
piecewise hyperbolic in a precise sense.

In Section~\ref{sec:propert-one-dimens} we describe crucial properties
of the one-dimensional quotient map of the global Poincar\'e map over
the leaves of the stable foliation, and associate to each ergodic
physical measure of the flow an hyperbolic skew-product semiflow.

In Section~\ref{cha:exponent-mixing-sing} we prove our Main
Theorem~\ref{mthm:decay-of-correlations-singular-hyperbolic} and
Corollary~\ref{mcor:convergence-equilibrium} using all the previous
results.

Finally, in
Section~\ref{chapter:decay-correlations-suspension-general} we prove a
technical result, Theorem~\ref{theorem:decay-onedimensional-weaker},
which is crucial to the previous arguments.

\subsection*{Acknowledgements}

This is based on the PhD thesis of E. Trindade at the Instituto de
Matematica e Estatistica-Universidade Federal da Bahia (UFBA) under a
CAPES scholarship. E.T.  thanks the Mathematics and Statistics
Institute at UFBA for the use of its facilities and the financial
support from CAPES during his M.Sc. and Ph.D. studies.  We thank
A. Castro; Y. Lima; D. Smania and P. Varandas for many comments and
suggestions which greatly improved the text. We also thank the
anonymous referee for the useful suggestions that improved the text.

\section{Strategy of the proof}
\label{sec:strategy-comments-co}

As in previous works on robust exponential
mixing for geometric Lorenz attractors \cite{ArVar,AMV15,ArMel16},
the proof relies on finding a convenient conjugation between the flow in a
neighborhood of the attracting set $\Lambda$ and a skew-product
semiflow satisfying strong dynamical and ergodic properties.

We present this semiflow in what follows and then state the main
technical result which is behind
Theorem~\ref{mthm:decay-of-correlations-singular-hyperbolic} and
Corollary~\ref{mcor:convergence-equilibrium}.

\subsection{Hyperbolic skew product semiflow}
\label{section:hyperbolic-skew-product-semiflow}

The main strategy of this work is to take a flow admitting a
singular-hyperbolic set (with some assumptions that will be presented
along the text) and reduce it to the setting that we present in this
section. After obtaining the results for hyperbolic skew product
semiflows, we explain how to take them to the original flow.

\subsubsection{Uniformly expanding maps}
\label{subsection:uniformly-expanding-maps}

Let $\alpha \in (0,1]$ and $\Delta$ be a compact interval of
$\RR$. Without loss of generality we assume that $\Delta = [0,1]$ in
this section. Let $\cP = \{(c_m,d_m):\ m \in \NN\}$ be a countable
partition ($\leb\bmod0$) of $\Delta$. Let $F: \Delta \to \Delta$ be
$C^{1 + \alpha}$ on each element $J$ of the partition $\cP$ with
$F(J) = \Delta$ and $F$ extends to a homeomorphism from $\overline{J}$
to $\Delta$, for every $J \in \cP$. Given $J \in \cP$, we say that a
map $h: \Delta \to \overline{J}$ is an \emph{inverse branch} of $F$ if
$F \circ h = id$. We denote by $\cH$ and $\cH_n$ the set of all
inverse branches of $F$ and $F^n$, respectively, for all $n \geq 1$.

Given a function $\psi: \Delta \to \RR$ we denote
$|\psi|_\infty := \sup_{x\in\Delta}|\psi(x)|$ and
$|\psi|_\alpha := \sup_{x \neq y}|\psi(x) - \psi(y)| / |x -
y|^\alpha$.

We say that $F$ is a $C^{1 + \alpha}$ \emph{uniformly expanding map}
if there exist constants $C > 0$ and $\rho \in (0,1)$ such that

\begin{enumerate}
\item $|h'|_\infty \leq C\rho^n$ for all $h \in \cH_n$, \label{item:backward-contraction}
\item $|\log |h'||_\alpha \leq C$ for all $h \in \cH$.
\end{enumerate}

\begin{remark}\label{rmk:sum-derivative-inverse-branch}
  It follows from (1) and (2) that 
$
    \sum_{h \in \cH_n}|h'|_\infty < \infty.
$ 
\end{remark}

It is standard that $C^{1 + \alpha}$ uniformly expanding maps have a
unique absolutely continuous $F$-invariant ergodic measure with
$\alpha$-H\"older positive density function bounded from above and
bellow away from zero. We denote this measure by $\mu_F$.

\subsubsection{$C^{1 + }$ expanding semiflows}
\label{subsec:expanding-semiflows}

Consider a function $r: \Delta \to (0, +\infty)$ which is $C^1$ on each
element of the partition $\cP$. We assume the following conditions on $r:$

\begin{enumerate}
  \setcounter{enumi}{2}
\item $|(r \circ h)'|_\infty \leq C$ for all $h \in \cH$; \label{item:derivative-r-h}
\item $r$ \emph{has exponential tail}: there exists $\varepsilon > 0$ such that $\sum_{h \in
    \cH}e^{\varepsilon |r \circ h|_\infty}|h'|_\infty < \infty$; \label{item:expontial-tail}
\item \emph{uniform non-integrability (UNI)}: it is not possible to
  write $r= \psi + \varphi \circ F - \varphi$ with
  $\psi: \Delta \to \RR$ constant in elements of the partition $\cP$
  and $\varphi: \Delta \to \RR$ a $C^1$ function.
\end{enumerate}

Let
$\Delta^r = \{(x,u) \in \Delta \times \RR: \ 0 \leq u \leq r(x)\} /
\sim$ be a quotient space, where $(x, r(x)) \sim (F(x), 0)$, and
define the \emph{suspension semiflow} $F_t: \Delta^r \to \Delta^r$
with \emph{roof function $r$} by $F_t(x,u) = (x, u + t)$, for all
$t \geq 0$, computed modulo the given identification. The semiflow
$F_t$ has an ergodic invariant probability measure
$\mu^r_F = (\mu_F \times \Leb) / \int_\Delta r\, d\mu_F$. If
conditions (1)-(4) hold, then we say that $F_t$ is a $C^{1 + \alpha}$
\emph{expanding semiflow}.

\subsubsection{Decay of correlations for $C^{1 + }$ expanding semiflows}
\label{subsec:decay-expanding-semiflows}

We define $C^\alpha_\loc(\Delta^r)$ to consist of $L^\infty$ functions
$\psi: \Delta^r \to \RR$ such that
$\|\psi\|_\alpha = |\psi|_\infty + |\psi|_{\alpha, \loc} < \infty$,
where
\begin{align*}
  |\psi|_{\alpha, \loc} =
  \sup_{h \in \cH}\sup_{(x_1, u) \neq (x_2, u)}\dfrac{|\psi(hx_1,u) -
  \psi(hx_2,u)|}{|x_1 - x_2|^\alpha}.
\end{align*}
Given an integer $k \geq 1$, define
$C^{\alpha,k}_\loc(\Delta^r)$ to consist of $C^\alpha_\loc(\Delta^r)$ functions
$\psi$ with $\|\psi\|_{\alpha, k} = \sum^k_{j = 0}\|\partial^j_t
\psi\|_\alpha < \infty$, where $\partial_t$ denotes the differentiation along
the semiflow direction.

\begin{theorem}[Decay of correlations for expanding
  semiflows]\label{theorem:decay-onedimensional-weaker} If conditions (1)-(5)
  hold, then there are constants $c, C > 0$ so that for all
  $\varphi \in L^\infty(\Delta^r),$
  $\psi \in C^{\alpha, 2}_\loc(\Delta^r)$,
  \begin{align*}
\left| \int (\varphi \circ F_t) \psi\, d\mu^r_F - \int \varphi\,
    d\mu^r_F \int \psi\, d\mu^r_F \right| \leq
  Ce^{-ct}|\varphi|_\infty\|\psi\|_{\alpha, 2}, \quad \forall t > 0.
  \end{align*}
\end{theorem}

Theorem \ref{theorem:decay-onedimensional-weaker} is a generalization
of \cite[Theorem 2.1]{ArMel16}. The original result was proved for
$\alpha$-Hölder observables. We extend to the more general class of
observables presented above; see
Section~\ref{chapter:decay-correlations-suspension-general} for a
proof. This generality is needed to transfer the results obtained for
semiflows to the original singular-hyperbolic flow, as will become
clear in Section~\ref{cha:exponent-mixing-sing}. This is analogous to
the introduction of ``dynamical observables'' in a similar setting to
study rapid mixing; see \cite{Melbourne_2018}.

\subsubsection{Hyperbolic skew products}
\label{subsection:hyperbolic-skew-products}

Let $F: \Delta \to \Delta$ be a $C^{1 + \alpha}$ expanding map, as in
Subsection~\ref{subsection:uniformly-expanding-maps}, and $\Omega$ a
compact Riemannian manifold inside $\RR^N$, for some integer
$N \geq 1$. Let $\widehat{\Delta} = \Delta \times \Omega$ be a direct
product endowed with the distance given by
$|(x_1,y_1) - (x_2, y_2)| = |x_1 - x_2| + |y_1 - y_2|$.  Consider also
$G: \widehat{\Delta} \to \Omega$ a $C^{1 + \alpha}$ map and define
$\widehat{F}: \widehat{\Delta} \to \widehat{\Delta}$ by
$\widehat{F}(x,y) = (F(x), G(x,y))$. We say that $\widehat{F}$ is a
\emph{uniformly hyperbolic skew product} if it satisfies

\begin{enumerate}
  \setcounter{enumi}{5}
\item (uniform contraction along $\Omega$) there exist constants
  $C > 0$ and $\gamma \in (0,1)$ such that
  $|\widehat{F}^n(x,y_1) - \widehat{F}^n(x, y_2)| \leq C\gamma^n|y_1 -
  y_2|,$ for all $x \in \Delta$ and $y_1, y_2 \in \Omega$.
\end{enumerate} For each integer $n \geq 1$, we denote the iterates of
$\widehat{F}$ by $\widehat{F}^n(x,y) = (F^n(x), G_n(x,y))$ for all
$(x,y) \in \widehat{\Delta}$. Hence, item (6) above becomes 
$|G_n(x,y_1) - G_n(x,y_2)| \leq C \gamma^n|y_1 - y_2|,$ for all $(x,
y_i) \in \widehat{\Delta}$, $i = 1, 2$.

Let $\pi: \widehat{\Delta} \to \Delta$ be the projection
$\pi(x,y) = x$, for all $(x,y) \in \widehat{\Delta}$. Note that
$\pi \circ \widehat{F} = F \circ \pi$, that is, $\pi$ is a
\emph{semiconjugacy} between $\widehat{F}$ and $F$. Moreover, the
property (4) says that the leaf $\pi^{-1}(x)$ is exponentially
contracted by the skew product $\widehat{F}$, for all $x \in \Delta$.

\subsubsection*{Invariant probability measure for the skew
  product}\label{sec:invari-probab-measur}
  
In the following proposition we recall how to obtain a
$\widehat{F}$-invariant probability measure using the (absolutely
continuous) invariant probability measure $\mu_F$ for the map $F$.

\begin{proposition}{\cite[Section 6]{APPV}}
  \label{prop:invariant-measure-for-skew-products}
  Let $\varphi: \widehat{\Delta} \to \RR$ be a continuous function and
  define $\varphi_\pm: \Delta \to \RR$ by
  $\varphi_+(x) = \sup_{y \in \Omega}\varphi(x,y)$ and
  $\varphi_-(x) = \inf_{y \in \Omega}\varphi(x,y)$. Then the limits
  $\lim_{n \to +\infty}\int_{\Delta}(\varphi \circ \widehat{F}^n)_+\,
  d\mu_F$ and
  $\lim_{n \to +\infty}\int_{\Delta}(\varphi \circ \widehat{F}^n)_-\,
  d\mu_F$ exist, are equal, and define a $\widehat{F}$-invariant
  probability measure $\mu_{\widehat{F}}$ such that
  $\pi_*\mu_{\widehat{F}} = \mu_F$.
\end{proposition}

\subsubsection{Hyperbolic skew product semiflow}
\label{sec:hyperb-skew-product}

Let $F:\Delta \to \Delta$ be a $C^{1 + \alpha}$ uniformly expanding
map with partition $\cP$;
$\widehat{F}: \widehat{\Delta} \to \widehat{\Delta}$ a
$C^{1 + \alpha}$ hyperbolic skew product with
$\pi \circ \widehat{F} = F \circ \pi$ as in the previous
Subsections~\ref{subsec:expanding-semiflows}
and~\ref{subsection:hyperbolic-skew-products}; and
$r: \Delta \to (0, +\infty)$ be $C^1$ on elements of the partition
$\cP$ with $\inf r > 0$. \emph{We extend the definition of $r$ to
  $\widehat{\Delta}$ by setting}\footnote{Note that here we are
  assuming that the return time to the base of the semiflow is
  constant on stable leaves.} $r(x,y) = r(x)$ for all
$(x,y) \in \widehat{\Delta}$. Considering the quotient space
$\widehat{\Delta}^r = \{(z,u) \in \widehat{\Delta} \times \RR:\ 0 \leq
u \leq r(z)\}/\sim$, where $(z,r(z)) \sim (\widehat{F}(z), 0)$, we
define the \emph{suspension semiflow} $\widehat{F}_t$ with \emph{roof
  function} $r$ by $\widehat{F}_t(z,u) = (z, u + t)$, for all
$t \geq 0$, computed modulo the given identification. This semiflow
has an ergodic invariant probability measure
$\mu^r_{\widehat{F}} = \mu_{\widehat{F}} \times \Leb /
\int_{\widehat{\Delta}} r\, d\mu_{\widehat{F}}$. If $r$ satisfies the
conditions \eqref{item:derivative-r-h} and
\eqref{item:expontial-tail}, then we say that the $\widehat{F}_t$ is a
\emph{$C^{1 + \alpha}$ hyperbolic skew product semiflow}.

\subsubsection*{Exponential mixing for hyperbolic skew product semiflows}
\label{subsec:decay-hyperbolic-skew-product-flow-statement}  

Let $C^\alpha_\loc(\widehat{\Delta}^r)$ denote the subset of
$L^\infty$ functions $\psi: \widehat{\Delta}^r \to \RR$ such that
$\|\psi\|_\alpha = |\psi|_\infty + |\psi|_{\alpha, \loc}$,
where
$$|\psi|_{\alpha, \loc} = \sup_{h \in \cH}\sup_{(x_1,y_1,u) \neq
  (x_2,y_2,u)}\frac{|\psi(hx_1,y_1,u) - \psi(hx_2,y_2,u)|}{|x_1 -
  x_2|^\alpha + |y_1 - y_2|}$$ and let
$C^{\alpha,k}_\loc(\widehat{\Delta}^r)$ be the subset of
$C^\alpha_\loc(\widehat{\Delta}^r)$ functions
$\varphi:\widehat{\Delta}^r \to \RR$ such that
$\|\psi\|_{\alpha,k} := \sum^k_{j = 0}|\partial^j_t \psi|_\alpha <
\infty$, where $\partial_t$ denotes the differentiation along the
semiflow direction and $k\ge1$ is a given integer.

\begin{theorem}\label{prop:decay-correlations-skew-product-weaker}
  Suppose that
  $\widehat{F}_t:\widehat{\Delta}^r \to \widehat{\Delta}^r$ is a
  $C^{1 + \alpha}$ hyperbolic skew product with roof function $r$
  satisfying the UNI condition (5). Then there exist constants $c, C > 0$
  such that
  $\big|\int(\varphi \circ \widehat{F}_t)\cdot \psi\,
    d\mu^r_{\widehat{F}} - \int \varphi\, d\mu^r_{\widehat{F}}\int
    \psi\, d\mu^r_{\widehat{F}}\big| \leq Ce^{-ct}\|\varphi\|_\alpha
  \|\psi\|_{\alpha,2},$ for all
  $\varphi \in C^\alpha_\loc(\widehat{\Delta}^r)$,
  $\psi \in C^{\alpha, 2}_\loc (\widehat{\Delta}^r)$ and $t > 0$.
\end{theorem}

Theorem \ref{prop:decay-correlations-skew-product-weaker} is a
generalization of \cite[Theorem 3.3]{ArMel16}. As already noted (after
the statement of Theorem~\ref{theorem:decay-onedimensional-weaker}),
here we also need to relax the conditions on the observables
(obtaining ``dynamical observables'') of the original theorem to fit
our needs. The proof of this theorem can be found in Section
\ref{section:decay-hyperbolic-skew-product-semiflows}.

\subsection{The main technical result}
\label{sec:main-technic-result}

We present now our main technical result at the core of
Theorem~\ref{mthm:decay-of-correlations-singular-hyperbolic} and
Corollary~\ref{mcor:convergence-equilibrium}.  We construct a $C^2$
open set of vector fields that are semiconjugated to a
$C^{1 + \alpha}$ hyperbolic skew product semiflow and have the
necessary properties that allow us to transfer the decay of
correlations obtained in Theorem
\ref{prop:decay-correlations-skew-product-weaker} to the original
flow.

\begin{theorem}\label{theorem:conjugacy-skew-product-hyperbolic}
  There exists an open subset $\cU \subset \fX^2(M)$ such that each
  vector field $X\in\cU$ admits a non-trivial connected
  singular-hyperbolic attracting set with $U$ as trapping region and
  $\alpha\in(0,1)$ so that, for all small enough $\epsilon>0$ the
  following holds.  We can
  find a $C^\infty$ function $\varrho: M\to(1/2,3/2)$, which is
  $\epsilon$-$C^2$-close to $1$, and such that $Y=\varrho\cdot X$
  admits a $C^2$-neighborhood $\cV\subset\cU$ satisfying: for each
  ergodic physical measure $\mu$ of $Z\in\cV$ supported in $U$, there
  exists a $C^{1+\alpha}$ hyperbolic skew product semiflow
  $\widehat{F}_t: \widehat{\Delta}^r \to \widehat{\Delta}^r$ with roof
  function $r$ satisfying the UNI condition and a map
  $p: \widehat{\Delta}^r \to U$ satisfying:
  \begin{enumerate}[(i)]
  \item $Z_t\circ p = p \circ \widehat{F}_t$, for all $t > 0$ and
    $p_*\mu^r_{\widehat{F}} = \mu$;
    \item there exists a constant $C > 0$ such that
      $\|\varphi \circ p\|_\alpha \leq C|\varphi|_{C^1}$ for all
      $\varphi \in C^1(U)$ and
      $\|\psi \circ p\|_{\alpha,2} \leq C|\psi|_{C^3}$ for all
      $\psi \in C^3(U)$. 
  \end{enumerate}
\end{theorem}
Here and in what follows we write $|\cdot|_{C^k}$ for the $C^k$-norm
$\|\cdot\|_k$ of real functions on a manifold.  The proof of this
result is the content of the following sections.

\begin{remark}
  \label{rmk:topeq}
  Theorem~\ref{theorem:conjugacy-skew-product-hyperbolic} can be
  interpreted as: \emph{every singular-hyperbolic attracting set is
    robustly exponentially mixing} with respect to its physical
  measures \emph{modulo an arbitrary small perturbation of the speed
    of the vector field}.
\end{remark}

\subsection{$q$-dissipativity}
\label{sec:q-dissipativity}

We recall the following consequence of the Whitney Embedding and
Tubular Neighborhood Theorems: if $\Lambda$ is an attracting set of a
vector field $X$ of a compact finite-dimensional manifold $M$ then,
after embedding the manifold into some Euclidean space $\RR^N$, we may
extend $X$ to a neighborhood of $M$, so that $M$ and $\Lambda$ become
attracting sets of the extended vector field, with the same
smoothness. Hence we assume without loss of generality in what follows
that $X$ is a smooth vector field on a compact region of some
Euclidean space.

\subsubsection{$q$-dissipativity and smooth stable foliation}
\label{sub:dissipativity}

Let $\mathcal{M}$ denote the set of $X_t$-invariant ergodic
probability measures on $\Lambda\subset M\subset\RR^N$.  If
$A = [a_{ij}]$ is an $N \times N$ real matrix, we denote
$\|A\|_2 = \left(\sum^N_{i = 1}\sum^N_{j = 1}a^2_{ij}\right)^{1/2}$.
For each $m\in\cM$, we label the Lyapunov exponents
\[ \chi_1(m)\le \chi_2(m)\le\dots\le\chi_d(m),\] it follows that
$ \Theta \leq \chi_1(m)\le \chi_2(m)\le\dots\le\chi_{d_s}(m) \leq
-\lambda < 0,$ where $\lambda$ comes from the definition of partially
hyperbolic set (see Subsection~\ref{sec:PH}) and $\Theta$ is given by
$\log\inf_{x \in U} \|(DX(x))^{-1}\|^{-1}_2$. Because $\lambda$ and
$\Theta$ are independent of the measure $m \in \mathcal{M}$, it is
possible to choose $\ell > 1$ such that
\begin{align}
  \label{eq:ell-dissipative}
  \ell \geq
  \frac{1}{d_s\chi_{d_s}(m)}\sum^{d_s}_{j=1}\chi_j(m)
  \quad \text{for all}\quad m \in \mathcal{M}.
\end{align}

\begin{definition}\label{def:SD}
  Let $X \in \fX^1(M)$ admitting a partially hyperbolic attracting
  $\Lambda \subset M$ be given and let $\ell, q > 1$ such that
  $\ell q > 1/d_s$ and $\ell$ satisfy~\eqref{eq:ell-dissipative}. We
  say that $\Lambda$ is \emph{$q$-strongly dissipative}\footnote{This
    definition was first given in \cite{ArMel17} but its statement was
    only valid for $3$-flows. We present here a corrected proof for
    completeness.} if
  \begin{enumerate}[(a)]
  \item for every equilibrium $\sigma \in \Lambda$ (if any), the
    eigenvalues $\lambda_j$ of $DX(\sigma)$, ordered so that
    $\mathfrak{R}\lambda_1 \leq \mathfrak{R}\lambda_2 \leq \cdots \leq
    \mathfrak{R}\lambda_N$, satisfy $\mathfrak{R}(\lambda_{d_s} -
     \lambda_{d_s + 1} + q\lambda_N) < 0$;
  \item $\sup_{x \in \Lambda}\{\Div X(x) + (\ell d_sq - 1)\|DX(x)\|_2 \}
    < 0$
  \end{enumerate}
\end{definition}
The stable foliation of a singular-hyperbolic attracting set $\Lambda$
is $C^q$ on a neighborhood of $\Lambda$ if this set is $q$-strongly
dissipative for a $C^q$ vector field. As we will see in the next
subsection, this allows us to consider smooth cross-sections of $X$ to
be composed by stable discs $W^s_x$.

\begin{theorem}\label{thm:smooth-foliation}
  Let $\Lambda$ be a sectional-hyperbolic attracting set with respect
  to $X\in\fX^2(M)$ with a trapping
  region $U_0$. Suppose that $\Lambda$ is $q$-strongly dissipative for
  some $q \in (1/d_s, 2]$.  Then there exists a neighborhood $U_0$ of
  $\Lambda$ such that the stable manifolds $\{W^s_x:\ x \in U_0\}$
  define a $C^q$ foliation of $U_0$.
\end{theorem}

\begin{proof}
  For each $t \in \RR$, let
  $\eta_t(x)=\log\big\{\|DX_t|E^s_x\|\cdot
  \|DX_{-t}|E^{cu}_{X_tx}\|\cdot \|DX_t|E^{cu}_x\|^q\big\}$. Note that
  $\{\eta_t:\ t\in\RR\}$ is a continuous family of continuous
  functions each of which is subadditive, that is,
  $\eta_{s+t}(x)\le \eta_s(x)+\eta_t(X_sx), s,t\ge0, x\in\Lambda$.

We claim that for each $m\in\cM$, the limit
$\lim_{t\to\infty}t^{-1} \eta_t(x)$ exists and is negative for
$m$-almost every $x\in\Lambda$. It then follows
from~\cite[Proposition~3.4]{arbieto2010} that there exists constants
$C,\beta>0$ such that $\exp\eta_t(x)\le Ce^{-\beta t}$ for all $t>0$,
$x\in\Lambda$. In particular, for $t$ sufficiently large,
$\exp\eta_t(x)<1$ for all $x\in\Lambda$. Hence, for such $t$, we
obtain
$ \|DX_t|_{E^s_x}\|\cdot \|DX_{-t}|_{E^{cu}_{X_tx}}\| \cdot\|DX_t
|_{E^{cu}_x}\|^q <1$ for all $x\in \Lambda$.  From this last
inequality the result follows from \cite[Theorem 4.12]{ArMel17} and
\cite[Remark 4.13]{ArMel17}.

It remains to verify the claim. Since $\Lambda$ is partially
hyperbolic, the Lyapunov exponents $\chi_j(m)$, $j=1,\dots,d_s$ are
associated with $E^s$ and are negative, while the remaining exponents
are associated with $E^{cu}$.

We have $\lim\log \|DX_t|E^s_x\|^{1/t}=\chi_{d_s}(m)$ and
$\lim\log \|DX_{-t}|E^{cu}_{X_tx}\|^{1/t}=-\chi_{d_s+1}(m)$, for
$m$-a.e. $x\in\Lambda$ as $t\nearrow\infty$, and also
  \begin{align*}
\lim_{t\to\infty}\log \|DX_t|E^{cu}_x\|^{1/t}= \lim_{t\to\infty} \log
\|DX_t\mid T_xM\|^{1/t}= \chi_d(m).
  \end{align*}
  Hence
$
\lim_{t\to\infty}t^{-1}\eta_t(x)=\chi_{d_s}(m) -
\chi_{d_s+1}(m)+q\chi_d(m),
m$-almost everywhere,

  If $m$ is a Dirac delta at an equilibrium $\sigma\in \Lambda$, then
$\chi_j(m)=\Re\lambda_j$ for $j=1,\dots,d$, where $\lambda_j$ are the
eigenvalues of $DG(\sigma)$. Hence, it is immediate from
Definition~\ref{def:SD}(a) that $\lim_{t\to\infty}t^{-1}\eta_t(\sigma)<0$.

  If $m$ is not supported on an equilibrium, then there is a zero Lyapunov
exponent in the flow direction. Sectional expansion ensures that
$\chi_{d_s+1}(m)=0$ and $\chi_j(m)>0$ for $j=d_s+2,\dots,d$. Hence using
inequality~\eqref{eq:ell-dissipative}, $m$-almost
everywhere, 
\begin{align*}
  \lim_{t\to\infty}t^{-1}\eta_t(x)
  &\textstyle= \chi_{d_s}(m)+q\chi_d(m) \le (\ell
    d_s)^{-1}\sum_{j=1}^{d_s}\chi_j(m)+q\chi_d(m)
  \\
  &\textstyle = (\ell d_s)^{-1}\bigl(\sum_{j=1}^{d_s}\chi_j(m)+\ell
    d_sq\chi_d(m)\bigr)
  \\
  &\textstyle\le
    (\ell d_s)^{-1}\bigl(\sum_{j=1}^d\chi_j(m)+(\ell
    d_sq-1)\chi_d(m)\bigr) \\
  &\textstyle=
    (\ell d_s)^{-1}\lim_{t\to\infty}t^{-1} \bigl(\log |\det DX_t(x)|+(\ell d_sq-1)
    \log\|DX_t(x)\|\bigr)
  \\
  &\textstyle\le
    (\ell d_s)^{-1}\lim_{t\to\infty}t^{-1}\int_0^t \bigl(\Div DX(X_ux)+(\ell d_sq-1)
    \|DX(X_u(x))\|_2\bigr)\,du
  \\
  &\textstyle\le
    (\ell d_s)^{-1}\sup_{x\in\Lambda}\big\{\Div DX(x)+(\ell d_sq-1) \|DX(x)\|_2\bigr\}.
\end{align*}
By Definition~\ref{def:SD}(b), we again have that
$\lim_{t\to\infty}t^{-1} \eta_t(x)<0$ for $m$-almost every
$x\in\Lambda$. This completes the proof of the claim.
\end{proof}


\subsection{Examples of $q$-dissipative singular-hyperbolic attracting
  sets}
\label{sec:exampl-q-dissip}

We present some open classes of examples of vector fields satisfying
the assumptions of the Main Results.

\begin{example}
  Let $X: \RR^3 \to \RR^3$ defined by the classical Lorenz equations below
  \begin{align}
                   \label{eq:lorenz-equations}
  \begin{cases}
        \dfrac{dx_1}{dt} & =  10 (x_2 - x_1),\\
        \dfrac{dx_2}{dt} & =  28x_1 - x_2 - x_1 x_3,\\
        \dfrac{dx_3}{dt} & =  x_1x_2 - \dfrac{8}{3}x_3.
      \end{cases}
  \end{align}

  It is known that there exists an ellipsoid $E$ such that every
  positive trajectory of $X$ crosses $E$ transversely and never leaves
  it. In particular, we have that $E$ is a trapping region for
  $X$. Moreover, there exist three singularities for $X$ inside $E$,
  two with complex expanding eigenvalues and one Lorenz-like. See
  Figure \ref{fig:lorenz-original} and check, e.g., \cite[Section
  3.3]{AraPac07} for more details.

  \begin{figure}[!htb]
    \centering
    \includegraphics[height=4cm]{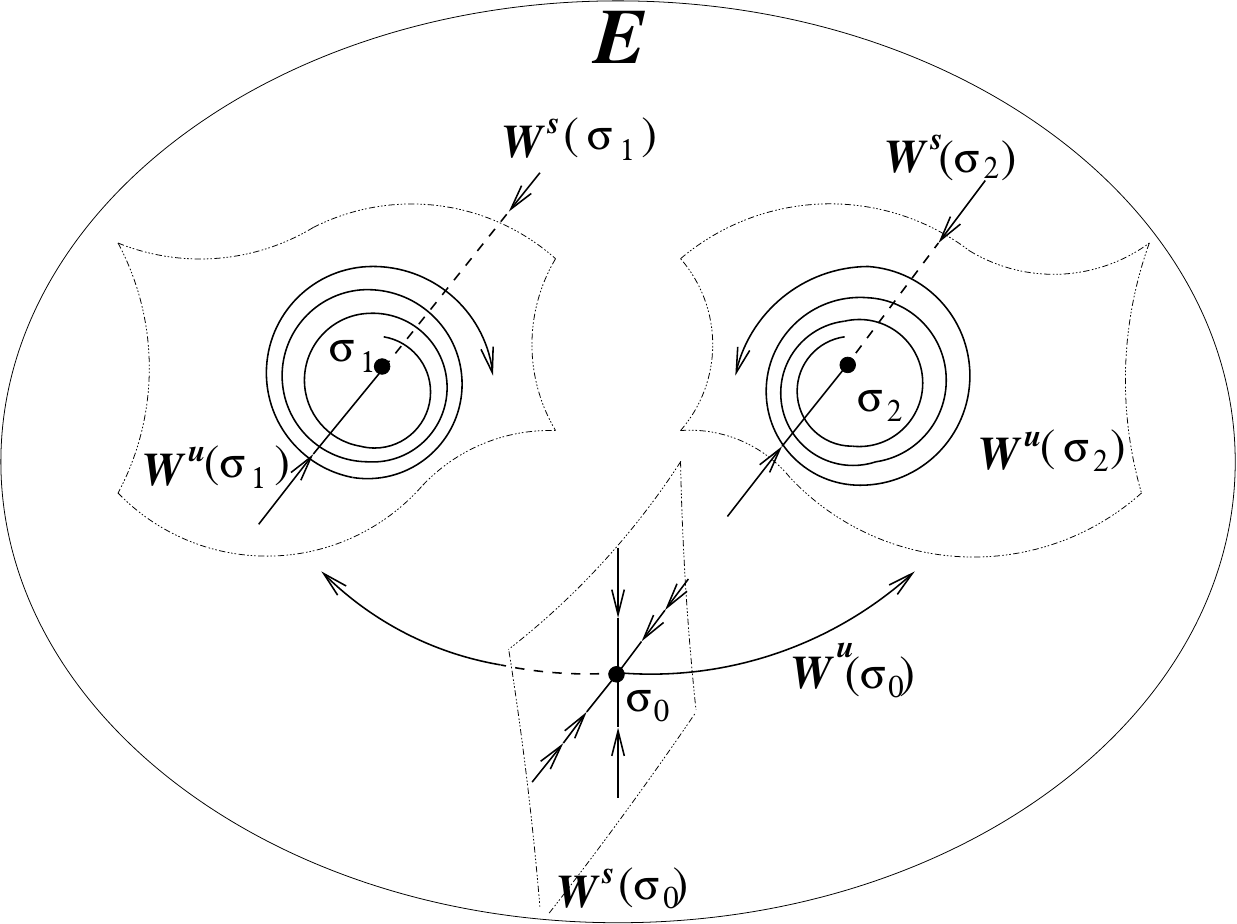}
    \caption{Trapping region for the original Lorenz equations with a
      Lorenz-like singularity $\sigma_0$ and two singularities non-Lorenz-like
      $\sigma_1$ and $\sigma_2$.}
    \label{fig:lorenz-original}
  \end{figure}
  
  The authors in \cite{ArMel17} proved that the flow of $X$ is
  $1.278$-strongly dissipative. As we see in Subsection
  \ref{sec:uni-condition} we may need to perturb $X$ to get the UNI
  condition. Thus we can apply our results in a neighborhood of a
  vector field arbitrarily close to $X$.
\end{example}

\begin{example}
  In \cite{AraPac07} the authors construct a singular-hyperbolic
  attracting set with three Lorenz-like singularities by modifying the
  geometric Lorenz attractor in the following way: first add two
  singularities $\sigma_1$ and $\sigma_2$ for the flow inside
  $W^u(\sigma)$ as in the left-hand side of
  Figure~\ref{fig:geo-lorenz-modified}.
  \begin{figure}[!htb]
    \centering
    \includegraphics[height=4cm]{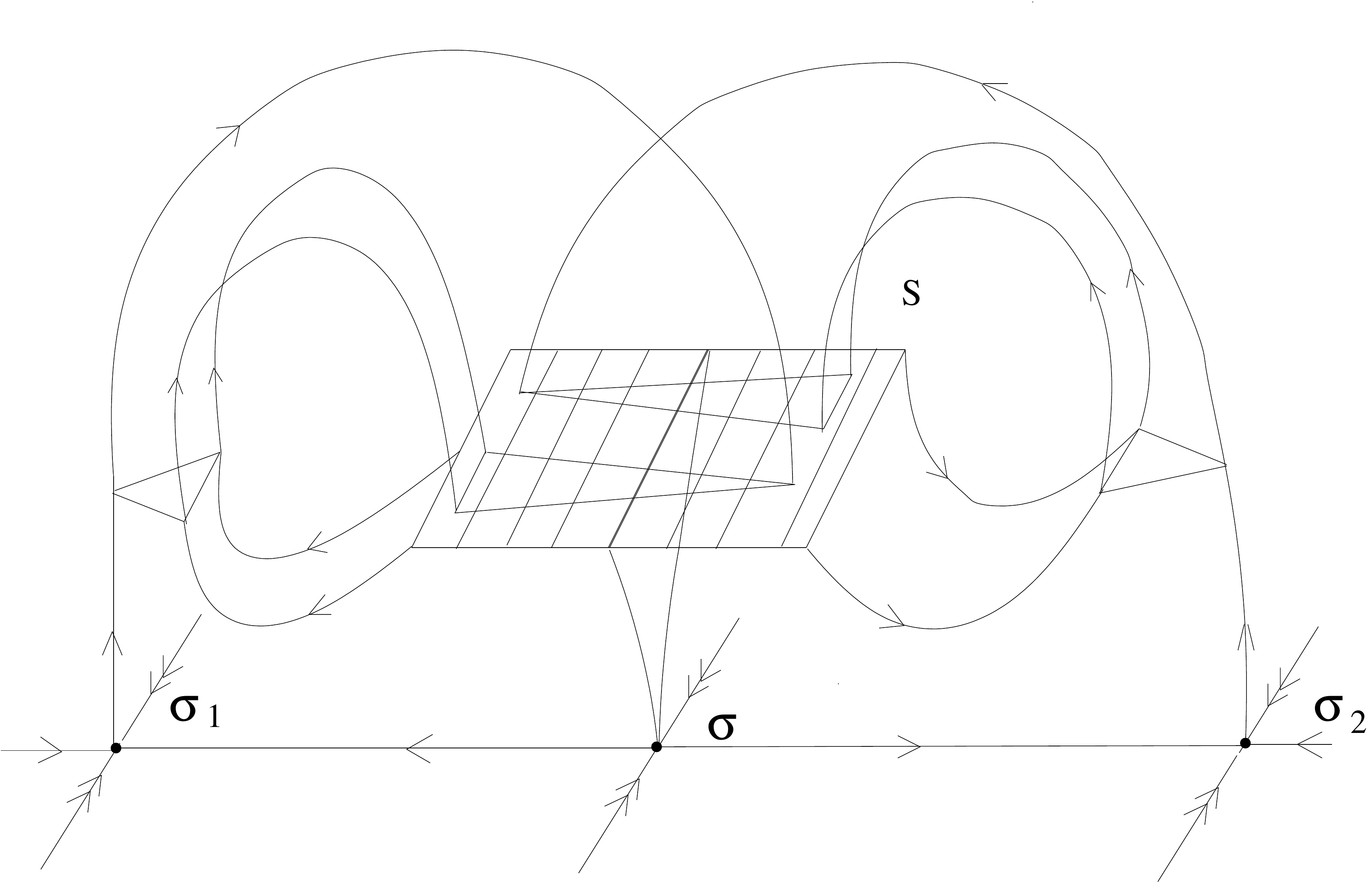}
      \includegraphics[height=5cm]{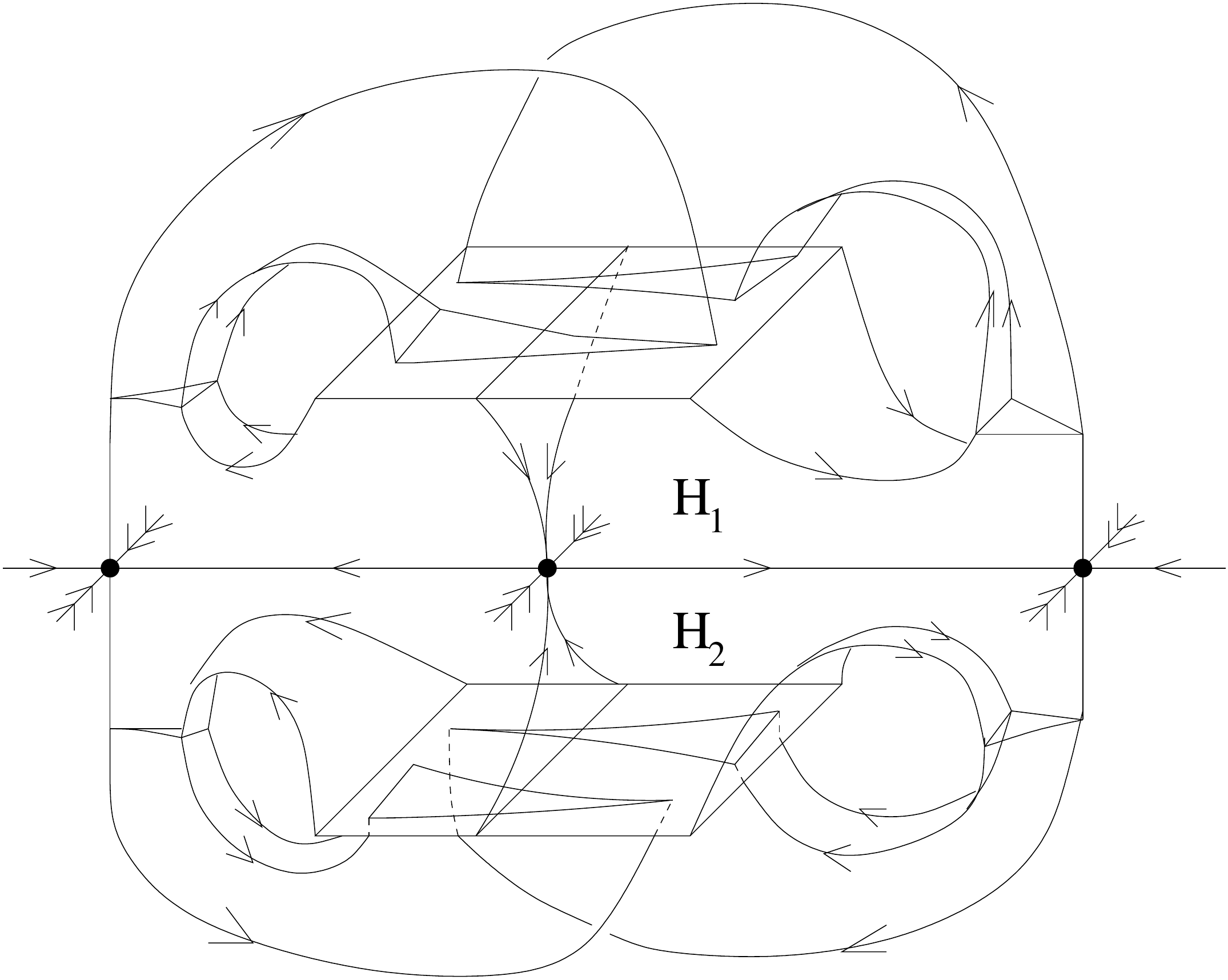}
      \caption{Construction of a singular-hyperbolic attracting set by
        modifying the Geometric Lorenz Attractor, obtaining an example
        of a singular-hyperbolic attracting set, non-transitive and
        containing three Lorenz-like singularities.}
      \label{fig:geo-lorenz-modified}
  \end{figure}
  As result of this construction we get a singular-hyperbolic
  attracting set, non-transitive, with three Lorenz-like
  singularities. The singularities can be chosen in the construction
  to satisfy the $q$-strongly dissipative condition. Moreover, the
  sets $H_1$ and $H_2$ in the right-hand side of
  Figure~\ref{fig:geo-lorenz-modified} are closed, invariant and
  transitive. It follows that each of them support a unique SRB
  measure for the flow. For more details of this construction check
  \cite[Section 9.1]{AraPac07}.
\end{example}

\begin{remark}\label{rmk:no-Lorenz}
  We could also include four complex expanding singularities on the
  ``lobes'' of Figure \ref{fig:geo-lorenz-modified} and transform this
  example in one containing non-Lorenz-like singularities.

  There are examples of singular-hyperbolic attracting sets whose
  singularities are all non-Lorenz-like; see e.g. \cite{Morales07} and
  references therein. Note that these examples become ``trivial''
  according to our definitions.
\end{remark}

\begin{example}
  Now we explain how to obtain an example of $q$-dissipative singular flow in
  higher dimension with $d_s > 1$. Let $X: \RR^3 \to \RR^3$ be given by the
  Lorenz equations \eqref{eq:lorenz-equations} and let $Y: \RR^k \to \RR^k$ be a
  smooth vector field admitting a singularity $\sigma$ which all its eigenvalues are
  negative (attractor). Denoting by $W$ the topological basin of this
  singularity and $U$ the topological basin for Lorenz attractor. Then,
  defining $Z: \RR^3\times \RR^k \to \RR^3 \times \RR^k$, by $Z(x,y) = (X(x), Y(y))$, we have
  that $U \times W$ is the topological basin for $\widehat{\Lambda} = \Lambda \times \{\sigma\}$,
  where $\Lambda$ is the Lorenz attractor.

  Denoting by $\lambda_3 < \lambda_2 < 0 < \lambda_1$ the eigenvalues for the
  singularity $0$ of $X$, we know that $\ell$ can be taken equal to $1$ (see the
  proof that $X$ is $1.278$-strongly dissipative in \cite[Section 5]{ArMel17}).
  Thus, if we choose the eigenvalues $\lambda_{k + 3} < \lambda_{k + 2} < \cdots
  < \lambda_4$ of $\sigma$ for $Y$ all close to $\lambda_3$, it follows that $\widehat{\Lambda}$
  it strongly dissipative singular-hyperbolic attracting with $d_s = 1 + k$ and
  $\ell$ arbitrarily close to $1$.
\end{example}


\subsection{Conjectures}
\label{sec:conjectures}

We propose some conjectures of results that may be obtained by
extending the techniques used in this text.

\subsubsection{No need for smoothness of the strong stable foliation}
\label{sec:no-smoothness-strong}

The assumption of constant return times along stable leaves, implicit
in Subsection~\ref{sec:hyperb-skew-product}, seems to be a feature of
the specific technical tools used in the proof.

We note that according to Lemma~\ref{le:hC1+} and
Theorem~\ref{thm:props-poincare-maps} the one-dimensional quotient map
$f$ is piecewise $C^{1+}$ smooth, independent of the smoothness of the
stable foliation.  This might be a starting point to an alternate
strategy to find a skew-product semiflow with the needed properties
and conjugated to the original flow, without assuming that the roof
function is constant on stable leaves.

\begin{conjecture}
   \label{conj:UNI1d}
   There exists a exponential mixing skew-product semiflow built over
   the expanding semiflow with a roof function which is non-constant
   on stable leaves, and semiconjugated to the original flow.
\end{conjecture}

\subsubsection{Uniform non-integrability holds for all
  singular-hyperbolic attracting sets}
\label{sec:uniform-non-integr}

Since, by Theorem~\ref{theorem:conjugacy-skew-product-hyperbolic}, we
obtain a finite collection of skew-product semiflows which are
semiconjugated to the flow on a neighborhood of the support of each
ergodic physical probability measure of our singular-hyperbolic
attracting set, we might obtain in general an attracting set having an
ergodic physical measure which mixes exponentially and another ergodic
physical measure with slow rate of mixing.

We conjecture that this is not possible. We note that the UNI
condition was obtained in~\cite{ArMel16} for Lorenz-like attractors
with a unique Lorenz-like singularity and ergodic physical probability
measure without perturbing the vector field -- in particular,
obtaining the exponential mixing property for the flow of the original
Lorenz equations. This should extend to the general case with finitely
many singularities.

\begin{conjecture}
  \label{conj:allUNI}
  The Uniform Non-Integrability (UNI) condition holds for all ergodic
  physical probability measures supported on each non-trivial
  singular-hyperbolic attracting set.
\end{conjecture}

\subsubsection{Exponential mixing for other equilibrium states}
\label{sec:exponent-mixing-othe}

We recall that Dolgopyat \cite{Do98}, in the work which first provided
the technical path to proving exponential decay for Anosov flows,
obtained exponential mixing for the physical/SRB measure under strong
assumptions on the smoothness of both the stable and unstable
foliations. In the same work, fast decay (in the sense of Schwarz,
that is, superpolynomial) was obtained for equilibrium states with
respect to H\"older continuous potentials with respect to
topologically mixing $C^\infty$ Anosov flows.

Recently Tsujii and Zhang~\cite{tsujii2020zhang} proposed a proof of
exponential mixing for all equilibrium states with respect to any
H\"older continuous potential of topological mixing $C^\infty$ Anosov
flows on $3$-manifolds.

\begin{conjecture}
   \label{conj:eqstates}
   The techniques from~\cite{tsujii2020zhang} can be adapted to
   singular flows to extend the results on this text for
   equilibrium state associated to H\"older continuous potentials.
\end{conjecture}

This naturally leads to extend the main tools of exponential mixing
for expanding semiflows to cover all such equilibrium states instead
of dealing only with absolutely continuous invariant measures.

Recently, in \cite{DaltroVar2021}, exponential mixing has been
obtained for all Gibbs measures (of which the absolutely continuous
invariant measure is a particular example) in the simplified setting
of suspension semiflows over full branch piecewise expanding
$C^{1+\alpha}$ maps with \emph{finitely} many branches. This was
extended in~\cite{DaltroVar2021a} to Markov $C^{1+\alpha}$ piecewise
expanding maps to obtain exponential mixing for each equilibrium state
of Axiom A attractors for $C^2$ flows with respect to any H\"older
continuous potential.

\subsubsection{Exponential mixing for higher dimensional
   sectional-hyperbolic attracting sets}
\label{sec:exponent-mixing-high}

Open examples of Anosov flows with exponential mixing
physical/SRB measures in arbitrary finite dimensional compact
manifolds were obtained by Butterley and War~\cite{BuWar20} exploring
the same techniques presented in this text.

If we relax the codimension $2$ condition on the stable bundle of
singular-hyperbolic attracting sets, that is, the assumption $\dim
E^{cu}=2$, then we have sectional-hyperbolic systems -- introduced by
Metzger and Morales in~\cite{MeMor08}.

It has been show \cite{LeplYa17} that sectional-hyperbolic attractors
have a unique physical measure and that, removing the transitivity
assumption, we still have finitely many ergodic physical measures
whose basins cover a full Lebesgue measure subset of the trapping
region; see \cite{Araujo19}.

In general the holonomies of the stable foliation in cross-sections
are no longer smooth, but only H\"older continuous in all higher
dimensional cases -- although these holonomies are still absolutely
continuous maps: this is a consequence of partial hyperbolicity for
sufficiently smooth ($C^2$) flows.

More specifically, a concrete example of a sectional-hyperbolic
attractor was provided by Bonatti, Pumari\~no and Viana in
\cite{BPV97}, also known as the multidimensional Lorenz attractor.

\begin{conjecture}
   \label{conj:BPVexpmix}
   The multidimensional Lorenz attractor is exponentially
   mixing. Moreover, this conclusion holds for an open and dense
   subset of all sectional-hyperbolic attracting sets.
\end{conjecture}


\section[Skew product semiflow from singular-hyperbolicity]{Global
  Poincar\'e return map for Singular-hyperbolic attracting sets}
\label{chapter:decay-correlations-suspension}

We recall some results from~\cite{ArMel17}.  These results hold for
general partially hyperbolic attracting sets with $d_{cu}\ge2$ and do
not depend on the existence of a dense forward orbit (transitivity).

\subsection{Properties of partially hyperbolic attracting sets}
\label{sec:propert-partially-hy}

In what follows we write $X_t$ for the flow generated by a $C^1$
vector field $X$ on a compact finite-dimensional manifold $M$ having
an attracting set $\Lambda$ with isolating neighborhood $U_0$:
$\Lambda=\cap_{t>0} X_t(U_0)$ and $\overline{X_t(U_0)}\subset U_0$ for
all $t\ge T_0$ for some $T_0>0$.

\begin{proposition}{\cite[Proposition~3.2 and Remark~3.3]{ArMel17}}
  \label{pr:Es} Let $\Lambda$ be a partially
  hyperbolic attracting set.  The stable bundle $E^s$ over $\Lambda$
  extends to a continuous uniformly contracting $DX_t$-invariant
  bundle $E^s$ over an open neighborhood of $\Lambda$.
\end{proposition}


We assume without loss of generality that $E^s$ extends as in
Proposition~\ref{pr:Es} to $U_0$.

Denoting by $B^k$ the $k$-dimensional open unit disk of $\RR^k$
endowed with the Euclidean distance induced by the Euclidean norm
$\|\cdot\|_2$.  Let $\mathrm{Emb}^2(B^k,M)$ denote the set of $C^2$
embeddings $\phi:B^k\to M$ endowed with the $C^2$ distance. Given
$\phi \in \mathrm{Emb}^2(B^k, M)$ we denote by
$\Lip(\phi) = \sup_{x \neq y} \big(\dist(\phi(x),\phi(y))/\|x -
  y\|_2\big)$ the Lipschitz constant of $\phi$.  We say that a subset
$D \subset M$ is a $C^2$ \emph{embedded $k$-dimensional disk} if there
exists $\phi \in \mathrm{Emb}^2(B^k, M)$ such that $\phi(B^k) = D$.

\begin{proposition}{\cite[Theorem~4.2 and Lemma~4.8]{ArMel17}}
  \label{pr:Ws}
  Let $\Lambda$ be a partially hyperbolic attracting set.  There
  exists a positively invariant neighborhood $U_0$ of $\Lambda$, and
  constants $C>0$, $\lambda\in(0,1)$, such that the following are
  true:
  \begin{enumerate}
  \item 
    For every point $x \in U_0$ there is a $C^r$ embedded
    $d_s$-dimensional disk $W^s_x\subset M$, with $x\in W^s_x$, such
    that $T_xW^s_x=E^s_x$ and for all $t>0$: $X_t(W^s_x)\subset
    W^s_{X_tx}$ and $d(X_tx,X_ty)\le C\lambda^t d(x,y)$ for all $y\in W^s_x$.

  \item The disks $W^s_x$ depend continuously on $x$ in the $C^0$
    topology: there is a continuous map
    $\gamma:U_0\to {\rm Emb}^0(\cD^{d_s},M)$ such that
    $\gamma(x)(0)=x$ and $\gamma(x)(\cD^{d_s})=W^s_x$.  Moreover,
    there exists $L>0$ such that $\Lip\gamma(x)\le L$ for all
    $x\in U_0$.
  \item The family of disks $\{W^s_x:x\in U_0\}$ defines a topological
    foliation of $U_0$.
  \end{enumerate}
\end{proposition}


The splitting $T_\Lambda M=E^s\oplus E^{cu}$ extends continuously to a
splitting $T_{U_0} M=E^s\oplus E^{cu}$ where $E^s$ is the invariant
uniformly contracting bundle in Proposition~\ref{pr:Es}.  (In
general, $E^{cu}$ is not invariant.)  Given $a>0$, we define the {\em
  center-unstable cone field},
\[
  \cC^{cu}_x(a)=\{v= v^s+v^{cu}\in E^s_x\oplus E^{cu}_x:\|v^s\|\le
  a\|v^{cu}\|\}, \quad x\in U_0.
\]

\begin{proposition}{\cite[Proposition~3.1]{ArMel17}} \label{pr:Ccu}
  Let $\Lambda$ be a partially hyperbolic attracting set.  There
  exists $T_0>0$ such that for any $a>0$ (after possibly shrinking
  $U_0$) we have $ DX_t\cdot \cC^{cu}_x(a)\subset \cC^{cu}_{X_tx}(a)$
  for all $t\ge T_0$, $x\in U_0$.
\end{proposition}


\begin{proposition}{\cite[Proposition 2.10]{ArMel18}} \label{pr:VE}
Let $\Lambda$ be a singular-hyperbolic attracting set.
After possibly increasing $T_0$ and shrinking $U_0$, there exist
constants $K,\theta>0$ so that
$|\det(DX_t| E^{cu}_x)|\ge K e^{\theta t}$ for all
$x\in U_0$, $t\geq 0$.
\end{proposition}

\subsubsection{The stable lamination is a topological foliation}
\label{sec:stable-laminat-topol}

Proposition \ref{pr:Ws} ensures the existence of an $X_t$-invariant
\emph{stable lamination} $\cW^s_\Lambda$ consisting of smoothly
embedded disks $W^s_x$ through each point $x\in\Lambda$.  Although not
true for general partially hyperbolic attractors, for
singular-hyperbolic attractors in our setting $\cW^s_\Lambda$ indeed
defines a topological foliation in an open neighborhood of $\Lambda$.

\begin{theorem}{\cite[Theorem 5.1]{ArMel18}} \label{thm:stabltopfol}
  Let $\Lambda$ be a singular-hyperbolic attracting set.  Then the
  stable lamination $\cW^s_\Lambda$ is a topological foliation of an
  open neighborhood of~$\Lambda$.
\end{theorem}


From now on, we refer to $\cW^s_\Lambda=\{W^s_x:x\in\Lambda\}$ as the
\emph{stable foliation}.

\subsubsection{Absolute continuity of the stable foliation}
\label{sec:abscontWs}

From now on we assume that the vector field $X$ is of class $C^2$.
Let $Y_0,\,Y_1\subset U_0$ be two smooth disjoint $d_{cu}$-dimensional
disks that are transverse to the stable foliation $\cW^s_\Lambda$.
Suppose that for all $x\in Y_0$, the local stable leaf $W^s_x$ intersects
each of $Y_0$ and $Y_1$ in precisely one point. The {\em stable
  holonomy} $H:Y_0\to Y_1$ is given by defining $H(x)$ to be the
intersection point of $W^s_x$ with $Y_1$.

A key fact for us is regularity of stable holonomies.

\begin{theorem}{\cite[Theorem 6.3]{ArMel18}}\label{thm:H}
  The stable holonomy $H:Y_0\to Y_1$ is absolutely
  continuous. That is, $m_1\ll H_*m_0$ where $m_i$ is
  Lebesgue measure on $Y_i$, $i=0,1$.  Moreover, the
  Jacobian $JH:Y_0\to\RR$ given by
  \begin{align*}
    JH(x)=\frac{dm_1}{dH_*m_0}(Hx)=
    \lim_{r\to0}\frac{m_1(H(B(x,r)))}{m_0(B(x,r))},\quad x\in Y_0,
  \end{align*}
  is bounded above and below and is $C^\alpha$ for some
  $\alpha > 0$.
\end{theorem}


Hence, we can assume without loss of generality, that there
exists a foliation $\cW^s_\Lambda$ of $U_0$, which continuously
extends the stable lamination of $\Lambda$ together with a
positively invariant field of cones
$(\cC^{cu}_x)_{x\in U_0}$ on $T_{U_0}M$.  Moreover, the
Jacobian of holonomies along contracting leaves on
cross-sections of singular-hyperbolic attracting sets in our
setting is a H\"older function.
It is well-known that the $C^2$ smoothness of $X$ is crucial to these
properties since the work of Anosov~\cite{An67}.

\subsection{Global Poincar\'e return map}
\label{sec:global-poincare-retu}

In \cite{APPV} the construction of a global Poincar\'e map for any
singular-hyperbolic attractor is carried out based on the existence of
``adapted cross-sections'' and $C^{1 + \alpha}$ stable holonomies on
these cross-sections. With the results just presented this
construction can be performed for any singular-hyperbolic attracting
set. This construction was presented in~\cite[Sections 3 and
4]{ArMel18}, so from there we obtain:
\begin{itemize}
\item a finite collection
  $\Xi=\Sigma_1+\dots+\Sigma_m$\footnote{We write $A+B$ the
    union of the disjoint subsets $A$ and $B$.} of (pairwise
  disjoint) cross-sections for $X$ so that
  \begin{itemize}
  \item each $\Sigma_i$ is diffeomorphically identified with
    $(-1,1)\times B^{d_s}$;
  \item the \emph{stable boundary}
    $\partial^s\Sigma_i\cong \{\pm1\}\times B^{d_s}$
    consists of two curves contained in stable leaves; and
  \item each $\Sigma_i$ is foliated by
    $W^s_x(\Sigma_i)= \bigcup_{|t|<\varepsilon_0}X_t(W^s_x)\cap
    \Sigma_i$ for a small fixed $\varepsilon_0 > 0$. We denote this
    foliation by $W^s(\Sigma_i), i=1,\dots,m$;
  \end{itemize}
\item a Poincar{\'e} map
  $P:\Xi\setminus\Gamma\to\Xi, P(x)=X_{\tau(x)}(x)$ with
  $\tau:\Xi\setminus\Gamma\to [T,+\infty)$ the associated return time,
  which is $C^2$ smooth in $\Sigma_i\setminus\Gamma, i=1,\dots,m$;
  preserves the foliation $W^s(\Xi)$ and a big enough time $T>0$,
  where $\Gamma=\Gamma_0\cup\Gamma_1$ is a finite family of stable
  disks $W^s_{x_i}(\Xi)$ so that
  \begin{itemize}
  \item
    $\Gamma_0=\{x\in\Xi:X_{T+1}(x)\in\bigcup_{\sigma\in
      S}\gamma^s_\sigma\}$ for
    $S=S(X,\Lambda)=\{\sigma\in\Lambda: X(\sigma)=0\}$
    and $\gamma_\sigma^s$ is the local stable manifold of
    $\sigma$ in a small fixed neighborhood of $\sigma\in S$;
    and
  \item
    $\Gamma_1=\{x\in\Xi:P(x)\in\partial^s\Xi=\cup^m_{i = 1}\partial^s\Sigma_i\}$;
  \end{itemize}
\item and open neighborhoods $V_\sigma$ for each $ \sigma \in S$ so that
  defining $V_0 = \cup_{\sigma \in S}V_\sigma$ we have that  every orbit of
  a regular point $z\in U_0\setminus V_0$ eventually hits $\Xi$ or
  else $z\in\Gamma$.
\end{itemize}

Having this, the same arguments from \cite{APPV} (see
\cite[Proposition 4.1 and Theorem 4.3]{ArMel18} and \cite[Section
2.3]{Araujo19}) show that $DP$ contracts $T_\Xi W^s(\Xi)$ and expands
vectors on the unstable cones
$\{C^u_x(\Xi) := C^{cu}_x(a)\cap
T_x\Xi\}_{x\in\Xi\setminus\Gamma}$. The stable holonomies for $P$
enable us to reduce its dynamics to a one-dimensional map, as follows.

Let $\Sigma$ be a cross-section in $\Xi$.  A smooth curve $\gamma: [0,1] \to \Sigma$ is called a $u$-\emph{curve} if
$\gamma'(t) \in C^u_{\gamma(t)}(\Xi)$ for all $t \in [0,1]$. We say that the
$u$-curve $\gamma$ \emph{crosses} $\Sigma$ if each leaf $W^s_x(\Sigma)$ of $\Sigma$ intersects $\gamma$
in a unique point.

Let $\gamma_i\subset\Sigma_i$\footnote{We also use the term \emph{curve} to
  denote the image of the curve.} be $u$-curves that cross $\Sigma_i$, $i = 1, 2,
\dots, m$. The \emph{(sectional) stable holonomy}
$\pi_\gamma:\Xi\to\gamma=\sum^m_{i = 1}\gamma_i$ is defined by setting $\pi_\gamma(x)$ to be the
intersection point of $W^s_x(\Sigma_i)$ with $\gamma_i$, for 
$x\in\Sigma_i$ and $i=1, 2, \dots, m$.

\begin{lemma}{\cite[Lemma 7.1]{ArMel18}}\label{le:hC1+}
  The stable holonomy $\pi_\gamma$ is $C^{1+\alpha}$ for some $\alpha>0$.
\end{lemma}

Following the same arguments in \cite{APPV} (see also \cite[Section
7]{ArMel18}) we obtain a one-dimensional piecewise $C^{1+\alpha}$
quotient map over the stable leaves
$f_\gamma:\gamma\setminus\Gamma\to \gamma$ for some $0<\alpha<1$ so that
$\pi_\gamma(P(x)) = f_\gamma(x)$ and $|f'_\gamma (x)|>2$, for all $x \in \gamma
\backslash \Gamma$.

Let $\gamma_i: I_i \to \Sigma_i$ be a smooth parametrization of a $u$-curve in $\Sigma_i$ for each $i = 1, 2,
\dots, m$. We assume that $\{I_i: i = 1, 2, \dots, m\}$ is a family of disjoint
intervals and define $I = I_1 \cup I_2 \cup \cdots I_m$. We define a
parametrization of $\gamma$ as $\gamma: I \to \Xi$ by $\gamma(t) = \gamma_i(t)$
if $t \in I_i$. Using the last parametrization we can identify $f_\gamma$ with
the one-dimensional map $f: I \backslash \cD \to I$ by $f(x) =
\gamma^{-1}(f_\gamma(\gamma(x)))$, where $\cD =
\gamma^{-1}(\pi_\gamma(\Gamma))$ is the \emph{critical set} for $f$. Moreover,
defining the \emph{singular set} $\cS =
\gamma^{-1}(\pi_\gamma(\Gamma_0))$ we get, as shown in \cite[Proof of Lemma
8.4]{ArMel18}, that $f'|_{I \backslash \cD}$ behaves like a power of the
distance near $\cS$ in the following sense: there exist constants
$\eta \in (0,1)$ and $C, q > 0$ such that

\begin{enumerate}[(C1)]
\item $C^{-1}d(x,\cS)^q \leq |f'(x)| \leq Cd(x,\cS)^{-q}$, for all $x
  \in I \backslash \cS$;
\item $|\log|f'(x)| - \log|f'(y)|| \leq C|x - y|^\eta(|f'(x)|^{-q} +
  |f'(x)|^q)$, for all $x, y \in I \backslash \cS$, with $|x - y| <
  d(x,\cS) / 2$.
\end{enumerate}

\begin{remark}
  With the identifications above and in order to simplify notations, we
  sometimes make no distinction between $I$ and $\gamma$, $f$ and $f_\gamma$,
  $\cD$ and $\pi_\gamma(\Gamma)$ and $\cS$ and
  $\pi_\gamma(\Gamma_0)$. We assume in what follows that $I=[0,1]$.
\end{remark}

\begin{remark}[Quotient maps are conjugated]\text{}
  \begin{enumerate}[(a)]
  \item For $j = 1, 2$ let $\gamma^j = \sum_i \gamma^j_i$, where
    $\gamma^j_i$ is a $u$-curve in $\Sigma_i$. If
    $f_j: \gamma^j \backslash \Gamma \to \gamma^j$ are two quotients
    along stable leaves (as explained above), then they are
    $C^{1 + \alpha}$ conjugated. Indeed, let $\pi_j$ be the stable
    holonomy with respect to $\gamma^j$. Defining
    $g: \gamma^1 \to \gamma^2$ by $g = \pi_2|_{\gamma^1}$ it follows
    that $g$ is a $C^{1 + \alpha}$ diffeomorphism. We claim that $g$
    is a conjugacy between $f_1$ and $f_2$. By the invariance of the
    stable leaves under the Poincar\'e map, we have that
    $P(g(x)) \in W^s_{P(x)}(\Xi)$ for all $x \in \gamma^1$. Hence
    $g(f_1(x)) = g(\pi_1(P(x))) = \pi_2(P(g(x))) = f_2(g(x))$ for all
    $x \in \gamma^1$.
    
    
  \item Moreover, it follows from (a) that there exists a constant $C > 0$,
    depending only on the holonomy map $g$, such that
    \begin{eqnarray}\label{eq:comp-unidimensional-maps}
      C^{-1}|f_2(g(x_1)) - f_2(g(x_2))| \leq |f_1(x_1) - f_1(x_2)| \leq C|f_2 (g(x_1)) - f_2(g(x_2))|,
    \end{eqnarray} for all $x_1, x_2 \in \gamma^1_i$, $i = 1, 2, \dots, m$.
  \end{enumerate}
\end{remark}

For $0 < \delta < 1$ we define the \emph{smooth $\delta$-truncated
  distance of $x$ to $\cD$ on $I$} by

$$
\dist_\delta(x, \cD), = 
\left\{\begin{array}{lcl}
         \dist(x, \cD), & \text{if} & 0 < \dist(x, \cD) \leq \delta\\
         \left( \frac{1 - \delta}{\delta}\right)\dist(x, \cD) + 2\delta - 1, &
                                                                               \text{if}
                                    & \delta < \dist(x, \cD) < 2\delta\\
         1, & \text{if} & \dist(x, \cD) \geq 2\delta,
       \end{array}\right.
     $$
     where $\dist$ denotes the Euclidean distance in the interval $I$ here.

Given $\delta > 0$, let $B(\Gamma, \delta) = \{x \in \Xi: \dist(x, \Gamma) <
\delta\}$ and $\chi_{B(\Gamma, \delta)}: \Xi \to \{0,1\}$ be the
\emph{characteristic function} of $B(\Gamma, \delta)$.
We say that a function $\varphi: \Xi_0 = \Xi \setminus \Gamma
\to \RR$ has \emph{logarithmic growth near $\Gamma$} if there is a constant
$C=C(\varphi) > 0$
such that for every small $\delta$ it holds $|\varphi(x)| \cdot \chi_{B(\Gamma,\delta)}(x) \leq C
|\log\dist_\delta(\pi_\gamma(x), \cD)|$, for all $x \in \Xi_0$.

The construction outlined above can be summarized as in \cite[Theorem
2.8]{ArSzTr} as follows:

\begin{theorem}{\cite[Theorem 2.8]{ArSzTr}}\label{thm:props-poincare-maps}
  Let $X \in \fX^2(M)$ be a vector field admitting a non-trivial
  connected singular-hyperbolic attracting set $\Lambda$. Then there
  exists $\alpha > 0$, a finite family $\Xi$ of cross-sections and a
  global Poincar\'e map $P: \Xi_0 \to \Xi$, $P(x) = X_{\tau(x)}(x)$
  such that
  \begin{enumerate}
  \item \label{item:global-cross-section} the domain
    $\Xi_0=\Xi\setminus\Gamma$ is the entire cross-sections with a
    family $\Gamma$ of finitely many smooth arcs removed and
    \begin{enumerate}
    \item
      \label{subitem:logarithmic-growth}$\tau:\Xi_0\to[\tau_0,+\infty)$
      is a smooth function with logarithmic growth near $\Gamma$ and
      bounded away from zero by some uniform constant $\tau_0>0$;
    \item there exists a constant $C > 0$ so that
      $|\tau(x)-\tau(y)| < C\dist(x,y)$ for all points $y \in W^s_x(\Xi)$;
    \end{enumerate}
  \item We can choose coordinates on $\Xi$ so that the map $P$ can be
    written as $F:\widetilde Q \to Q$, $F(x,y)=(f(x), g(x,y))$, where
    $Q= I\times B^{d_s}$, $ I=[0,1]$ and
    $\widetilde Q = Q \setminus \widetilde\Gamma$ with
    $\widetilde\Gamma = \cD\times B^{d_s}$ and
    $\cD=\{c_1,\dots,c_n\}\subset I$ a finite set of points.
  \item The map $f: I\setminus\cD\to I$ is a piecewise $C^{1+\alpha}$
    map with finitely many branches, defined on the connected
    components of $I\setminus\cD$, with finitely many ergodic
    absolutely continuous invariant probability measures
    $\mu^i_f, i=1, \dots, k$, whose ergodic basins cover $I$ Lebesgue
    modulo zero. Also
    \begin{enumerate}
    \item $\inf\{|f'x|:x\in I\setminus\cD\}>2$;
    \item each $c\in\cD$ has a well-defined \emph{one-sided critical
        order}: there exist $\delta > 0$ and numbers
      $0 < \alpha^\pm(c) \leq 1, \kappa_\pm(c)>0$
      satisfying:
      $|f(x) - f(c)| \le\kappa_+(c) |x - c|^{\alpha^+(c)},\ \ |f'(x)|
      \le\kappa_+(c) |x - c|^{\alpha^+(c) - 1}$ for
      $x \in (c, c + \delta)$; and
      $|f(x) - f(c)| \le\kappa_-(c) |x - c|^{\alpha^-(c)},\ \ |f'(x)|
      \le\kappa_-(c) |x - c|^{\alpha^-(c) - 1}$ for
      $x \in (c, c - \delta)$;
    \item $1/|f'|$ has universal bounded $p$-variation\footnote{See
        \cite{Ke85} for the definition of $p$-variation.};
      and $d\mu^i_f/dm$ has bounded $p$-variation for some $p>0$.
    \end{enumerate}
  \item The map $g:\widetilde Q\to B^{d_s}$ preserves and uniformly
    contracts the vertical foliation $\cF=\{\{x\}\times B^{d_s}\}_{x\in
      I}$ of $Q$: there is $\lambda\in(0,1)$ so that $
    \dist(g(x,y_1),g(x,y_2)) \leq \lambda \cdot |y_1-y_2|$ for each $x \in I
    \backslash \cD$ and $y_1,
    y_2 \in B^{d_s}$.
  \item The map $F$ admits a finite family of physical ergodic
    probability measures $\mu^{i}_{F}$ which are induced by $\mu^i_f$
    in a standard way\footnote{See
      Proposition~\ref{prop:invariant-measure-for-skew-products} and
      \cite[Section 6.1]{APPV} where it is shown how to get
      $(\pi_\gamma)_*\mu^i_F=\mu^i_f$.}. Moreover, the Poincar\'e time
    $\tau$ is integrable both with respect to each $\mu^{i}_{f}$ and
    with respect to the two-dimensional Lebesgue area measure of $Q$.
  \item The subset\footnote{The subset $\cS$ can be identified with
      $h(\Gamma_0)$ while $\cD\setminus\cS$ can be identified with
      $h(\Gamma_1)$} $\cS = \{c \in \cD: 0 \leq \alpha^\pm(c) < 1\}$
    (of singular points) is nonempty and satisfies:
    \begin{enumerate}
    \item there exists $T_0 \in \NN$ such that for all $c \in \cD$
      there is $T = T(c) \leq T_0$ such that $f^T(c) \in \cS$;
    \item there exists $\delta > 0$ such that given $c \in \cD$, for
      all $0 < j < T=T(c)$ there exists $d \in \cD$ such that
      $f^j(c) \in (d - \delta, d + \delta) \backslash \{d\}$;
    \item there exist $\varepsilon, \delta > 0$ such that
      $f|_{(c, c + \delta)}$ is a diffeomorphism into the interval
      $(f^{T(c)}(c), f^{T(c)}(c) + \varepsilon)$ and the same holds
      true for the left neighborhoods $(c - \delta, c)$ and
      $(f^{T(c)}(c) - \varepsilon, f^{T(c)}(c))$;
    \item $c\in\cS\ \iff \lim_{t\to c}|t-b|^{1-\alpha(c)}\cdot|f'(t)|$
      exists and is finite;
    \item $c\in\cD\setminus\cS \iff$ the limit $\lim_{t\to c}|f'(t)|$
      exists and is finite.
    \end{enumerate}
  \end{enumerate}
\end{theorem}

\subsection{Constant Poincar\'e return time on stable leaves}
\label{sec:ensuring-that-poinca}

Now we explain how to ensure that the Poincar\'e return time of the
previous construction is constant on stable leaves.

\subsubsection{$q$-dissipativity and cross-sections}
\label{sec:foliat-cross-section}

\begin{figure}[!htb]
  \centering
  \includegraphics[height=5cm]{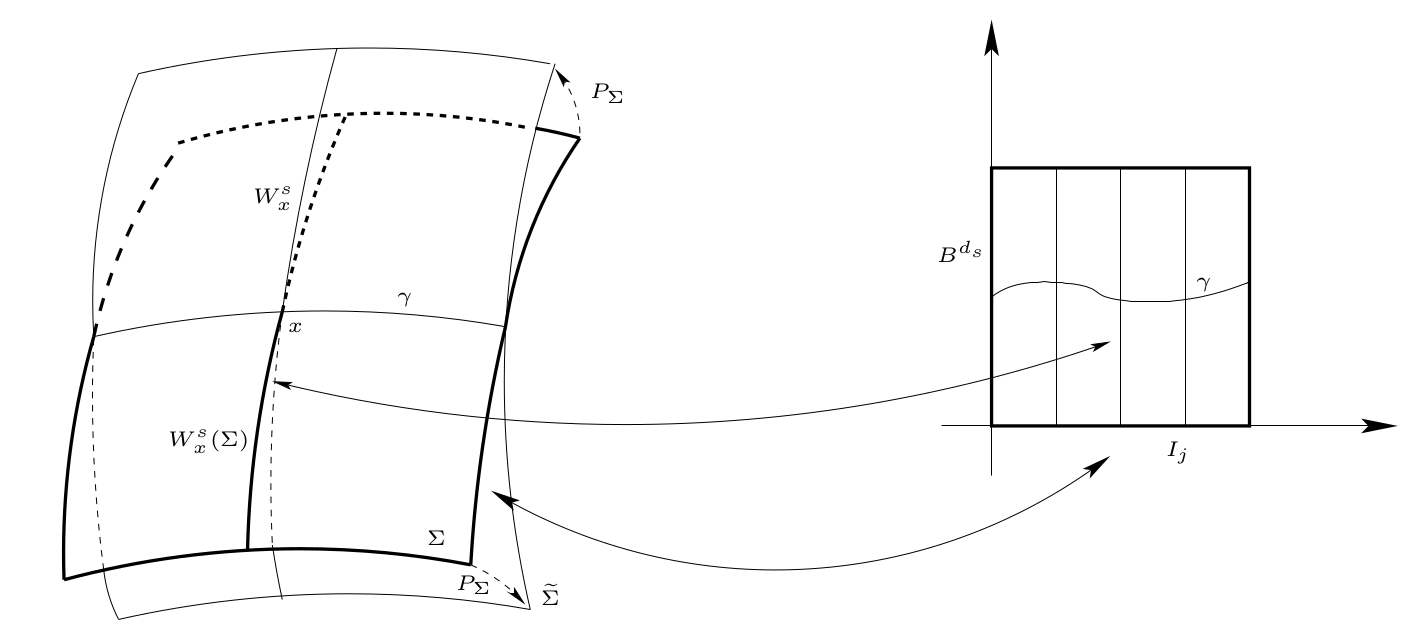}
  \caption[Stable foliation.]{On the right we have the identification
    of each cross-section $\Sigma$ with $I_j \times B^{d_s}$, for some
    $j = 1, 2, \dots, m$. On the left, we see how to change from
    $\Sigma$ to $\widetilde\Sigma$, obtaining a new smooth
    cross-section close to the previous one and foliated by the discs
    $W^s_x$.}
  \label{fig:section-stable-leaves}
\end{figure}

Using Theorem~\ref{thm:smooth-foliation} we may assume that $W^s_x(\Sigma)
\subset W^s_x$ for all cross-sections $\Sigma \subset \Xi$ and all $x \in
\Sigma$. Indeed, letting $\gamma \subset \Sigma$ be a $u$-curve the
cross-section $\widetilde{\Sigma} := \cup_{x \in \gamma}W^s_x$ is a submanifold
of class $C^q$. Moreover, if the disks $W^s_x$ have diameter small enough we have that the
Poincaré map $P_\Sigma: \Sigma \to \widetilde{\Sigma}$ between $\Sigma$ and
$\widetilde{\Sigma}$ is a diffeomorphism and has Poincaré time close to zero
(see Figure \ref{fig:section-stable-leaves}). With these considerations we
change each cross-section $\Sigma$ in $\Xi$ by the cross-section
$\widetilde{\Sigma}$ constructed as above.

\begin{remark}[Constant Poincar\'e time on stable leaves]
  As consequence of the change in the cross-sections, we now have that
  the Poincaré time $\tau$ is smooth and constant on stable leaves of
  $\Xi$.
\end{remark}

\subsubsection{Properties of the global Poincar\'e return time}
\label{sec:propert-global-poinc}

In what follows we state the linearization result of \cite{Newhouse16}
in the particular case of a saddle singularity for a $2$-dimensional
flow.

\begin{lemma}\cite[Theorem 1.5]{Newhouse16}
  \label{lemma:newhouse-linearization}
  Let $M$ be a surface and $X \in \fX^{1 + \alpha}(M)$, with
  $0 < \alpha < 1$. If $\sigma \in M$ is a singularity of saddle type
  for $X$ and $L = DX(\sigma)$, then there are a neighborhood
  $V \subset M$ of $\sigma$, a real number $\beta \in (0, \alpha)$ and
  a $C^{1 + \beta}$ diffeomorphism $h$ from $V$ onto its image such
  that $h(\sigma) = 0$ and $h(X_t(x)) = L_t(h(x)),$ for all
  $t \in \RR$ such that $V \cap X_{-t}(V) \neq \varnothing$ and all
  $x \in V \cap X_{-t}(V)$.
\end{lemma}

Every Lorenz-like singularity $\sigma$ admits a local central-unstable
invariant manifold $W=W_\sigma$ in a neighborhood of $\sigma$, as
smooth as the vector field $X$, such that
$T_\sigma W=E^{cu}_\sigma=E^u_\sigma\oplus E^s_\sigma$, where
$E^u_\sigma$ and $E^s_\sigma$ are the eigenspaces of $DX(\sigma)$
corresponding to the positive and least negative eigenvalues of
$DX(\sigma)$; see e.g.  We may assume without loss of generality that
$TW\subset C^{cu}$, that is, $W$ is a central-unstable two-dimensional
submanifold. Hence, we may apply
Lemma~\ref{lemma:newhouse-linearization} to $X\mid_W$ where the
singularity $\sigma$ becomes a two-dimensional hyperbolic saddle
singularity.

We now deduce some properties of the Poincar\'e return time
which will be useful in what follows.

\begin{lemma}
  \label{lemma:property-poincare-time}
  Let $x, y \in I\setminus\cD$ such that there
  is no element of $\cD$ between $x$ and $y$.  Then
  there exist constants $\alpha, C > 0$ so that
  \begin{align*}
|\tau(x) - \tau(y)|
  \leq C \left( \frac{|x - y|}{\min\{\dist(x, \Gamma), \dist(y,
    \Gamma)\}} + |x - y|^\alpha \right)
    \;\&\;
    |\tau'(x)| \le \frac{C}{\dist(x, \Gamma)}.
  \end{align*}
\end{lemma}

\begin{proof}
  By assumption, $X_{[0, \tau(x)]}(x)$ and $X_{[0, \tau(y)]}(y)$ hit
  the same cross-sections because there is no singular point between
  $x$ and $y$. Thus there are $\Sigma_1,\Sigma_2\in\Xi$ so that
  $x,y\in\Sigma_1$ and $P(x), P(y)\in\Sigma_2$.

  If $X_{[0, \tau(x)]}(x)$ and $X_{[0, \tau(y)]}(y)$ do not intersects
  $V_0$, then it follows that we can find a constant $C>0$ so that
  $|\tau(x) - \tau(y)| \le C |x - y|$ and we are done with the first
  inequality of the statement.

  Otherwise, $\xi_x=X_{[0, \tau(x)]}(x)$ and
  $\xi_y=X_{[0, \tau(y)]}(y)$ intersect $V_\sigma$, for some
  Lorenz-like singularity $\sigma \in \Lambda$.  We choose ingoing and
  outgoing cross-sections $\Sigma_i^\sigma, i=1,2$ of the flow inside
  $V_\sigma$ so that the trajectories $\xi_x,\xi_y$ cross
  $\Sigma_1^\sigma$ after leaving $\Sigma_1$ and $\Sigma_2^\sigma$
  before arriving at $\Sigma_2$; see
  Figure~\ref{fig:poincare-time-2}. We can construct these
  cross-sections as unions of strong-stable leaves as explained in
  Subsection~\ref{sec:foliat-cross-section}. We write
  $\gamma_\sigma^s$ for the local stable manifold of the equilibrium
  $\sigma$; see Figure~\ref{fig:poincare-time-2}.

  Let $\lambda_2 < 0 < -\lambda_2 < \lambda_1$ be the
  eigenvalues of $DX(\sigma)|_{E^{cu}_\sigma}$ and fix some
  central-unstable manifold $W=W_\sigma$ to which we will apply
  Lemma~\ref{lemma:newhouse-linearization}.
  We define
  $\alpha = \min\{-\lambda^\sigma_2 / \lambda^\sigma_1:\ p \in S(X,
  \Lambda)\} < 1$ a lower bound on the contraction/expansion ratio
  over all Lorenz-like singularities of the attracting set.

  Using Lemma~\ref{lemma:newhouse-linearization}, we smoothly
  linearize the flow in a neighborhood of $\sigma$ inside
  $V_\sigma \cap W$. Through the corresponding coordinate change, we
  can find $u$-curves $\gamma_1 = \{(x_1, 1):\ |x_1| \leq 1\}$ and
  $\gamma_2 = \{(\pm 1, x_2):\ |x_2| \leq 1\}$ inside
  $V_\sigma \cap W$ so that the Poincar\'e map
  $P_1 : \gamma_1 \backslash \gamma^s_\sigma \to \gamma_2$ is
  explicitly given by
  $P_1(x_1, 1) = (\sgn(x_1), |x_1|^{-\lambda_2 / \lambda_1})$.

  Let $g_1 : \gamma_1 \to [0, +\infty)$ be the first-hitting time
  function between the $u$-curves $\gamma$ and $\gamma_1$ and
  $g_2: \gamma_2 \to [0, +\infty)$ be the first-hitting time function
  between the $u$-curves $\gamma_2$ and $\gamma$. We assume without
  loss of generality that there is no singularity for the flow between
  $\gamma$ and $\gamma_1$ and between $\gamma_2$ and $\gamma$. We have
  that $g_1$ and $g_2$ are bounded functions of class $C^2$.

  Let $\tau_1: \gamma_1 \to [0, +\infty)$ be the first-hitting time
  between the $u$-curves $\gamma_1$ and $\gamma_2$. Using the choice
  of coordinates inside $V_\sigma \cap W$ we have that
  $\tau_1(x_1, 1) = - (\lambda_1)^{-1} \log |x_1|$; see Figure
  \ref{fig:poincare-time-2}. Since $\Sigma_i,\Sigma_i^\sigma, i=1,2$
  are unions of stable leaves, then \emph{the Poincar\'e times are
    constant on stable leaves and so we can deduce properties of the
    global Poincar\'e return time through the functions $g_1,g_2$ and
    $\tau_1$}.   For $z \in \{x, y\}$ we have
  \begin{align}\label{eq:g_12tau}
  \tau(z) = g_1(z) + \tau_1(X_{g_1(z)}(z)) +
  g_2(X_{\tau_1(z)}(X_{g_1(z)}(z))).
  \end{align}
   \begin{figure}[!htb]
    \centering
    \includegraphics[height=8cm]{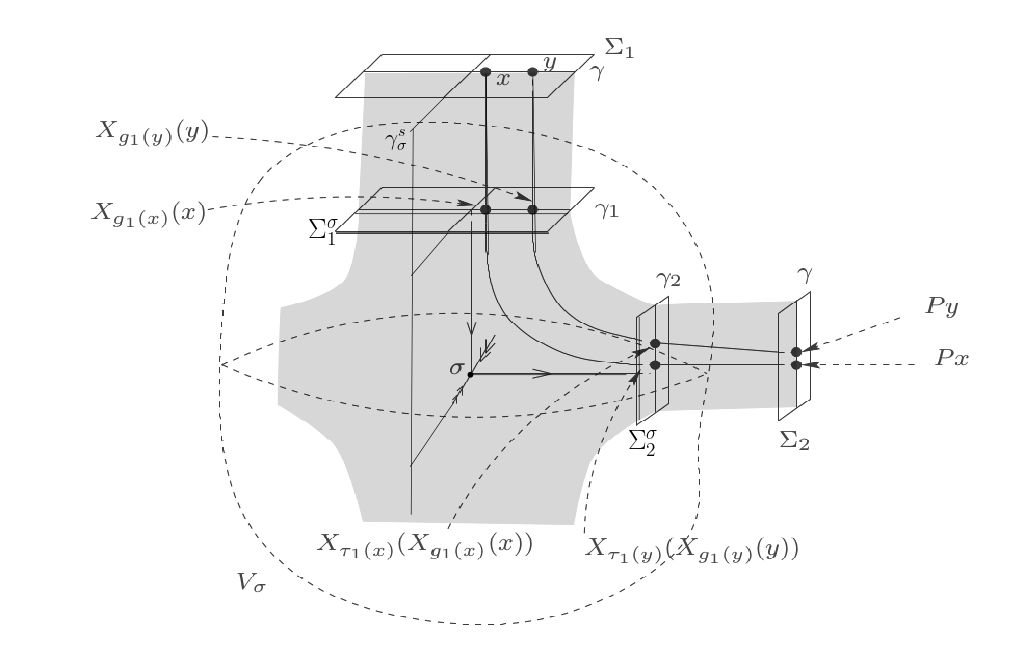}
    \caption{Linearization along a central-unstable surface on a
      neighborhood of a Lorenz-like singularity.}
    \label{fig:poincare-time-2}
  \end{figure}
  Because $g_1$ is $C^2$ we have that $|g_1(x) - g_1(y)|$ is bounded
  above by a constant times $|x-y|\le|x-y|^\alpha$.  It also follows
  that $|\tau_1(X_{g_1(x)}(x)) - \tau_1(X_{g_1(y)}(y))|$ equals
  \begin{align*}
    (\lambda_1)^{-1}
    &\Big|\log\dist(X_{g_1(x)}(x), \gamma^s_\sigma \cap \gamma_1)
  - \log\dist(X_{g_1(y)}(y), \gamma^s_\sigma \cap\gamma_1)\Big|
  \\
  & \le
  \displaystyle\frac{C_0|X_{g_1(x)}(x)- X_{g_1(y)}(y)|}
  {\min\{\dist(X_{g_1(x)}(x),\gamma^s_\sigma\cap\gamma_1),
    \dist(X_{g_1(y)}(y), \gamma^s_\sigma \cap \gamma_1)\}}
  \end{align*}
  for some constant $C_0>0$.  Because $g_1$ is smooth and bounded, we
  have that the Poincar\'e map between $\gamma$ and $\gamma_1$
  distorts distances by at most by a constant factor. Hence, we obtain
  $|X_{g_1(x)}(x) - X_{g_1(y)}(y)| \le C_0 |x - X_{g_1(x) -
    g_1(y)}(y)| \le C_0 |x - y|$ and
  $C_0 \dist(X_{g_1(z)}(z), \gamma^s_\sigma \cap \gamma_1) \ge
  \dist(z, \pi_\gamma(\Gamma))$. It follows that
  \begin{align}
    \label{eq:poincare_time_aux_1}
    |\tau_1(X_{g_1(x)}(x)) - \tau_1(X_{g_1(y)}(y))|
    \le C_0  \dfrac{|x - y|}{\min\{\dist(x,\pi_\gamma(\Gamma)),
      \dist(y, \pi_\gamma(\Gamma))\}}.
  \end{align}
  Finally, using the expression of $P_1$ it follows that
  $g_2 \circ P_1$ is $\alpha$-H\"older on $\gamma_1$, and so
  $|g_2(X_{\tau_1(x)}(X_{g_1(x)}(x))) -
  g_2(X_{\tau_1(y)}(X_{g_1(y)}(y)))|$ can be written
  \begin{align}
    \label{eq:poincare_time_aux_2}
    |g_2(P_1(X_{g_1(x)}(x))) - g_2(P_1(X_{g_1(y)}(y)))|
    \le C_1 |X_{g_1(x)}(x) - X_{g_1(y)}(y)|^\alpha
    \le C_2 |x - y|^\alpha
  \end{align}
  for come constants $C_1,C_2>0$.  Using the inequalities for $g_1$
  together with~\eqref{eq:poincare_time_aux_1}
  and~\eqref{eq:poincare_time_aux_2}, from \eqref{eq:g_12tau} we
  arrive at the first inequality of the statement in the singular case
  as well.  This completes the proof of the first inequality in the
  statement.

  Now for the proof of the second inequality: in the case that
  $X_{[0, \tau(x)]}(x)$ never enters $V_0$, the result follows because
  $\tau|_{\Sigma_1}$ is of class $C^2$.  Otherwise, 
  $X_{[0,\tau(x)]}$ enters a neighborhood $V_\sigma$ of some
  singularity $\sigma \in S(\Lambda, X)$ and from~\eqref{eq:g_12tau}
  we write $D\tau(x)$ as
  \begin{align*}
    Dg_1(x) + D\tau_1(X_{g_1(x)}(x))DX_{g_1(x)}(x)
    +  Dg_2(P_1(X_{g_1(x)}(x)))DP_1(X_{g_1(x)}(x))DX_{g_1(x)}(x).
  \end{align*}
  Because $g_1$ and $g_2$ are $C^2$ and bounded on the $u$-curves
  where they are defined, we get that there exist a constant $C_3 > 0$
  such that $|Dg_1|_\infty$, $|Dg_2|_\infty \leq C_3$. From the
  expression of $\tau_1$ we get that
  $D\tau_1(X_{g_1}(x)) = \dist(X_{g_1(x)}(x), \gamma^s_\sigma \cap
  \gamma_1)^{-1} \le C_4 \dist(x, \cD)^{-1}$. Identifying $P_1$ with
  $P_1(x_1) = |x_1|^{-\lambda_2 / \lambda_1}$ we have that
  $|DP_1(x_1)| \le C_5|x_1|^{\alpha - 1} \leq |x_1|^{-1}$. Using this
  we get
  $|DP(X_{g_1(x)}(x))| \le C_6 \dist(X_{g_1(x)}(x), \gamma^s_\sigma
  \cap \gamma_1)^{-1} \le C_6 C_5 \dist(x, \pi_\gamma(\Gamma))^{-1}$
  and the result follows.
\end{proof}

\begin{remark}[Horizontal lines are $u$-curves]
\label{subsec:horizontal-line-u-curve}
We choose an identification $L: I \times B^{d_s} \to \Xi$ such that for all
$c \in B^{d_s}$ the curve $\gamma_{i,c}:I_i \to \Sigma_i$ defined by
$\gamma_{i,c}(t) = L(t,c)$ is a $u$-curve in $\Sigma_i$ for all $i = 1, 2,
\dots, m$. In other words, with the identification given by $L$ we may assume
that each horizontal line $\{(t, c):\ t \in I_i\}$ is a $u$-curve for all $c \in
B^{d_s}$.
\end{remark}




\section{Properties of the one-dimensional quotient dynamics}
\label{sec:propert-one-dimens}

We need some specific consequences of the construction and properties
of the one-dimensional quotient map $f$ obtained in
Section~\ref{chapter:decay-correlations-suspension}.

\subsection{Topological properties of the one-dimensional dynamics}
\label{sec:topological-dynamics}

The following provides the existence of a special class of
periodic orbits for $f$.

\begin{proposition}{\cite[Lemma 6.30]{AraPac2010}}
  \label{pr:intervalocresce}
  Let $f:\sum_j I_j\to I$ be a piecewise $C^1$ expanding map with finitely many
  branches $I_1, I_2 \dots, I_m$ such that each $I_j$ is a nonempty open
  interval, $|f'|_{I_j}|\ge\sigma>2$ and $I\setminus(\sum_j I_j)$ is finite.
  Then, for each small $\delta>0$ there exists $n = n(\delta)$ such that, for
  every nonempty open interval $J \subset \sum_j I_j$ with $|J| \geq \delta$, we
  can find $0 \leq k \leq n$, a sub-interval $\widehat J$ of $J$ and $1 \leq j \leq m$
  satisfying
  \begin{align*}
    f^k|_{\widehat J} : \widehat J\to I_j \quad\text{is a diffeomorphism}.
  \end{align*}
  In addition, $f$ admits finitely many periodic orbits
  $\cO(p_1),\dots,\cO(p_k)$ contained in $\sum_j I_j$ with
  the property that every nonempty open interval
  $J\subset\sum_j I_j$ admits an open sub-interval $\widehat J$,
  a periodic point $p_j$ and an iterate $n$ such that
  $f^{n}|_{\widehat J}$ is a diffeomorphism onto a neighborhood
  of $p_j$.
\end{proposition}

\begin{remark}\label{rmk:denstableF}\text{}
  \begin{enumerate}
  \item For the bidimensional map $F$ this shows that there
    are finitely many periodic orbits
    $\cO(P_1), \dots,\cO(P_k)$ for $F$ so that
    $\pi(\cO(P_i))=\cO(p_i), i=1,\dots,k$, where
    $\pi:Q\to I$ is the projection on the first
    coordinate. Moreover, the union of the stable manifolds
    of these periodic orbits is dense in $Q$. See
    \cite[Section 6.2]{AraPac2010} for details.
  \item This also implies that the stable manifolds of the
    periodic orbits $P_i$ obtained above are dense in a
    neighborhood $U_0$ of $\Lambda$.

  \end{enumerate}
\end{remark}

\subsection{Exponential slow recurrence to the critical set}

As a subtle consequence of Theorem~\ref{thm:props-poincare-maps} in
\cite{ArSzTr} it was proved that the quotient map along stable leaves
has \emph{exponentially slow recurrence to the critical/singular set
  $\cD$} as follows.

\begin{lemma}\cite[Theorem C]{ArSzTr}\label{lemma:slow-recurrence}
  For each $\varepsilon>0$ we can find $\delta>0$ and $\xi > 0$ so
  that
  \begin{align}\label{eq.expslowrecurrence}
    \limsup_{n\to\infty}\frac{1}{n}\log\leb
    \left\{ x\in I: \frac1n \sum_{j=0}^{n-1}
    -\log\dist_{\delta}(f^j(x), \cD)>\varepsilon\right\} <-\xi.
  \end{align}  
\end{lemma}

\begin{remark}
  Exponential slow recurrence implies a weaker condition: \emph{the
    slow recurrence to $\cD$}, that is, for all $\varepsilon > 0$
  there exists $\delta > 0$ such that
  \begin{eqnarray}\label{eq:slowrecurrence}
    \limsup_{n \to +\infty}\frac{1}{n}\sum^{n - 1}_{j =
    0}-\log\dist_\delta(f^j(x), \cD) \leq \varepsilon \ \text{for} \ \Leb-a.e.\ 
    x \in I.
  \end{eqnarray}
\end{remark}

\subsection{Ergodic properties of $f$}
\label{sec:ergodic-properties}

The map $f$ is piecewise expanding with H\"older
derivative which enables us to use strong results on
one-dimensional dynamics.

\subsubsection{Existence and finiteness of acim's}
\label{sec:existence-finiten-ac}

It is well-known \cite{HK82} that $C^1$ piecewise expanding maps $f$
of the interval such that $1/|f'|$ has bounded variation have
absolutely continuous invariant probability measures whose basins
cover Lebesgue almost all points of $I$.

Using an extension of the notion of bounded variation this result was
extended in \cite{Ke85} to $C^1$ piecewise expanding maps $f$ such
that $1/|f'|$ is $\alpha$-H\"older for some $\alpha\in(0,1)$.  In
addition, from \cite[Theorem 3.3]{Ke85}, there are finitely many ergodic
absolutely continuous invariant probability measures
$\upsilon_1,\dots,\upsilon_l$ of $f$ and every absolutely continuous
invariant probability measure $\upsilon$ decomposes into a convex
linear combination $\upsilon=\sum_{i=1}^l a_i\upsilon_i$. From
\cite[Theorem 3.2]{Ke85} considering any subinterval $J\subset I$ and
the normalized Lebesgue measure $\leb_J=(\Leb|_{J})/\leb(J)$ on $J$,
every weak$^*$ accumulation point of
$n^{-1}\sum_{j=0}^{n-1} f_*^j(\leb_J)$ is an absolutely continuous
invariant probability measure $\upsilon$ for $f$ (since the
characteristic function of $J$ is of generalized $1/\alpha$-bounded
variation). Hence, the basins of the $\upsilon_1,\dots,\upsilon_l$ cover
$I$ Lebesgue modulo zero:
$ \leb\big(I\setminus(B(\upsilon_1)+\dots+ B(\upsilon_l)\big) =0. $
Note that from \cite[Lemma 1.4]{Ke85} we also know that \emph{the
  density $\varphi$ of any absolutely continuous $f$-invariant
  probability measure} \emph{is bounded from above}.

\subsubsection{Absolutely continuous measures and
  periodic orbits}
\label{sec:absolut-contin-measu}

Now we relate some topological and ergodic properties.

\begin{lemma}\cite[Lemma 2.12]{ArSzTr}\label{le:persuppacim}
  For each periodic orbit $\cO(p_i)$ of $f$ given by
  Proposition~\ref{pr:intervalocresce}, there exists a unique
  ergodic absolutely continuous $f$-invariant probability
  measure $\upsilon_j$ such that $p_i\in\inter(\supp\upsilon_j)$,
  and vice-versa.
\end{lemma}

\subsection{Consequences for the flow dynamics in the trapping region}
\label{sec:conseq-flow-dynamics}

Combining the previous properties we can deduce the following useful result.

\begin{theorem}\cite[Theorem 2.14]{ArSzTr}
  \label{thm:densestablesing}
  The union of the stable manifolds of the singularities in a
  non-trivial connected singular-hyperbolic attracting set $\Lambda$
  is dense in the topological basin of attraction, that is
$
    U_0\subset \overline{\bigcup_{\sigma\in S(X,\Lambda)}W^s(\sigma)}.
$
\end{theorem}

This in particular implies the following.

\begin{proposition}
  \label{pr:nobasicpiece}
  The support of every ergodic physical measure of a non-trivial
  connected singular-hyperbolic attracting set contains some
  Lorenz-like singularity.
\end{proposition}

\begin{proof}
  Arguing by contradiction, let $\mu$ be an ergodic physical measure
  such that $\supp(\mu)$ does not contain any Lorenz-like
  singularity. Hence it does not contain any singularity by
  Remark~\ref{rmk:no-recur-non-Lorenz}(3).

  Therefore, $\Lambda_0=\supp(\mu)$ is a uniformly hyperbolic
  transitive subset (by Theorem~\ref{thm:hyplemma}) and unstable
  manifolds are well-defined and contained in $\Lambda_0$; see
  e.g. \cite{PS82}. Thus $\Lambda_0$ is a connected hyperbolic
  attractor which is a closed and open subset of $\Lambda$; see
  e.g. \cite{PM82}. By connectedness of $\Lambda$ we must have
  $\Lambda_0=\Lambda$. This contradicts the non-trivial assumption on
  $\Lambda$.
\end{proof}

\subsection[Construction of an induced Markov map]{Construction of an
  induced piecewise expanding Markov map for the one-dimensional
  quotient transformation}
\label{section:markov-induced}

Here we explain how to obtain a $C^{1 + \alpha}$ expanding map induced
by $f$ with inducing time having exponential tail, as defined in
Subsection \ref{subsection:uniformly-expanding-maps}, for each
$f$-invariant ergodic absolutely continuous probability measure. A
general reference containing the main results and complete detailed
proofs is~\cite{Alves2020b}.

\subsubsection{Hyperbolic times}
\label{sec:hyperbolic-times}

Let $B > 1$ and $\beta$ be as in the non-degeneracy
conditions (C1) and (C2). Let $0 < \sigma < 1$, $0 < b < 1/2$ and
$\delta > 0$. We say that a natural number $n$ is a
$(\sigma, \delta)$-hyperbolic time for $x \in I$ if for all
$1 \leq k \leq n$, we have
\begin{align*}
  |(f^k)'(f^{n - k}(x))| \geq \sigma^{-k} \ \ \ \text{and}\ \ \
  \dist_\delta{(f^{n - k}(x),\cD)} \geq \sigma^{kb}.
\end{align*}
\begin{lemma}{\cite[Lemma 5.2 \& Corollary 5.3]{ABV00}}
  \label{le:hyp-time1} Given $\sigma\in(0,1)$ and
  $\delta>0$, there exist $\delta_1, D_{1}$ (depending only on
  $\sigma,\delta$ and on the map $f$) such that for any $x\in I$ and
  \( n\geq 1 \) a $(\sigma,\delta)$-hyperbolic time for \( x \), there
  exists an open interval \( V_n(x) \) containing \( x \) with the
  following properties:
  \begin{enumerate}[(a)]
  \item $f^{n}$ maps $V_n(x)$ diffeomorphically onto the interval
    $(f^n(x) - \delta_1, f^n(x) + \delta_1)$;
  \item for $1\le k <n$ and $y, z\in V_n(x)$, $
    |f^{n-k}(y)-f^{n-k}(z)| \le \sigma^{k / 2}|f^{n}(y)-f^{n}(z)|$;
  \item $f^n$ has distortion bounded by $ D_1$ on $V_n(x)$, that is,
    $|(f^n)'(y)| / |(f^n)'(z)| \leq D_1$, for all $y, z \in V_n(x)$;
  \item 
    \(V_n(x) \subset (x - 2\delta_1\sigma^n, x + 2\delta_1 \sigma^n) \).
  \end{enumerate}
\end{lemma}


The sets $V_n(x)$ in the last lemma are called \emph{hyperbolic pre-intervals}
and their images $f^n(V_n(x))$ are called \emph{hyperbolic intervals}.

\begin{lemma}{(\cite[Lemma 5.4]{ABV00} and \cite[Lemma
    1.8]{APP14})}\label{le:hyp-time2}
  There exist \( \theta>0 \) and $0<\delta<1$ (depending only on $f$
  and on its expanding rate) such that, for Lebesgue almost every
  \( x \in I \), we can find $n_0\ge1$ satisfying: for each $n>n_0$
  there are $(\sigma,\delta)$-hyperbolic times
  $1 \le n_1 < \cdots < n_l \le n$ for \( x \) with $l\ge\theta
  n$. Moreover, each hyperbolic time $n_i$ satisfies
  \begin{eqnarray}\label{eq:prop-app14}
    \sum_{j = n_i - k}^{n_i - 1}\log\dist_\delta(f^j(x), \cD) \geq b k \log
    \sigma,\ \text{for all } 0 \leq k \leq n_i, \ 1 \leq i \leq l.
  \end{eqnarray}
\end{lemma}

\subsubsection{Inducing the one-dimensional map}
\label{sec:inducing-one-dimens}

We have all the conditions to perform the construction of an induced
map from $f$ as in \cite[Main Theorem
1]{alves-luzzatto-pinheiro2005}.

Let $\delta_1$ be given by Lemma \ref{le:hyp-time1}. \emph{For each
  $f$-invariant ergodic absolutely continuous probability measure
  $\nu$,} we fix a point $p \in \inter(\supp\nu)$ and an integer
$N \geq 1$ such that $\cup^N_{j = 0}f^{-j}(p)$ is $\delta_1/3$-dense
in $\supp\nu$ and does not contain any element of $\cD$.

\begin{theorem}\cite{alves-luzzatto-pinheiro2004,ArVar, gouezel}\label{thm:markov}
  There exists a neighborhood
  $\Delta \subset \inter(\supp\nu) \setminus \cD$, a countable Lebesgue
  modulo zero partition $\cP$ of \( \Delta \) into sub-intervals; a
  function \( R: \Delta \to \NN \) defined almost everywhere, constant
  on elements of the partition \( \cP \); and constants \( C >0 \),
  $0 < \rho < 1$ such that, for all \( J\in\cP \) and \( R=R(J) \),
  the map \( F:=f^{R}:J \to \Delta \) is a \( C^{1 + \alpha} \)
  diffeomorphism, satisfies the bounded distortion property and is
  uniformly expanding: for each \( x,y\in J \)
  \begin{align*}
    \left| \frac{F'(x)}{F'(y)} -1 \right| \le C
    |F(x)-F(y)|^\alpha
    \quad\text{and}\quad
    |F(x)-F(y)|>\rho^{-1}|x-y|.
  \end{align*}
  Moreover, for each $J\in\cP$ there exists $0<k\le N$
  such that $n:=R(J)-k$ is a
  $(\sigma,\delta_1)$-hyperbolic time for each
  $x\in J$; $J\subset V_{n}(x)$ and, in addition, there is some $\delta
  > 0$ such that
  $\dist(f^j(J), \cD) \geq \delta$ for all $n\le j <
  R(J)$.
\end{theorem}

\begin{remark}\label{remark:R-hyperbolic-time}
  Without loss of generality we assume that $\rho < \sigma$ and use
  $\sigma^{-1}$ as the expansion rate for $F$.  Making $\delta_!$
  smaller if necessary, we can take $\delta=\delta_1$ and we may
  assume that $R(J)$ is a $(\sigma, \delta_1)$-hyperbolic time for all
  $x \in J$ (instead of $R(J) - k$), because the map $f$ is expanding
  and the iterates $f^j(J)$ are $\delta$-distant from $\cD$. From now
  on we make this assumption.
\end{remark}

Let $\cP^{(n)} = \bigvee^{n - 1}_{j = 0}F^{-j}(\cP)$ and denote by
$\cP^{(n)}(x)$ the element of the partition $\cP^{(n)}$ that contains $x$. Letting
$R_n(x) = \sum_{j = 0}^{n - 1}R(F^j(x))$ we get that $F^n(x) = f^{R_n(x)}(x)$. It follows from Theorem
\ref{thm:markov} that, for all $x, y \in \Delta$ such that $y \in \cP^{(n)}(x)$
\begin{eqnarray}\label{eq:iterate-markov}
  |F^n(x) - F^n(y)| \geq \sigma^{-(R_n(x) - j)}|f^j(x) - f^j(y)|,\ j = 0, 1, \dots,
  R_n(x) - 1,
\end{eqnarray}
and since $R_n(x)$ is a $(\sigma, \delta_1)$-hyperbolic time for every
$y \in \cP^{(n)}(x)$
\begin{eqnarray}\label{eq:slow-recurrence-hyperbolic-time}
  \dist(f^j(y), \cD) \geq \sigma^{b(R_n(x) - j)},\ j = 0, 1, \dots, R_n(x) - 1.
\end{eqnarray}

\subsection{The $C^1$ expanding semiflow}
\label{sec:c1-expand-semifl}

Now we check the conditions (3) and (4) from
Subsection~\ref{section:hyperbolic-skew-product-semiflow} to obtain an
expanding semiflow associated to each ergodic absolutely continuous
$f$-invariant probability measure.
  
\subsubsection{The induced roof function}

Let $r: \Delta \to \RR$ be defined as
$r(x) = \sum^{R(x) - 1}_{j = 0} \tau(f^j(x)),$ for all $x \in
\Delta$. Next we prove that $r$ satisfies condition (3) of
Subsection~\ref{subsection:uniformly-expanding-maps}.
\begin{lemma}
  \label{lemma:prop-iii-armel} For all $h \in \cH$ it holds that
  $|(r \circ h)'|_\infty < + \infty.$
\end{lemma}

\begin{proof}
  Let $h \in \cH$, $h: \Delta \to J$ be an inverse branch
  for $F = f^R$ with inducing time $\ell = R(J)$. Fixing  $x \in
  J$, we have
  \begin{align*}
  |(r \circ h)'(x)| = \frac{|r'(h(x))|}{|F'(h(x))|} =
  \left|\displaystyle\sum^{\ell - 1}_{j = 0}\frac{(\tau' \circ f^j)
      \cdot (f^j)'}{F'} \circ h(x)\right|.
  \end{align*}
  By Lemma~\ref{lemma:property-poincare-time} and the fact that $\ell$
  is a hyperbolic time for all $x \in J$ we have that
  $|\tau'(f^j(h(x)))| \leq C\dist(f^j(h(x)), \cD)^{-1} \leq
  C\sigma^{-(\ell - j) b}$, for some constant $C > 0$ and for all
  $j = 0, 1, \dots, \ell - 1.$ It follows from inequality
  \eqref{eq:iterate-markov} that
  $|(f^j)'/F'| \circ h(x) \leq \sigma^{\ell - j}$. Hence
  $|(r \circ h)'(x)| \leq C\sum^{\ell - 1}_{j = 0}\sigma^{-(\ell -
    j)b} \cdot \sigma^{\ell - j} \leq C \sum_{j =
    0}^{\infty}\sigma^{(1 - b)j} < \infty.  $
\end{proof}

For each $x \in \Delta$ and $\varepsilon, \delta > 0$ we define the
\emph{recurrence time} of $x$ by
\begin{align*}
\cR(x) = \cR_{\varepsilon, \delta}(x) = \min\left\{N \geq 1:\,
  \dfrac{1}{N}\sum^{N - 1}_{j = 0} - \log\dist_\delta(f^j(x), \cD) \leq 2
  \varepsilon\right\}.
\end{align*}
The slow recurrence given by the Lemma \ref{lemma:slow-recurrence} can
be translated as: for all $\varepsilon > 0$, there exists $\delta > 0$
such that
$\limsup_{n \to +\infty}\frac{1}{n}\log \Leb(\{x \in I:\, \cR(x) >
n\}) < 0.$ We say that the sets $\{x \in I:\, \cR(x) > n\}$ are the
\emph{tail of hyperbolic times}. Hence, the tail of hyperbolic times
converges exponentially fast to zero. In~\cite{gouezel} the
construction of the induced Markov map
from~\cite{alves-luzzatto-pinheiro2004} was improved so that the tail
of $R$ converges to zero at the same speed as the tail of hyperbolic
times. In particular, in the same setting of Theorem~\ref{thm:markov},
the inducing time function \( R \) has exponential tail.
\begin{proposition}\label{thm:exp-tail} There
  exist constants $c, C > 0$ such that for all $n \geq 1$
  $$
  \leb(\{x\in \Delta : R(x)>n\}) \leq C e^{-c n}.
  $$ 
\end{proposition}

\begin{remark}\label{rmk:eslami}
  Recently a similar result has been stated in \cite{Eslami20} under
  weaker assumptions.
\end{remark}

It follows from Proposition \ref{thm:exp-tail} that

\begin{proposition}\label{pr:rexptail}
  The function $r$ has exponential tail.
\end{proposition}

\begin{proof}
  Fixed $J \in \cP$, using inequality \eqref{eq:prop-app14} of
  Lemma \ref{le:hyp-time2} and the fact that $\tau$ has logarithmic growth
  (item \eqref{subitem:logarithmic-growth} of
  Theorem \ref{thm:props-poincare-maps}) we get

  $$r(x) = \sum_{j = 0}^{R(x) - 1}\tau(f^j(x)) \leq \sum_{j = 0}^{R(x) - 1}-\log
  \dist_\delta(f^j(x), \cD) \leq - bR(x)\log \sigma .$$ Thus, there
  exists constants $C, \tilde c > 0$ such that
  $\Leb(\{x \in \Delta:\ r(x) > n\}) \le C e^{-\tilde cn}$. Using
  the bounded distortion of $F'$ is straightforward to get that $r$
  satisfies the exponential tail condition as stated in
  Subsection~\ref{subsec:expanding-semiflows}.
\end{proof}

\subsection{$C^{1+}$ skew product semiflow}
\label{sec:c1+-skew-product}

Now we note that from Theorem~\ref{thm:props-poincare-maps} we already
have all that is needed to obtain a $C^{1+}$ skew product semiflow, with
the exception of the UNI condition (5), which we focus on the
Subsection \ref{sec:uni-condition}.

Indeed, let $\widehat{\Delta} = \cup_{x \in \Delta}W^s_x$ and define
$\widehat{F}: \widehat{\Delta} \to \widehat{\Delta}$ by
$\widehat{F}(x) = P^{R(x)}(x)$ and the suspension semiflow
$\widehat{F}_t: \widehat{\Delta}^r \to \widehat{\Delta}^r$ with base
map $\widehat F$ and roof function $r$.

Through the identification of $\widehat{\Delta}$ with
$\Delta \times B^{d_s}$, it follows from item 4 of
Theorem~\ref{thm:props-poincare-maps} that $\widehat{F}_t$ satisfies
condition (4) of a hyperbolic skew product semiflow, from
Subsection~\ref{section:hyperbolic-skew-product-semiflow}.


\section{Exponential mixing for singular-hyperbolic attracting sets}
\label{cha:exponent-mixing-sing}

Throughout this section we denote by $\fX^s_{dsh}(M)$ the subset of
$\fX^s(M)$ that admits a $2$-strong dissipative singular-hyperbolic
attracting set, $s\ge1$.

Here we construct a $C^{2}$ open subset
$\cU \subset \fX^2_{dsh}(M)$ where the UNI condition holds. This
enable us to construct a suspension semiflow with exponential decay of
correlations as in Theorem
\ref{prop:decay-correlations-skew-product-weaker}.

We also prove Theorem \ref{theorem:conjugacy-skew-product-hyperbolic}
showing that smooth observables for the original flow lie on the right
function spaces when composed with the conjugacy.

Finally, we finish the section by deducing the exponential convergence
to the equilibrium for the original flow. This is a by product of all
the work made to prove exponential decay of correlations for the
physical measures.

\subsection{The UNI condition after small perturbations}
\label{sec:uni-condition}

In this section we construct an open and dense subset $\cU$ of
$\fX^2_{dsh}(M)$ where all vector fields have a roof
function that satisfies the UNI condition. In particular, we construct
a family of suspension semiflows, one for each ergodic physical
measure of the attracting set, which satisfy the conditions of Theorem
\ref{prop:decay-correlations-skew-product-weaker} and we get
exponentially mixing for them.

Recall that there exists a one-to-one correspondence between the
periodic points of the Poincaré map and its quotient along the stable
leaves.

\begin{lemma}
  \label{lemma:relation-periodic-orbit-dim1-dim2}
  A point $x \in I$ is periodic for $f$ if and only if there exists a
  periodic point $z \in \pi^{-1}(x)$ for the Poincaré map $P$.
\end{lemma}


\begin{remark}
  \label{remark:convention-periodic-points}
  Using Lemma \ref{lemma:relation-periodic-orbit-dim1-dim2}, from now
  on we make no distinction between a periodic point of the Poincar\'e
  map and its quotient along stable leaves.
\end{remark}

Since the strong dissipative condition is open in the $C^1$ topology
and singular-hyperbolic attracting sets persist by $C^1$-small
perturbations of the vector field, we have that $\fX^s_{dsh}(M)$ is
open in $\fX^s(M)$ (with the $C^s$ topology), for all $s\ge1$.

For a vector field $X \in \fX^2_{dsh}(M)$ we can repeat the
constructions of Chapter
\ref{chapter:decay-correlations-suspension}. We need to perform the
constructions for more than one vector field so, where necessary, we
make the dependence on the vector field explicit in what follows. For
instance, $P_X: \Xi \backslash \Gamma_X \to \Xi_X$ denotes the
Poincaré map of $X$ with Poincaré time given by $\tau_X$; and
$f_X:I\setminus\cD\to I$ is the corresponding one-dimensional quotient
map.

\begin{definition}\label{def:uni-condition-after}
  Let $\Lambda$ be a singular-hyperbolic attracting set for
  $X \in \fX^2_{dsh}(M)$.  Let $P: \Xi \backslash \Gamma \to \Xi$ and
  $f:I \backslash \cD\to I$ be the global Poincaré map and its
  quotient along stable leaves, respectively, for $X$ as in Subsection
  \ref{sec:global-poincare-retu}.

  We say that $X$ satisfies the UNI condition if, for each ergodic
  physical measure $\mu$ of $X$ corresponding to an ergodic physical
  measure $\hat\nu$ of $P$ given by an ergodic $f$-invariant
  absolutely continuous probability measure $\nu$, there exists an
  open interval $\Delta \subset \inter(\supp\nu)\subset I$ and an
  induced function $R: \Delta \to \NN$ (as in Theorem
  \ref{thm:markov}) such that the induced roof function
  $r: \Delta \to \RR$ given by
  $r(x) = \sum^{R(x) - 1}_{j = 0}\tau(f^j(x))$ satisfies the UNI
  condition.
\end{definition}

Let us fix an ergodic physical measure $\hat \nu$ for $P$ and $\nu$
for $f$. Letting $F=f^R:\Delta\to\Delta$ be an induced
full branch Markov map constructed for $X \in \fX^2(M)$ in Theorem
\ref{thm:markov}, for a function $\varphi: \Delta \to \RR$ and
$n \geq 1$ we denote
$S^F_n\varphi = \sum^{n - 1}_{j = 0} \varphi \circ F^j$.

Now we describe the open set $\cU$ where Theorem
\ref{theorem:conjugacy-skew-product-hyperbolic} holds. We set $\cU$ to
be the subset of vector fields in $\fX^2_{dsh}(M)$ such that, for all
$X \in \cU$, each physical measure $\nu$ of $f$ and each
corresponding induced Markov map $F:\Delta\to\Delta$, there exist two
distinct periodic points $x_1, x_2 \in \Delta\subset\inter(\supp\nu)$
for the induced Markov map $F$ with the same period $p$ and
satisfying:
\begin{enumerate}[(i)]
\item the orbits are distinct; and they visit the interior of the same
  elements of the partition $\cP$ the same number of times as the
  other, but necessarily in some different order to each other; and
\item $S^{F}_pr(x_1) \neq S^{F}_pr(x_2)$.
\end{enumerate}

\begin{lemma}\label{lemma:U-UNI}
  The vector fields in $\cU$ satisfy the UNI condition.
\end{lemma}

\begin{proof}
  Let $X \in \cU$ and assume that $X$ does not satisfies the UNI
  condition.  Then, there exist an ergodic $f$-invariant absolutely
  continuous probability measure $\nu$, an induced map
  $F:\Delta\to\Delta$ with $\Delta\subset\inter(\supp\nu)$, a $C^1$
  function $\varphi:\Delta \to \RR$ and a function
  $\psi: \Delta \to \RR$ constant on elements of the induced partition
  $\cP$ of $\Delta$ such that $r = \psi + \varphi \circ F - \varphi$.

  Let $x_1, x_2 \in \Delta$ be two periodic points with period $p$ for
  the induced Markov map $F$ satisfying conditions (i) - (ii) of the
  definition of $\cU$. It follows that $S^F_pr(x_i) = S^F_p\psi(x_i)$.
  Since $\psi$ is constant on elements of the partition $\cP$, by
  condition (ii) above, we get that $S^F_p\psi(x_1) = S^F_p\psi(x_2)$.
  This is a contradiction with $S^F_pr(x_1) \neq S^F_pr(x_2)$.
\end{proof}

The proofs of the next two propositions follow the steps presented
in~\cite{ArVar-errata}.  In Proposition~\ref{prop:uni-is-dense} we
show that, if we start with a vector field $X \in \fX^2_{dsh}(M)$ that
does not satisfy the UNI condition and change slightly the velocity of
a well chosen periodic orbit, then the new vector field satisfies the
UNI condition and is arbitrarily $C^2$-close to the initial vector
field $X$.  In particular, we get that the subset of vector fields in
$\fX^2_{dsh}(M)$ that satisfies the UNI condition is dense in the
$C^2$ topology.  In Proposition \ref{prop:uni-open}, we show that the
inequality that we obtained in the previous propostion remains valid
for vector fields $C^{2}$-close to $X \in \cU$.

\begin{proposition}\label{prop:uni-is-dense}
  The set $\cU$ is $C^2$-dense in $\fX^2_{dsh}(M)$:
  for each $X \in \fX^2_{dsh}(M)$ there exists $\delta > 0$
  and a $\delta$-$C^2$-close vector field $Y \in \cU$ which is a
  multiple of $X$. 
\end{proposition}

In other words, any $X\in \fX^2_{dsh}(M)$ admits a time
reparametrization which lies in $\cU$.

\begin{proof}
  Let $X \in \fX^2_{dsh}(M)$ and assume that $X$ does not satisfy the
  UNI condition. Hence, there exist an ergodic $f_X$-invariant
  absolutely continuous probability measure $\nu$, an induced map
  $F_X:\Delta_X\to\Delta_X$ with $\Delta_X\subset\inter(\supp\nu)$, a $C^1$
  function $\varphi:\Delta \to \RR$, and a function
  $\psi: \Delta_X \to \RR$ constant on elements of the partition
  $\cP_X$ together with $\varphi: \Delta_X \to \RR$ of class $C^1$ such that
  $r_X = \psi + \varphi \circ F_X - \varphi$.

  Let $x_1$ and $x_2$ be two periodic points with the same period $p$
  for the map $F_X$. We may assume without loss of generality, because
  $F_X$ is a full branch Gibbs-Markov map, that the orbits are
  distinct and visit the same elements of the partition $\cP_X$ the
  same number of times as the other, but in a different order.  If
  $J_1,J_2$ are two disjoint elements of the partition $\cP_X$, we can
  choose the period $p=4$, $x_1$ and $x_2$ such that
  \begin{align*} x_1, F_X(x_1) \in J_1, F_X^2(x_1), F_X^3(x_1) \in J_2
    \quad\text{and}\quad x_2, F_X^2(x_2) \in J_1, F_X(x_2), F_X^3(x_2) \in J_2.
  \end{align*} Note that
  $
    S^{F_X}_pr(x_1)  =  S^{F_X}_p\psi(x_1) = S^{F_X}_p\psi(x_2) =
    S^{F_X}_pr(x_2)
  $ since $x_1$ and $x_2$ visit the same elements of $\cP_X$
  an equal number of times and $\psi$ is constant on each element of
  $\cP_X$. We have that $S^{F_X}_pr(x_i)$ is a multiple of the period
  of $x_i$ with respect to the flow of $X$, and does not depends on
  the functions $\varphi$ and $\psi$.  Indeed
  \begin{equation}\label{eq:period-for-flow}
    S^{F_X}_4r(x_i)
    = \sum^{3}_{j = 0}r(F_X^j(x_i))
    = \sum^{3}_{j = 0}
    \sum^{R(F_X^j(x_i)) - 1}_{k =0}
    \tau(f_X^{k}(F_X^j(x_i)))
    = \sum^{R_4(x_i) - 1}_{k = 0} \tau(f^k(x_i))
  \end{equation}
  Because $F_X^4(x_i) = x_i$, it follows that
  $f_X^{R_4(x_i)}(x_i) = x_i$.  Letting $t_0 = S^{F_X}_4r(x_i)$ and
  using~\eqref{eq:period-for-flow}, we have $X_{t_0}(x_i) = x_i$
  (recall the convention that we are using for periodic points on
  Remark \ref{remark:convention-periodic-points}). Hence, it follows
  that $t_0$ is a multiple of the period of $x_1$ and $x_2$ by the
  action of the flow.

  We modify the roof function $\tau_X:\Xi\to\RR^+$ in a small
  neighborhood of $x_2$ that does not intersect the orbit of $x_1$ to
  ensure that the induced roof function $r_X$ satisfies
  $S^{F_X}_{p}r(x_1) > S^{F_X}_{p}r(x_2)$.  Let $V_0$ and $V_1$ be
  open small neighborhoods of $x_2$ that do not intersect the orbit of
  $x_1$ with $\overline{V_0} \subsetneq V_1$ and consider a $C^\infty$
  bump function $\xi: M \to [0,1]$ such that $\xi|_{V_0} \equiv 1$ and
  $\xi_{M \backslash V_1} \equiv 0$. For all $\delta > 0$ define the
  vector field $X_\delta = X + \delta \xi X$ which is $\delta$-close
  to $X$ in the $C^1$-topology. Inside $V_0$ the vector field
  $X_\delta$ is equal to $(1+ \delta)X$ and outside of $V_1$ it equals
  $X$. Thus $x_2$ is still a periodic orbit for $X_\delta$ but with a
  smaller period than before. Thus, he have that
  $S^{F_{X_\delta}}_4r(x_1) > S^{F_{X_\delta}}_4r(x_2)$ as we desired.

  Now it follows from Lemma \ref{lemma:U-UNI} that $X_\delta$
  satisfies the UNI condition. Moreover
  $\|X-X_\delta\|_2=\delta\|\xi X\|_2$ and so we can make the
  perturbation arbitrarily close to the original vector field in $C^2$
  topology. 
\end{proof}

\begin{remark}
  \label{rmk:Cs}
  If we start with a $C^s$ vector field, for some $s\ge2$, then the
  same argument gives a $\delta$-$C^s$-close vector field
  $X_\delta\in\cU$.
\end{remark}

\begin{proposition}\label{prop:uni-open}
  The set $\cU$ is $C^{2}$-open in $\fX^2_{dsh}(M)$.
\end{proposition}

\begin{proof} Let $X \in \cU$ and for a fixed ergodic physical measure,
  let $x_1, x_2 \in \Delta_X$ as conditions (i)-(ii) of the
  definition of $\cU$ above. We are going to show that these
  conditions persist for all $C^2$-close enough vector fields
  $Y$. Since we have only finitely many ergodic physical measures
  supported on the attracting set $\Lambda$ of $X$, it is enough to
  argue for one such ergodic physical measure.
  
  Let $L \geq 1$ be big enough so that
  $\cO_{F_X}(x_i) \subset \{x \in M:\ R_X(x) \leq L\}$, for
  $i = 1, 2$. Recall that because $x_1, x_2$ are also periodic points
  for $X$ inside a singular-hyperbolic attracting set, then they are
  hyperbolic periodic orbits of saddle type and admit smooth
  continuations to all $C^{1}$ nearby vector fields $Y$.

  Because the construction of the induced Markov map from the
  one-dimensional map is made outside a neighborhood of the
  critical/singular set, we can also control the distance of
  $F_X|_{\{R_X \leq L\}}$ and $F_Y|_{\{R_Y \leq L\}}$ (since we are
  only working with finitely many iterates of the maps $f_X$ and
  $f_Y$). Moreover, because the construction is inductive, in each
  step of the construction we can ensure that the open intervals of
  the partition $\cP_X$ inside $\{R_X \leq L\}$ are arbitrarily close
  to their correspondent open intervals of $\cP_Y$ inside
  $\{R_Y \leq L\}$. (To keep all the ingredients of the inductive
  construction preserved by $Y$ here, we need $Y$ and $X$ to be
  $C^{2}$-close. For example, this is needed to control the
  size of the hyperbolic balls. Check the outline of the construction
  in Section~\ref{section:markov-induced} and for more details check
  \cite[Sections 3 and 4]{alves-luzzatto-pinheiro2004} and
  \cite{ArVar-errata}).

  Thus, we have that the continuation $y_i$ of $x_i$, $i = 1, 2$, for
  vector fields $Y$ close to $X$ has the same combinatorics as before,
  that is, $y_1$ and $y_2$ has the same period $p$ and visit the same
  elements of the partition $\mathcal{P}_Y$ with respect to the map
  $F_Y$. Finally, because $S^{F_Y}_pr_Y(y_1)$ is a multiple of the
  period of $y_1$ for $Y$, we have that
  $S^{F_Y}_pr_Y(y_1) \neq S^{F_Y}_pr_Y(y_2)$. Now it follows that
  $r_Y$ cannot be written as
  $r_Y = \psi_Y + \varphi_Y \circ F_Y - \varphi_Y$ with
  $\psi_Y: \Delta_Y \to \RR$ constant on elements of the partition
  $\mathcal{P}_Y$ and $\varphi_Y:\Delta_Y \to \RR$ with class $C^1$,
  otherwise following the same argument of the Proposition
  \ref{prop:uni-is-dense} we would get that
  $S^{F_Y}_pr_Y(y_1) = S^{F_Y}_pr_Y(y_2)$.
\end{proof}

\subsection{Proof of the main technical result }
\label{sec:proof-theorem-refthe}

Here we prove Theorem \ref{theorem:conjugacy-skew-product-hyperbolic}.

We show that the original flow is semiconjugated to a suspension
semiflow $\widehat{F}_t$ and that, given observables with certain
amount of regularity for the original flow, we get observables in the
right space for the suspension semiflow, and the measure in the
original flow is the pushforward of the measure for the suspension
semiflow. This provides what is needed to transfer the results about
decay of correlations from the suspension semiflow to the original
flow.

Let $X \in \cU$ and for each ergodic physical measure supported on the
attracting set, let $F: \Delta \to \Delta$ be the induced Markov map
for $X$ with inducing function given by $R$ and roof function given by
$r(x) = \sum^{R(x) - 1}_{j = 0}\tau(P^j(x))$, as before. We consider also
$\widehat{\Delta} = \cup_{x \in \Delta}W^s_x$ and
$\widehat{F}: \widehat{\Delta} \to \widehat{\Delta}$ defined by
$\widehat{F}(x) = P^{R(x)}(x)$ together with the suspension semiflow
$\widehat{F}_t: \widehat{\Delta}^r \to \widehat{\Delta}^r$.

Using the identification of $\widehat{\Delta}$ with
$\Delta \times B^{d_s}$, it follows that $\widehat{F}_t$ and
$\mu^r_{\widehat{F}}$ satisfy Theorem
\ref{prop:decay-correlations-skew-product-weaker}, that is, we have
exponentially fast decay of correlations in the function spaces
$C^\alpha_{\loc}(\widehat{\Delta}^r)$ and
$C^{\alpha, 2}_{\loc}(\widehat{\Delta}^r)$ for the skew product
semiflow associated to each physical measure of the global Poincar\'e
return map.

\subsubsection{From the suspension flow to the original flow}
\label{sec:from-suspension-flow}

The harder part of Theorem
\ref{theorem:conjugacy-skew-product-hyperbolic} is item (ii). We
obtain this using a H\"older bound on the semiconjugation between the
skew product semiflow and the original flow, given by Theorem
\ref{thm:conjugation}.  In the rest of the section we prove this
bound.

From now on we work with a fixed vector field $X \in \cU$ and a fixed
ergodic physical measure and its corresponding skew product semiflow.

The next result enables us to pass from the ambient manifold $\RR^N$ using the
map $p: \widehat{\Delta}^r \to \RR^N$ given by $p(x,y,u) = X_u(x,y)$.

\begin{theorem}\label{thm:conjugation}
  There is a constant $C > 0$ so that for all
  $(x_1,y_1,u_1), (x_2,y_2,u_2) \in \widehat{\Delta}^r$, we have
  $|p(x_1,y_1,u_1) - p(x_2,y_2,u_2)| \leq C(|F(x_1) - F(x_2)|^\alpha + |y_1 - y_2|+
  |u_1 - u_2|).$
\end{theorem}

\begin{remark}
  \label{rmk:noHolder}
  The map $p$ cannot be H\"older with globally bounded H\"older
  constant since the expansion rate of $F$ is unbounded over all the
  atoms of the induced partition. This detail, which demands the
  extension of the space of admissible observables, was missed in
  previous works on exponential mixing for singular-hyperbolic
  attracting sets.
\end{remark}

\begin{proof}[Proof of Theorem
  \ref{theorem:conjugacy-skew-product-hyperbolic}]
 
  \begin{enumerate}[(i)]
  \item The semiconjugacy property follows directly from the
    definition of $p$ and of the skew product semiflow.

    For the push-forward property, it is enough to check that
    $p_*\mu^r_{\widehat{F}}$ is an ergodic physical measure for $X_t$
    and use the finiteness of such measures and the decomposition of
    any physical measure provided by item 3 of
    Theorem~\ref{theorem:existence-physical-measures}. By construction
    of $\mu_{\widehat{F}}$ we have that this measure is a ergodic
    physical measure for the map $\widehat{F}$ (consult
    \cite[Subsection 6.2]{APPV}). Using the same arguments of
    \cite[Subsection 6.4]{APPV} it follows that $\mu^r_{\widehat{F}}$
    is a physical measure for the suspension flow $\widehat{F}$. Using
    the fact that $p$ is a semiconjugacy between $\widehat{F}_t$ and
    $X_t$, it follows that
    $p(B(\mu^r_{\widehat{F}})) \subset
    B(p_*\mu^r_{\widehat{F}})$. Since $p$ is a local diffeomorphism it
    follows that $B(p_*\mu^r_{\widehat{F}})$ has positive Lebesgue
    measure. Hence $p_*\mu^r_{\widehat{F}}$ is an ergodic invariant
    physical probability measure for $X$ supported in $U_0$, and so in
    $\Lambda$. It follows that $p_*\mu^r_{\widehat{F}}=\mu_i$ for some
    ergodic physical measure supported on $\Lambda$.
  \item The result follows from Theorem~\ref{thm:conjugation} and successive use
    of the Mean Value Inequality Theorem.
    \end{enumerate}
\end{proof}

We are left to prove Theorem~\ref{thm:conjugation} in what follows.

\subsubsection{Proof of the H\"older estimate}
\label{sec:proof-conjug-result}

In the next lemma, given the Poincar\'e time function $\tau$ for $X$,
we denote $\tau_k(x) = \sum^{k - 1}_{j = 0}\tau(P^j(x))$, for all
$x \in \Xi \backslash \Gamma$. Also recall that $\tau$ is constant on
each stable leaf, so $\tau(x) = \tau(\pi(x))$ for all
$x \in \Xi \backslash \Gamma$.

\begin{lemma}\label{lemma:conjug-aux1}
  There exists $C > 0$ such that for all $x_1, x_2 \in \Delta$ with $x_2 \in
  \cP(x_1)$ and $0 \leq k \leq R(x_1) = R(x_2)$ we have that $|\tau_k(x_1) -
  \tau_k(x_2)| \leq C |F(x_1) - F(x_2)|^\alpha$.
\end{lemma}

\begin{proof}
  Using Lemma \ref{lemma:property-poincare-time} and inequalities
  \eqref{eq:iterate-markov} and
  \eqref{eq:slow-recurrence-hyperbolic-time} we get a constant $K>0$
  so that
  $ |\tau_k(x_1) - \tau_k(x_2)|\leq \sum_{j = 0}^{k -
    1}|\tau(f^j(x_1)) - \tau(f^j(x_2)| $ is bounded by
  \begin{align*}
      K&\sum_{j = 0}^{k - 1}\left[\dfrac{|f^j(x_1) -
               f^j(x_2)|}{\min\{\dist(f^j(x_1), \cD), \dist(f^j(x_2),
        \cD)\}} + |f^j(x_1) - f^j(x_2)|^\alpha \right]
    \\
    & \leq
      K \sum_{j = 0}^{k - 1}\left[
      \dfrac{\sigma^{R(x_1) - j}}{\sigma^{b(R(x_1) - j)}}|F(x_1) -
      F(x_2)| + \sigma^{\alpha(R(x_1) - j)}|F(x_1) - F(x_2)|^\alpha\right].
  \end{align*}
    Recalling that $0 < b < 1/2$, we get $\alpha / 2 < 1 - b$ and
    $
      |\tau_k(x_1) - \tau_k(x_2)|$ is bounded above by $ \sum_{j = 0}^{k -
      1}\sigma^{\frac{\alpha}{2}(R(x_1) - j)}|F(x_1) -
      F(x_2)|^\alpha
      \le \textrm{Const.}
      |F(x_1) - F(x_2)|^\alpha.
  $
  \end{proof}

  In the next Lemma we use the partition
  $\widehat{\cP} = \{\cup_{x \in Q}W^s_s: \ Q \in \cP\}$ for
  $\widehat{\Delta}$ which is the same as
  $\widehat{Q} = \{Q \times B^{d_s}:\ Q \in \cP\}$ using the
  identification fixed in
  Section~\ref{chapter:decay-correlations-suspension}.

\begin{proposition}\label{lemma:conjug-aux2}
  There exists a constant $C > 0$ so that for all
  $(x_1, y_1), (x_2, y_2) \in \widehat{\Delta}$ with
  $x_2 \in \cP(x_1)$ and for all $u \in (0, \min\{r(x_1), r(x_2)\})$,
  we have 
  $$|X_u(x_1,y_1) - X_u(x_2,y_2)| \leq C (|F(x_1) - F(x_2)|^\alpha +
  |y_1 - y_2|).$$
\end{proposition}

\begin{proof}
  Let $C > 0$ be as in Lemma \ref{lemma:conjug-aux1} and let
  $\varepsilon > 0$ be such that $C\varepsilon^\alpha < \inf
  \tau$. Fix $Q \in \cP$. First we show the result for all
  $x_1, x_2 \in Q$ such that $|F(x_1) - F(x_2)| < \varepsilon$. There
  exist $\ell, k \in \{1, 2, \dots, R(x_1)\}$ such that
  $u \in [\tau_{\ell - 1}(x_1), \tau_\ell(x_1)]$ and
  $ u \in [\tau_{k - 1}(x_2), \tau_k(x_2)].$

  We get
  $|\tau_j(x_1) - \tau_j(x_2)| \leq C|F(x_1) - F(x_2)|^\alpha \leq
  \inf \tau,$ from Lemma \ref{lemma:conjug-aux1} and the choice of
  $\varepsilon$, for all $j = 0, 1, \dots, R(x_1) = R(x_2)$.  In
  particular, if follows that $|\ell - k| \leq 1$. Without loss we
  assume that $\ell \geq k$.

  Suppose initially that $\ell = k + 1$. In this case we have that
  $\tau_k(x_1) \leq u \leq \tau_k(x_2)$. Recall by Subsection
  \ref{subsec:horizontal-line-u-curve} that each horizontal line in a
  cross-section $\Sigma \subset \Xi$ is a $u$-curve up to identifications. Let
  $\gamma$ be the $u$-curve that contains $(x_1, y_1)$ and $\pi_\gamma$ be the
  projection along stable leaves to $\gamma$. Note that $P^k(x_1,y_1)$ and
  $P^k(x_2, y_2)$ are in the same element of the partition $\widehat{Q}$. In
  particular, there is no singular leaf between $P(x_1, y_1)$ and $P(x_2,
  y_2)$, otherwise $(x_1,y_1)$ and $(x_2, y_2)$ wouldn't be in the same
  element of $\widehat{\cP}$. Let $\eta$ be the strip determined by
  $W^s_{(x_1,y_1)}$ and $W^s_{(x_2,y_2)}$. It follows that $\gamma_k =
  P^k(\gamma \cap \eta)$ is a $u$-curve that crosses the strip $\eta_k$
  determined by $W^s_{P^k(x_1,y_1)}$ and $W^s_{P^k(x_2,y_2)}$ (see Figure
  \ref{fig:poincare-action}).
  \begin{figure}[!htb]
    \centering
    \includegraphics[height=5.5cm]{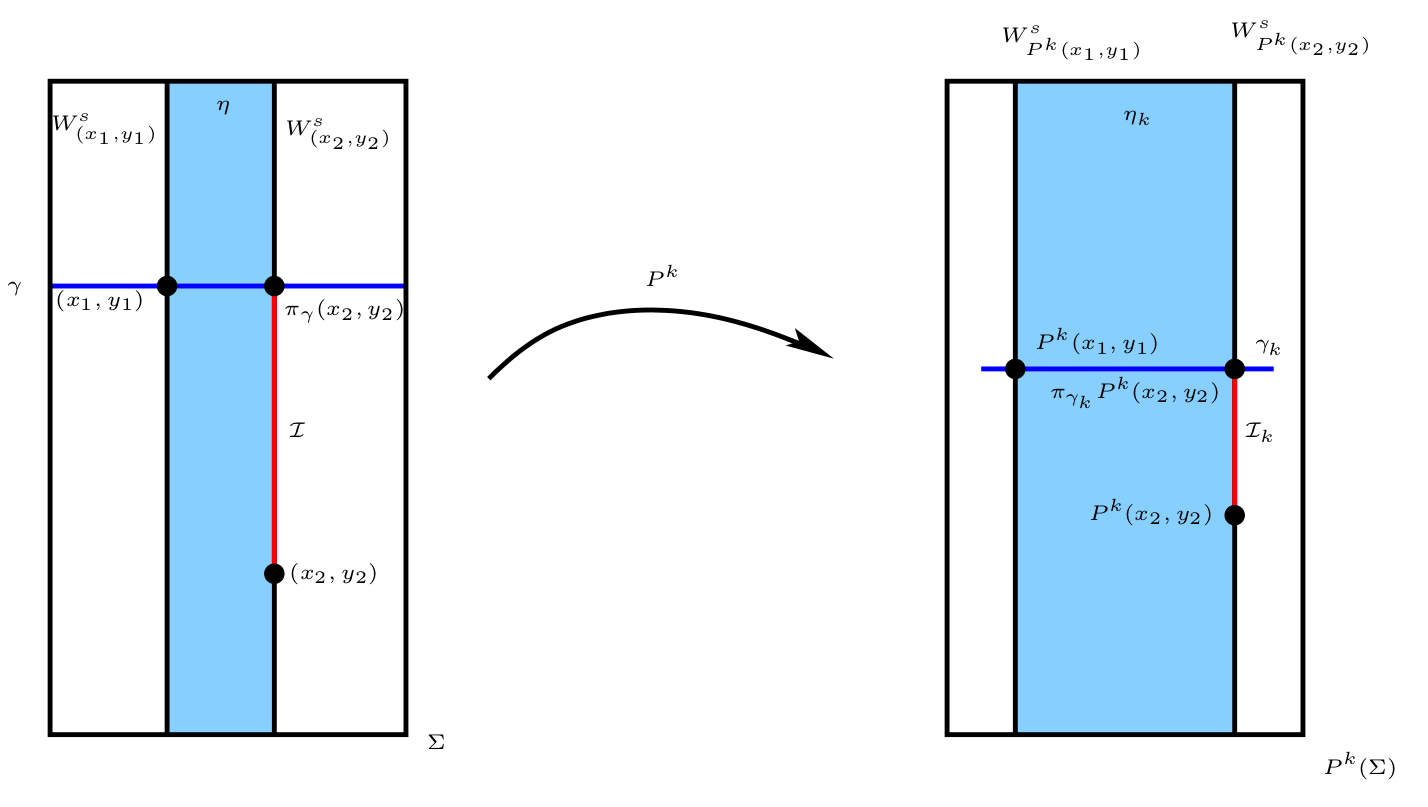}
    \caption{Action of $P^k$.}
    \label{fig:poincare-action}
  \end{figure}
  Hence, $P^k|_\eta$ is a diffeomorphism between the strips $\eta$ and
  $\eta_k$ that maps the interval $\cI$, bounded by $(x_2, y_2)$ and
  $\pi_\gamma(x_2,y_2)$ inside $W^s_{(x_2,y_2)}$, to the interval
  $\cI_k$ bounded by $P^k(x_2, y_2)$ and $\pi_{\gamma_k}P^k(x_2, y_2)$
  inside $W^s_{P^k(x_2,y_2)}$. Thus,
  \begin{align*}
    |X_u(x_1, y_1) - X_u(x_2, y_2)|
    & \leq |X_u(x_1,y_1) - P^k(x_1, y_1)| +
      |P^k(x_1, y_1) - P^k(x_2, y_2)|\\
    & \qquad +  |P^k(x_2, y_2) - X_u(x_2, y_2)|.
  \end{align*}
  We also have by uniform contractions of the stable foliation a
  constant $C>0$ so that
  \begin{align}
    |P^k(x_1, y_2) - P^k(x_2, y_2)|
    & \leq  \diam(\gamma_k\cap \eta_k) +  \diam(\cI_k) \nonumber\\
    & \le
      C |f^k_{\gamma_k}\pi_{\gamma_k}(x_1, y_1) - f^k_{\gamma_k}\pi_{\gamma_k}(x_2, y_2)| + \sigma^k\diam(\cI)\nonumber\\
    & \le C |f^k\pi_{\gamma}(x_1, y_1) - f^k\pi_{\gamma}(x_2, y_2)| + \sigma^k\diam(\cI)\nonumber\\
    & \le  C \sigma^{R(x_1) - k} |F(x_1) - F(x_2)| + \sigma^k|y_1 -
      y_2|,        \label{eq:ineq-aux-0}
  \end{align}
  where we have used \eqref{eq:comp-unidimensional-maps} (to compare
  $f_{\gamma_k}$ with $f$) and \eqref{eq:iterate-markov}.

  We also have that $|X_u(x_1,y_1) - P^k(x_1, y_1)| + |X_u(x_2,y_2) -
  P^k(x_2, y_2)|$ equals
  \begin{align*}
    |X_u(x_1,y_1)
    &-
      X_{\tau_k(x_1)}(x_1, y_1)| + |X_u(x_2,y_2) -
                   X_{\tau_k(x_2)}(x_2, y_2)|
    \\
    &\leq  |X|_\infty(u - \tau_k(x_1)) + |X|_\infty(\tau_k(x_2) - u)
    \\
    &\leq  C|X|_\infty|F(x_1) - F(x_2)|^\alpha,
  \end{align*}
  where in the last inequality we  used  Lemma~\ref{lemma:conjug-aux1}.

  Thus, it follows that there exists a constant $C > 0$ such that 
  \begin{align*}
  |X_u(x_1,y_1) - X_u(x_2, y_2)| \leq C(|F(x_1) - F(x_2)|^\alpha + |y_1 -
  y_2|).
  \end{align*}
  Now consider $\ell = k$. In this case, let $\xi = u - \tau_{k - 1}(x_1)$ and
  $\xi' = u - \tau_{k - 1}(x_2)$ and
  \begin{align}\label{eq:ineq-main}
    |X_u(x_1,y_1) - X_u(x_2, y_2)|
    & \leq
      |X_u(x_1, y_1) -  X_\xi(P^{k-1}(x_2,y_2))| \\
    & \qquad + |X_\xi(P^{k - 1}(x_2, y_2)) - X_u(x_2, y_2)|. \nonumber
  \end{align}
  We also have by Lemma \ref{lemma:conjug-aux1} that
  \begin{align}
    |X_\xi(P^{k - 1}(x_2, y_2)) - X_u(x_2, y_2)|
    & =  |X_\xi(P^{k - 1}(x_2,y_2)) -  X_{\xi'}(P^{k - 1}(x_2,y_2))|\nonumber
    \\
    &\leq \label{eq:ineq-aux-1}
      |X|_\infty|\xi - \xi'|
      =
      |X|_\infty |\tau_{k - 1}(x_1) - \tau_{k - 1}(x_2)|
    \\
    &\leq C|F(x_1) - F(x_2)|^\alpha.
    \nonumber
  \end{align}
  If $P^{k - 1}(x_i,y_i)$, $i = 1, 2$, are in a ingoing cross-section
  for a tubular neighborhood, then there exists a constant $C > 0$
  so that $|X_u(x_1,y_2) - X_\xi(P^{k - 1}(x_2,y_2))|$ equals
  \begin{align*}
      |X_\xi(P^{k - 1}(x_1,y_1)) - X_\xi(P^{k -  1}(x_2,y_2))|
     \leq  C|P^{k - 1}(x_1,y_1) - P^{k -1}(x_2,y_2)|
  \end{align*}
  Analogously to the case $\ell = k + 1$ (see inequality
  \eqref{eq:ineq-aux-0}) we have that
  $|P^{k - 1}(x_1,y_2) - P^{k - 1}(x_2,y_2)| \le C |F(x_1) -
  F(x_2)| + |y_1 - y_2|$ and the result follows.

  Now suppose that $P^{k - 1}(x_i,y_i)$, $i = 1, 2$, are in a ingoing
  cross-section for a flow-box around a singularity. Without loss of
  generality we also assume that
  $\tau(f^{k - 1}(x_1)) \leq \tau(f^{k - 1}(x_2))$. Note that
  $\xi = u - \tau_{k - 1}(x_1) \leq \tau(f^{k - 1}(x_1))$. Let
  $\eta_{k - 1}$ be the strip determined by
  $W^s_{P^{k - 1}(x_1, y_1)}$ and $W^s_{P^{k - 1}(x_2, y_2)}$ (see
  Figure \ref{fig:case-l-equal-k}). Again, there is no singular leaf
  in $\eta_{k - 1}$, otherwise $P^{k - 1}(x_1,y_1)$ and
  $P^{k - 1}(x_2,y_2)$ wouldn't be in the same element of
  $\widehat{\cP}$. Hence $P^{k - 1}(x_1,y_1)$ and $P^{k - 1}(x_2,y_2)$
  will hit the same cross-section in the future.

  Let $\Sigma_{k - 1}, \Sigma_k \in \Xi$ be such that
  $P^j(x_i, y_i) \in \Sigma_j$, for $i \in \{1, 2\}$ and
  $j \in \{k - 1, k\}$. We also define
  $\eta_k = P(\eta_{k - 1})\subset\Sigma_k$ and
  $\widehat{\eta}_k = X_{\tau(x_1)}(\eta_{k - 1})$ (see Figure
  \ref{fig:case-l-equal-k}).
  \begin{figure}[ht]
    \centering
    \includegraphics[height=5cm]{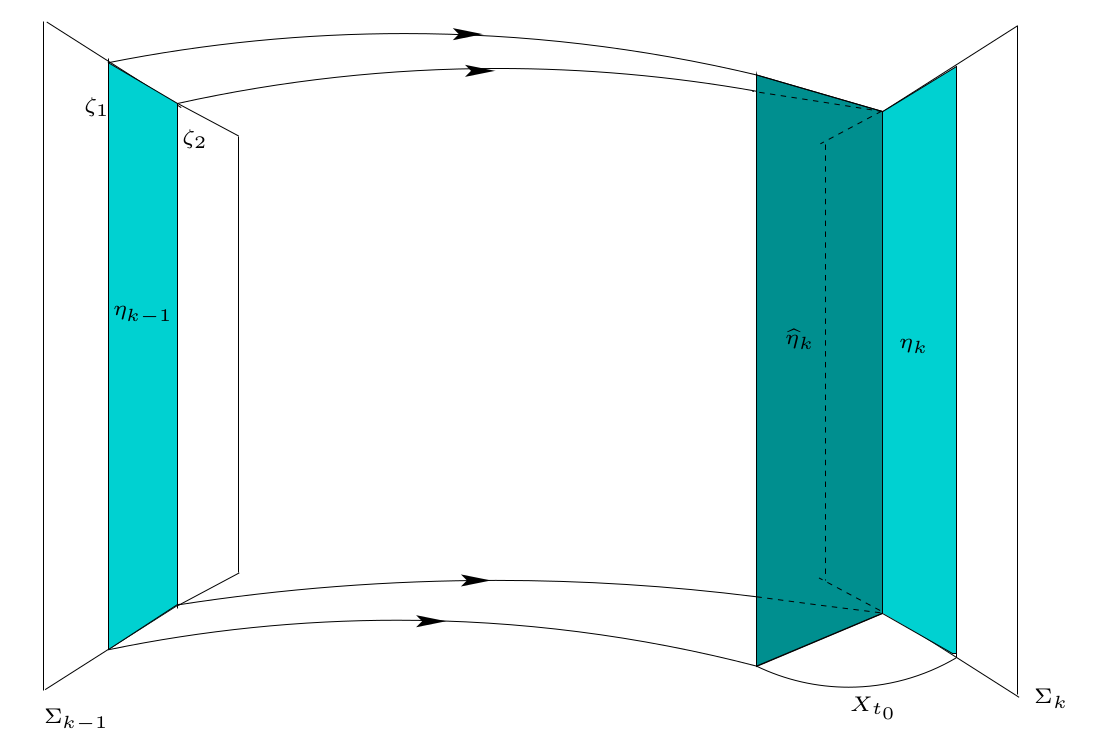}
    \caption{Case $\ell = k$ and traveling in a flow-box.}
    \label{fig:case-l-equal-k}
  \end{figure} Because the orbits are in a neighborhood of a
  Lorenz-like singularity we claim that
  $|\tau(f^{k - 1}(x_1)) - \tau(f^{k - 1}(x_2))|$ is bounded. Indeed,
  from Lemma \ref{lemma:property-poincare-time} there is a constant
  $C>0$ so that this expression is
  bounded from above by
  \begin{align*}
  C &\dfrac{|f^{k - 1}(x_1) - f^{k - 1}(x_2)|}{\min\{\dist(f^{k -
      1}(x_1), \cD), \dist(f^{k - 1}(x_2), \cD)\}}
  + |f^{k - 1}(x_1) - f^{k - 1}(x_2)|^\alpha
  \\
  &\leq C \left(\dfrac{\sigma^{R(x_1) - k + 1}}{\sigma^{b(R(x_1) - k +
        1)}} + \sigma^{\alpha(R(x_1) - k + 1)}\right)|F(x_1) -
  F(x_2)|^\alpha.
  \end{align*}
  Since $0 < b < 1/2$, it follows that the
  right hand side of the inequality above is bounded by a constant as
  we claimed.  Because
  $\tau(f^{k - 1}(x_1)) \leq \tau(f^{k - 1}(x_2))$, we have that
  $X_{\tau(f^{k - 1}(x_1))}(\zeta_1) \in \Sigma_{k}$ while
  $X_{\tau(f^{k - 1}(x_1))}(\zeta_2)$ is yet about to hit
  $\Sigma_{k}$; see Figure \ref{fig:case-l-equal-k}. Because
  $t_0:= |\tau(f^{k - 1}(x_1)) - \tau(f^{k - 1}(x_2))|$ is bounded we
  have that $\widehat{\eta}_{k}$ is diffeomorphic to $\eta_{k - 1}$ by
  a diffeomorphism that distorts distances at most by a constant
  factor. In particular, letting $\zeta_i = W^s_{f^{k - 1}(x_i)}$,
  $i = 1, 2$, there exists a constant $C > 0$ such that
  $\dist(X_{\tau(\zeta_1)}(\zeta_1), X_{\tau(\zeta_1)}(\zeta_2)) \leq
  C\dist(X_{\tau(\zeta_1)}(\zeta_1), X_{\tau(\zeta_2)}(\zeta_2)).$

  If we restrict the flow to a central-unstable invariant manifold
  $W=W_\sigma$ in a neighborhood of $\sigma$, as in the proof of
  Lemma~\ref{lemma:property-poincare-time} using a smooth
  linearization, then the points $\xi_1=W\cap\zeta_1$ and
  $\xi_2=W\cap\zeta_2$ move away from each other at a uniform rate,
  that is\footnote{Recall that $\tau$ is constant on stable leaves.}
  $\dist(X_{t}(\xi_1), X_{t}(\xi_2)) \leq Ce^{-\lambda (\tau(\zeta_1)
    - t)}\dist(X_{\tau(\zeta_1)}(\xi_1), X_{\tau(\zeta_2)}(\xi_2))$,
  for all $0<t<\tau(\zeta_1)$, where $\lambda>0$ is the expanding
  eigenvalue at the singularity.  Since the stable foliation is of
  class $C^2$ and transverse to $W$, we can write $V_\sigma$ as a
  product $V_\sigma=W_\sigma\times D$ where $D$ is a $d_s$-dimensional
  disk and the identification is given by a $C^{1+}$ diffeomorphism
  (the smoothness provided by
  Lemma~\ref{lemma:newhouse-linearization}). Hence we can extend the
  previous estimate for $\xi_1,\xi_2$ to the entire local stable leaf
  in $V_\sigma$ by at most a constant factor, that is
  $\dist(X_{t}(\zeta_1), X_{t}(\zeta_2)) \leq
  Ce^{-\lambda(\tau(\zeta_1)-t)}\dist(X_{\tau(\zeta_1)}(\zeta_1),
  X_{\tau(\zeta_2)}(\zeta_2))$ for $0<t<\tau(x_1)$.  In particular, we
  arrive at
  $ |X_\xi(P^{k - 1}(x_1,y_1)) - X_\xi(P^{k - 1}(x_2,y_2))| \leq
  C|P^k(x_1,y_1) - P^k(x_2,y_2)|.  $

  Finally
  $|X_u(x_1,y_2) - X_\xi(P^{k - 1}(x_2,y_2))| \leq |F(x_1) - F(x_2)| +
  |y_1 - y_2|$ by inequality \eqref{eq:ineq-aux-0} and the result
  follows.

 To conclude the proof we consider the case $x_1, x_2 \in Q$ with
 $|F(x_1) - F(x_2)| \geq \varepsilon$. Since $M$ is a compact
 we can set $K:=\sup_{(t, z) \in \RR \times M}|X_u(z)| < \infty$. Hence,
 \begin{align*}
   |X_u(x_1, y_1) - X_u(x_2,y_2)| 
   &\leq 2K \leq \dfrac{2K}{\varepsilon}|F(x_1)
   - F(x_2)| 
   \\
   &\leq C(|F(x_1) - F(x_2)|^\alpha + |y_1 - y_2|),
 \end{align*}
 letting $C > 0$ bigger so that $C > 2k / \varepsilon$ if necessary.
\end{proof}

We are finally ready to present 

\begin{proof}[Proof of Theorem~\ref{thm:conjugation}]
  By the Mean Value Theorem, there is a $u \in (u_1, u_2)$ such that
  \begin{align*}
    |p(x_2,y_2,u_1) - p(x_2,y_2,u_2)| = |X_{u_1}(x_2,y_2) - X_{u_2}(x_2,y_2)|
    \leq |X|_\infty|u_1 - u_2|.
  \end{align*}
  Now 
  $|p(x_1,y_1,u_1) - p(x_2,y_2,u_2)|$ is bounded above by
  \begin{align*}
  |p(x_1,y_1,u_1) - p(x_2,y_2,u_1)|
  &+ |p(x_2,y_2,u_1) - p(x_2,y_2,u_2)| 
  \\
  &\leq 
  |X_{u_1}(x_1,y_1) - X_{u_1}(x_2, y_2)| + |X|_\infty |u_1 - u_2|
  \\
  &\leq  C(|F(x_1) - F(x_2)|^\alpha + |y_1 - y_2| + |u_1 - u_ 2|),
  \end{align*}
  for some constant $C > 0$, as desired, after using
  Proposition~\ref{lemma:conjug-aux2}.
\end{proof}

\subsection{Exponential mixing}
\label{sec:exponential-mixing-1}

We deduce here
Theorem~\ref{mthm:decay-of-correlations-singular-hyperbolic} from
Theorem~\ref{theorem:conjugacy-skew-product-hyperbolic}.

\begin{proof}[Proof of Theorem~\ref{mthm:decay-of-correlations-singular-hyperbolic}]
  It follows from Theorem
  \ref{theorem:conjugacy-skew-product-hyperbolic} that if
  $\vfi\in C^1(U)$ and $\psi\in C^3(U)$, then
  $\varphi \circ p \in C^\alpha_\loc(\widehat{\Delta}^r)$ and
  $\psi \circ p \in C^{\alpha, 2}_\loc(\widehat{\Delta}^r)$ for each
  skew product semiflow associated to each ergodic physical measure
  supported in $U$.

  Hence, using item 3 of
  Theorem~\ref{theorem:existence-physical-measures} we get $s_i\ge0$,
  $\sum_{i=1}^ks_i=1$ and $\mu=\sum_is_i\mu_i$ where each $\mu_i$ is
  an ergodic physical measure for $X$ supported in the attracting set
  $\Lambda$. Theorem~\ref{theorem:conjugacy-skew-product-hyperbolic}
  ensures that
  $\mu=\sum_i s_i \cdot (p_i)_*\big(\mu_{\hat F_i}^{r_i}\big)$.  We
  normalize the observable $\vfi$ by defining
$
    \bar\vfi=\sum_{i=1}^k s_i(\vfi-\mu_i(\vfi))\chi_{B(\mu_i)}
$
  which satisfies $\mu_i(\bar\vfi)=0$ and also
  $\bar\vfi\circ p_i=s_i\big(\vfi\circ p_i-\mu_i(\vfi)) \in
  C^\alpha_{loc}(\widehat\Delta_i^{r_i})$ for all $i=1,\dots, k$.

  Combining this with
  Theorem~\ref{prop:decay-correlations-skew-product-weaker} we
  conclude that
  \begin{align*}
    \Big| \int (\varphi \circ X_t) \psi\, d\mu
    &- \int \varphi\, d\mu \int
      \psi\, d\mu \Big|
      =
      \left|\int\big(\varphi \circ X_t-\mu(\vfi))\psi\,d\mu\right|
    \\
    &=
      \left|
      \int\bar\vfi\circ X_t\cdot\psi\,d\mu
      \right|
      \le\sum_is_i
      \left|\int (\bar\vfi \circ p_i \circ \widehat{F_i}^t) \cdot \psi \circ
      p_i\, d\mu^{r_i}_{\widehat F_i} \right|
    \\
    &\leq
      \sum_i s_i C_ie^{-c_i t}
      \|\bar\vfi \circ   p_i\|_{\alpha} \|\psi \circ p_i\|_{\alpha, 2}
    \\
    &\le
      C e^{-ct}\sum_i s_i \|\vfi \circ   p_i\|_{\alpha} \|\psi \circ p_i\|_{\alpha, 2}
  \end{align*}
  for some $C,c>0$, since the number $k$ of ergodic physical measures
  is finite and $\vfi$ is bounded.  Using Theorem
  \ref{theorem:conjugacy-skew-product-hyperbolic} again we get
  \begin{align*}
  \left|\int (\varphi \circ X_t) \cdot \psi\, d\mu - \int \varphi\,
    d\mu\int \psi \, d\mu\right| \leq C^3e^{-c t}|\vfi|_{C^1}
  |\psi|_{C^3}.
  \end{align*}
  Finally, let us fix $\eta\in(0,1)$ and $\vfi,\psi\in C^\eta(U)$.
  Given $\delta>0$ we can choose $\tilde\vfi,\tilde\psi\in C^3(U)$
  such that $|\vfi-\tilde\vfi|_\infty<\delta^{\eta}\|\vfi\|_{\eta}$ and
  $|\tilde\vfi|_{C^1}\le\delta^{-1} \|\vfi\|_\eta$; and also
  $|\psi-\tilde\psi|_\infty<\delta^{\eta}\|\psi\|_{\eta}$ with
  $|\tilde\psi|_{C^3}\le\delta^{-3} \|\psi\|_\eta$.

  Then, if we denote
  $\rho(\vfi,\psi,t)=\mu (\psi\cdot\varphi \circ
  X_t)-\mu(\vfi)\mu(\psi)$ and the constants of the last estimate as
  $\tilde C,\tilde c>0$, we get
  \begin{align*}
    |\rho(\vfi,\psi,t)
    &-\rho(\tilde\vfi,\tilde\psi,t)|
    \le
    2|\vfi-\tilde\vfi|_\infty|\psi|_\infty
    +
    2|\tilde\vfi|_\infty|\psi-\tilde\psi|_\infty
    \\
    &\le
    2 \|\vfi\|_{\eta}\delta^{\eta}\|\psi\|_\eta
    +
      2|\tilde\vfi|_\infty \|\psi\|_{\eta}\delta^{\eta}
      \le
      2 \|\vfi\|_{\eta}\|\psi\|_\eta\delta^{\eta}
      +
      4 \|\vfi\|_\eta \|\psi\|_{\eta}\delta^{\eta}
  \end{align*}
  since $|\tilde\vfi|_\infty\le|\vfi|_\infty+|\vfi-\tilde\vfi|_\infty
  \le \|\vfi\|_\eta(1+\delta^\eta)\le 2\|\vfi\|_\eta$.
  Moreover, we also have
  $ |\rho(\tilde\vfi,\tilde\psi,t)| \le \tilde C e^{-\tilde
    ct}|\tilde\vfi|_{C^3}|\tilde\psi|_{C^3} \le \tilde C e^{-\tilde
    ct}\delta^{-4}\|\vfi\|_{\eta}\|\psi\|_{\eta} $, thus
  \begin{align*}
    |\rho(\vfi,\psi,t)|
    &\le
      \|\vfi\|_{\eta}\|\psi\|_{\eta}
      (\tilde C\delta^{-4} e^{-\tilde ct}+6\delta^\eta).
  \end{align*}
  Setting $\delta=e^{-\tilde c t / (4+\eta)}$ we obtain a constant
  $C>0$ so that
  \begin{align*}
    |\rho(\vfi,\psi,t)|
    &\le
      Ce^{-\eta \tilde c t / (4+\eta)}\|\varphi\|_\eta\|\psi\|_\eta, \quad t>0.
  \end{align*}
  This completes the proof after setting the exponent $c=\eta\tilde c/(4+\eta)$.
\end{proof}

\subsection{Exponential convergence to equilibrium}
\label{section:proof-convergence-equilibrium}

To prove exponential convergence to equilibrium for the flow, that is,
Corollary~\ref{mcor:convergence-equilibrium}, we need the following
corollary of Theorem \ref{prop:decay-correlations-skew-product-weaker}
whose proof we postpone to Subsection
\ref{sub:convergence-equilibrium-semiflow}.  We denote the Lebesgue
measure in $\widehat{\Delta}^r$ by $\Leb^r_3$, that is,
$\Leb^{r_i}_3:= (\Leb_{\widehat{\Delta_i}} \times \Leb_\RR) / \int
r_i\, d\Leb_{\widehat{\Delta}}$ corresponding to each one of the
ergodic physical measures $\mu_i$ supported on the attracting set.

\begin{corollary}[Exponential convergence to equilibrium for
  $\widehat{F}_t$]
  \label{corollary:convergence-to-equilibrium-skew-product-flow}
  In the same setting of Theorem \ref{prop:decay-correlations-skew-product-weaker} there
  exist constants $c, C > 0$ such that
  $$
  \left| \int (\varphi \circ \widehat{F}_t) \psi \, d\Leb^r_3 - \int
    \varphi d\mu^r_{\widehat{F}} \int \psi\, d\Leb^r_3 \right| \leq
  Ce^{-c t}\|\varphi\|_\alpha\|\psi\|_{\alpha,2},$$ for all
  $\varphi \in C^{\alpha}_{\loc}(\widehat{\Delta}^r)$,
  $\psi \in C^{\alpha, 2}_{\loc}(\widehat{\Delta}^r)$ and $t > 0$.
\end{corollary}

Now we use this to complete the proof of the remaining main result.

\begin{proof}[Proof of Corollary~\ref{mcor:convergence-equilibrium}]
  We argue similarly to the proof of
  Theorem~\ref{mthm:decay-of-correlations-singular-hyperbolic} using
  the decomposition $\Leb=\sum_i\vartheta_i\Leb_i$ and
  $\widetilde\mu=\sum_i\vartheta_i\mu_i$ to write for
  $\vfi,\psi\in C^1(U)$
  \begin{align*}
    \Big| \int (\varphi \circ X_t) \psi\, d\Leb
    &- \int \varphi\, d\mu \int
      \psi\, d\Leb \Big|
      =
      \left|\int\big(\varphi \circ X_t-\mu(\vfi))\psi\,d\Leb\right|
    \\
    &=
      \left|
      \int\bar\vfi\circ X_t\cdot\psi\,d\Leb
      \right|
      \le
      \sum_i\vartheta_i
      \left|
      \int\bar\vfi\circ X_t\cdot\psi\,d\Leb_i
      \right|
  \end{align*}
  where $\bar\vfi$ is just as in the proof of
  Theorem~\ref{mthm:decay-of-correlations-singular-hyperbolic}.
  Now we have $\Leb_i=\xi_i\cdot (p_i)_*\Leb_3^{r_i}$ where $\xi_i$ is the
  Jacobian of $p_i(w,u)=X_u(w)$, which depends on the Jacobian of the
  flow of $X$, which is of class $C^1$ since the vector field $X$ is
  of class $C^2$.

  Moreover, $\xi_i$ is strictly positive and uniformly bounded since,
  by strong dissipativeness, we have that the divergence of the vector
  field is strictly negative in a neighborhood of $\Lambda$: there
  exists $\vartheta>0$ so that $\Div X<-\vartheta<0$ on $U$. Hence
  \begin{align*}
    \xi_i(w,u)
    =|\det Dp_i(w,u)|
    &\le
      \frac{\|X(X_uw)\|}{\|X(w)\|}\cdot|\det DX_u(w)|
    \\
    &=
      \frac{\|X(X_uw)\|}{\|X(w)\|}\cdot\exp\int_0^u(\Div
      X)(X_s(w))\,ds
      \le
      C e^{-\vartheta u}
  \end{align*}
  where the constant $C>0$ depends only on the lenght of the vector
  field in a neighborhood of $\Lambda$. Thus $\xi_i\in C^1(U)$.
  We can therefore write
  \begin{align*}
      \sum_i\vartheta_i
      &\left|
      \int\bar\vfi\circ X_t\cdot\psi\,d\Leb_i
        \right|
        =
        \sum_i\vartheta_i\left|
        \int\bar\vfi\circ X_t\cdot(\psi\xi_i)\,d(p_i)_*\Leb_3^{r_i}
        \right|.
  \end{align*}
  At this point, we approximate $\psi\xi_i$ by a $C^3$ function: for a
  given $0<\delta<1$ we choose $\tilde\psi\in C^3(U)$ so that
  $ |\psi\xi_i-\tilde\psi|_\infty\le\delta|\psi\xi_i|_{C^1} $ and
  $ |\tilde\psi|_{C^3}\le\delta^{-3}|\psi\xi_i|_{C^1}.$
  On the one hand
  \begin{align*}
    \left|
    \int\!\! \tilde \psi\cdot\bar\vfi\circ X_t \,d(p_i)_*\!\Leb_3^{r_i}
    \!\!-\!\!
    \int\!\!(\psi\xi_i)\cdot\bar\vfi\circ X_t\,d(p_i)_*\!\Leb_3^{r_i}
    \right|
    \le
    |\bar\vfi|_\infty|\tilde\psi-\psi\xi_i|_\infty
    \le
    \delta|\bar\vfi|_{C^1}|\psi\xi_i|_{C^1}
  \end{align*}
  while on the other hand, by
  Corollary~\ref{corollary:convergence-to-equilibrium-skew-product-flow}
  and Theorem~\ref{theorem:conjugacy-skew-product-hyperbolic}
  \begin{align*}
    \left|
    \int\bar\vfi\circ X_t\cdot \tilde \psi \,d(p_i)_*\Leb_3^{r_i}
    \right|
    &=
      \left|\int (\bar\vfi \circ p_i \circ \widehat{F_i}^t) \cdot \tilde\psi \circ
      p_i\, d\Leb_3^{r_i} \right|
    \\
    &\le
      C_ie^{-c_it}\|\bar\vfi \circ   p_i\|_{\alpha} \|\tilde\psi \circ
      p_i\|_{\alpha, 2}
    \\
    &\le
      C^2C_ie^{-c_it}|\vfi |_{C^1} |\tilde\psi|_{C^3}
      \le
      C^2C_ie^{-c_it} \delta^{-3}|\vfi |_{C^1} |\psi\xi_i|_{C^1}.
  \end{align*}
  So we obtain
  \begin{align*}
    \left|
    \int\bar\vfi\circ X_t\cdot(\psi\xi_i)\,d(p_i)_*\Leb_3^{r_i}
    \right|
    &\le
    C^2C_ie^{-c_it} \delta^{-3}|\vfi |_{C^1} |\psi\xi_i|_{C^1}
    +
      \delta|\bar\vfi|_{C^1}|\psi\xi_i|_{C^1}
    \\
    &\le
      \tilde C_i e^{-c_it/4}|\vfi |_{C^1} |\psi|_{C^1}
  \end{align*}
  for some constant $\tilde C_i>0$, after setting
  $\delta=e^{-c_it/4}$.  Since this holds for each $i$, we get
  \begin{align*}
        \Big| \int (\varphi \circ X_t) \psi\, d\Leb
    &- \int \varphi\, d\mu \int
      \psi\, d\Leb \Big|
      \le
      \sum_is_i\left|
      \int\bar\vfi\circ X_t\cdot(\psi\xi_i)\,d(p_i)_*\Leb_3^{r_i}
      \right|.
  \\
    &\le
      \sum_i s_i C_ie^{-c_i t/4}
      |\vfi|_{C^1} |\psi|_{C^1}
    \le
      C e^{-ct}|\vfi|_{C^1} |\psi|_{C^1}
  \end{align*}
  for some constants $C,c>0$ by
  finiteness of the number $k$ of physical measures.

  Having established the result for smooth observables
  $\vfi,\psi\in C^1(U)$, we can now extend it to H\"older observables
  $\vfi,\psi\in C^\eta(U)$ for any $0<\eta<1$ using the exact same
  arguments as in the proof of
  Theorem~\ref{mthm:decay-of-correlations-singular-hyperbolic}.
\end{proof}


\section[Exponential mixing for semiflows]{Exponential mixing and
  convergence to equilibrium for hyperbolic skew product semiflows}
\label{chapter:decay-correlations-suspension-general}

In this chapter we present the proof of Theorems
\ref{theorem:decay-onedimensional-weaker} and
\ref{prop:decay-correlations-skew-product-weaker}. Because the proof
of these theorems follow the same steps as of \cite{ArMel16} we only
prove the parts that differ and refer to parts that are equal.  


From now on use the following convention: given two real sequences
$(a_n)_n$ and $(b_n)_n$, we write $a_n \lesssim b_n$ if there is a
constant $C > 0$ such that $a_n \leq Cb_n$, for all $n \geq 1$.

\subsection{Exponential mixing for {\boldmath $C^{1 + \alpha}$}
  expanding semiflows}

In this section we prove Theorem \ref{theorem:decay-onedimensional-weaker}. This
theorem is a generalization of \cite[Theorem 2.1]{ArMel16} to the function space
$C^{\alpha, 2}_\loc(\Delta^r)$ (recall the definition of this space on Subsection
\ref{subsec:decay-expanding-semiflows}). In the proof we use the results of
\cite{ArMel16} as much as possible and show the adaptations in the places where
they are required.

Throughout this section we consider $F:\Delta \to \Delta$ to be a
$C^{1 + \alpha}$ uniformly expanding map and
$r: \Delta \to (0, +\infty)$ a function satisfying conditions (iii) -
(v) (recall Subsection \ref{subsection:uniformly-expanding-maps}).
Setting $r_n = \sum^{n - 1}_{j = 0}r \circ F^j$, we note that we can
generalize the items (ii) and (iii) of Subsection
\ref{subsection:uniformly-expanding-maps} as: there exists a constant
$C > 0$ such that ($\text{ii}_1$) $|\log |h'||_\alpha \leq C$; and
$(iii)_1$ $|(r_n \circ h)'|_\infty \leq C$
for all $h \in \cH_n$ and all integer $n \geq 1$

We need to use an equivalent form of the  UNI condition (see~\cite[Proposition
7.4]{AvGoYoc}):
\begin{description}
\item[UNI - equivalent formulation] there exist $D > 0$, $n_0 \geq 1$
sufficiently large and $h_1, h_2 \in \cH_{n_0}$ so that
$\inf|(r_n \circ h_1 - r_n \circ h_2)'| \geq D.$
\end{description}

\subsubsection{Twisted transfer operator}
\label{sec:twisted-transf-opera}

Here we work with complex observables so we denote by $C^\alpha(\Delta)$ the
space of functions $\psi: \Delta \to \CC$ such that $\|\psi\|_\alpha =
|\psi|_\infty + |\psi|_\alpha < \infty$ and $C^\alpha_{\loc}(\Delta)$ the space
of functions $\psi: \Delta \to \CC$ such that $\|\psi\|_\alpha = |\psi|_\infty +
|\psi|_{\alpha, \loc} < \infty$, where
$$|\psi|_{\alpha, \loc} = \sup_{h \in \cH}\sup_{x \neq y}\dfrac{|\psi(hx) -
  \psi(hy)|}{|x -  y|^\alpha}.$$ It is also convenient to introduce the family
of equivalent norms: for all $b \in \RR$
\begin{align*}
\|\psi\|_b = \max\{|\psi|_\infty, |\psi|_{\alpha, \loc} / (1 +
|b|^\alpha)\}, \quad \psi \in C^\alpha_\loc(\Delta).
\end{align*}
For each $s \in \CC$ we denote by $P_s$ the non-normalized twisted
transfer operator, that is,
\begin{align*}
P_s = \sum_{h \in \cH}A_{s, h} \quad\text{where}\quad A_{s,h}\psi = e^{-s r \circ
  h}|h'|\psi \circ h.
\end{align*}
In what follows we present the result that guarantees that $P_s$ is
well defined.

\begin{proposition}\cite[Proposition
  2.5]{ArMel16}\label{prop:twisted-operator}
  Write $s=\sigma+ib$.
  There exists $\varepsilon\in(0,1)$ such that
  the family $s\mapsto P_s$ of operators on $C^\alpha(\Delta)$ is continuous on
  $\{\sigma>-\varepsilon\}$.
  Moreover, $\sup_{|\sigma|<\varepsilon}\|P_s\|_b<\infty$.
\end{proposition}

\begin{remark}\label{remark:decay-local-holder}
  Using that $|\psi(hx) - \psi(hy)| \leq
  |\psi|_{\alpha, \loc}|x - y|^\alpha$ in the proof of Proposition
  \ref{prop:twisted-operator} we note that we can define $P_s$ on
  $C^\alpha_\loc(\Delta)$ and its range remain in $C^\alpha(\Delta)$ when $\Re s
  > -\varepsilon$. Moreover, it also follows from the proof that
  $\|P_s\psi\|_b \leq \widetilde{C}\|\psi\|_b$, where $\widetilde{C}:=
  C(1 + |\sigma| + |b|^\alpha)(1 +
  |b|^\alpha)^{-1}\sum_{h \in \cH}e^{\varepsilon |r \circ h|_\infty}|h'|_\infty$.
\end{remark}

The unperturbed operator $P_0$ has a simple leading eigenvalue
$\lambda_0 = 1$ with strictly positive $C^\alpha$ eigenfunction $f_0$. By
Proposition \ref{prop:twisted-operator}, there exists $\varepsilon \in
(0, 1)$ such that $P_\sigma$ has a continuous family of simple
eigenvalues $\lambda_\sigma$ for $|\sigma| < \varepsilon$ with
associated $C^\alpha$ eigenfunctions $f_\sigma$. Shrinking
$\varepsilon$ if necessary, we can ensure that $\lambda_\sigma > 0$
and $f_\sigma$ is strictly positive for $|\sigma| < \varepsilon$.

\begin{remark}\label{remark:pertubation-eigenthings}
  By standard perturbation theory, for any $\delta > 0$ there exists
  $\varepsilon \in (0,1)$ such that $\sup_{|\sigma| <
    \varepsilon}|\lambda_\sigma - 1| < \delta$, $\sup_{|\sigma| <
    \varepsilon}|f_\sigma/f_0 - 1|_\infty < \delta$ and
  $\sup_{|\sigma| < \varepsilon}|f_\sigma / f_0 - 1|_\alpha < \delta$.
  In particular, we may assume that $1/2 \leq \lambda_\sigma \leq 2$,
  $f_0 / 2 \leq f_\sigma \leq 2f_0$ and $|f_0|_\alpha / 2 \leq
  |f_\sigma|\alpha \leq 2 |f_0|_\alpha$.
\end{remark}

For $s = \sigma + i b$ with $|\sigma| < \varepsilon$, we define the
normalized transfer operators
\begin{align*}
  \mathcal{L}_s\psi = (\lambda_\sigma f_\sigma)^{-1}P_s(f_\sigma
  \psi) = (\lambda_\sigma f_\sigma)^{-1}\sum_{h \in \cH}A_{s,h}(f_\sigma
  \psi).
\end{align*}
It also follows that $\mathcal{L}_s$ is defined on
$C^{\alpha}_\loc(\Delta)$ with range on $C^\alpha(\Delta)$ for all $s$
where $P_s$ is defined. Note that $\mathcal{L}_\sigma1 = 1$ for all
$\sigma$ and $|\mathcal{L}_s|_\infty \leq 1$ for all $s$ (where
defined).

Using Remark \ref{remark:decay-local-holder}, the strategy that we
follow now is to take the results of \cite{ArMel16} that depend on
observables $\psi \in C^\alpha(\Delta)$, change them to $P_s\psi$ or
$\cL_s\psi$, with $\psi \in C^\alpha_\loc(\Delta),$ and explain why
these changes are enough to prove Theorem
\ref{theorem:decay-onedimensional-weaker}.

\subsubsection{Lasota-Yorke inequality}

Note that $$\mathcal{L}^n_s\psi = \lambda^{-n}_\sigma f^{-1}_\sigma\sum_{h
  \in \cH_n}A_{s, h, n}(f_\sigma\psi), \text{  where  } A_{s,h,n}\psi := e^{-sr_n
  \circ h}|h'|\psi \circ h.$$

\begin{lemma}\cite[Lemma 2.7]{ArMel16}
  \label{lemma:lasota-yorke-original}
  There is a constant $C > 1$ such that
  $$|\cL^n_s\psi|_\alpha \leq C(1 + |b|^\alpha)|\psi|_\infty + C
  \rho^n|\psi|_\alpha \leq C(1 + |b|^\alpha)\{|\psi|_\infty +
  \rho^n\|\psi\|_b\},$$ for all $s = \sigma + i b$, $|\sigma| < \varepsilon$, $n
  \geq 1$ and all
  $\psi \in C^\alpha(\Delta)$.
\end{lemma}

\begin{lemma}\cite[Lemma 2.7
  adapted]{ArMel16} \label{lemma:lasota-yorke}
  There is a constant $C > 1$ such that
  $$|\mathcal{L}^{n + 1}_s\psi|_{\alpha} \leq C(1 + |b|^\alpha)|\psi|_\infty +
  C\rho^{n}|\psi|_{\alpha, \loc} \leq C(1 + |b|^\alpha)\{|\psi|_\infty
  + \rho^{n}\|\psi\|_b\},$$ for all $s = \sigma + i b$,
  $|\sigma| < \varepsilon$, $n \geq 1$ and all
  $\psi \in C^\alpha_\loc(\Delta)$.
\end{lemma}

\begin{proof}
  As $\cL_s \psi \in C^\alpha(\Delta)$ for all
  $\psi \in C^\alpha_\loc(\Delta)$, using the original result
  \cite[Lemma 2.7]{ArMel16} we get that there exists $C > 1$ so that
  $|\cL^{n + 1}_s\psi|_\alpha \leq C(1 + |b|^\alpha)|\cL_s\psi|_\infty
  + C\rho^n|\cL_s \psi|_\alpha.$ Using Remarks
  \ref{remark:decay-local-holder} and
  \ref{remark:pertubation-eigenthings} we get that
  $|\cL_s \psi|_\alpha \lesssim |\psi|_{\alpha, \loc}$ for all
  $\psi \in C^\alpha_\loc(\Delta)$. The results follows since we also
  have that $|\cL_s|_\infty \leq 1$.
\end{proof}

\begin{corollary}\label{corollary:bounded-normalized-transfer-operator}
  There exists $C > 1$ such that  $\|\cL_s^n\|_b \leq C$ for all $s =
  \sigma + i b$, with $|\sigma| < \varepsilon$ and all $n \geq 1$.
\end{corollary}

\begin{proof}
  It is clear that
  $|\cL^n_s\psi|_\infty \leq |\psi|_\infty \leq \|\psi\|_b$, for all
  $n \geq 1$. By Lemma \ref{lemma:lasota-yorke}, there exists $C > 1$
  such that $|\cL^n_s\psi|_\alpha \leq C(1 + |b|^\alpha)\|\psi\|_b$,
  for all $n \geq 2$, for all $\psi \in C^\alpha_\loc(\Delta)$. For
  $n = 1$ we have by Remark \ref{remark:decay-local-holder} that there
  exists $\widetilde{C} > 1$ such that
  $|\cL_s\psi|_\alpha \leq \widetilde{C}|\psi|_{\alpha, \loc} \leq
  \widetilde{C}\|\psi\|_b$, for all $\psi \in
  C^\alpha_\loc(\Delta)$. Hence, taking a bigger $C$ the result
  follows.
\end{proof}

In the following lemma the constant $C_0$ is the constant $C_4$ from the
definition of the family of cones $\cC_b$ in \cite[Subsection 2.3]{ArMel16}. For
our needs it is not necessary to enter in the details of cone invariance because we
only need to adapt the existing results of \cite{ArMel16}.

\begin{lemma}\cite[Corollary 2.15
  adapted]{ArMel16} \label{lemma:almost-contraction}
  There exist $\varepsilon, \beta \in (0,1)$ and $A, C > 0$ 
  such that $\|\cL^{4mn_0 + 1}_s\psi\|_b \leq C\beta^m\|\psi\|_b$ for
  all $m \geq A \log |b|$, $s = \sigma + ib$, $|\sigma| <
  \varepsilon$, $|b| \geq \max\{4\pi/ D, 1\}$ and all $\psi \in
  C^\alpha_\loc(\Delta)$ satisfying $|\cL_s\psi|_{\alpha} \leq
  C_0|b|^\alpha|\cL_s\psi|_\infty$.
\end{lemma}

\begin{proof}
  The result follows directly from \cite[Corollary 2.15]{ArMel16} using
  $\cL_s\psi$ in the place of $\psi$ for $\psi \in C^\alpha_\loc(\Delta)$.
\end{proof}

For the next theorem, if $C$ is the constant given by Lemma
\ref{lemma:lasota-yorke-original} and $\widetilde{C} > 0$ is such that $\|\cL_s\psi\|_b
\leq \widetilde{C}\|\psi\|_b$, for all $\psi \in C^\alpha_\loc(\Delta)$, we
require that $\max\{2C\widetilde{C}C^{-1}_0, 2\widetilde{C}C^{-1}_0\} < 1/3$. We also fix $n_0$
from UNI condition such that $C\widetilde{C}\rho^{n_0 - 1} < 1/3$.

\begin{theorem}\label{theorem:contraction-non-normalized-transfer-operator}
  Let $D' = \max\{4\pi / D, 2\}$. There exist $\varepsilon, \gamma \in (0, 1)$
  and $A > 0$ such that $\|P^n_s\|_b \leq \gamma^n$ for all $s = \sigma + ib$,
  $|\sigma| < \varepsilon$, $|b| \geq D'$ and $n \geq A \log|b|.$
\end{theorem}

\begin{proof}
  We claim that there exist constants $\varepsilon, \gamma_1 \in (0,1)$ and $A,
  C > 0$ such that $\|\cL^{4mn_0 + 1}_s\|_b \leq C\gamma^m_1$ for all $s = \sigma
  + i b$, $|\sigma| < \varepsilon$, $|b| \geq \max\{4\pi / D, 2\}$ and $m \geq A
  \log|b|$. Indeed, it is enough to prove that $\|\cL^{4mn_0 + 1}_s\psi\|_b \leq
  C \gamma^m_1\|\psi\|_b$ with $\psi \in C^\alpha_\loc(\Delta)$ satisfying
  $|\cL_s\psi|_{\alpha} > C_0|b|^\alpha|\cL_s\psi|_\infty$, otherwise the result
  follows from Lemma \ref{lemma:almost-contraction}. We have
  that
  $$
  \begin{array}{rcl}
  |\cL^{n_0}_s\psi|_\infty & \leq & |\cL_s\psi|_\infty \leq
  (C_0|b|^\alpha)^{-1}|\cL_s\psi|_{\alpha} \leq (C_0|b|^\alpha)^{-1}(1 +
                                    |b|^\alpha)\|\cL_s\psi\|_b \\
    & \leq & 2\widetilde{C}C^{-1}_0\|\psi\|_b \leq  \dfrac{1}{3}\|\psi\|_b.
  \end{array} $$ By Lemma
  \ref{lemma:lasota-yorke-original} there exists $C > 1$ such that 
  $$
  \begin{array}{rcl}
    |\cL^{n_0}_s \psi|_\alpha & \leq & C(1 +
                                       |b|^\alpha)\{|\cL_s\psi|_\infty +
                                       \rho^{n_0 - 1}\|\cL_s\psi\|_b\}\\
    & \leq & (1 +
                                       |b|^\alpha)\{2CC^{-1}_0\|\cL_s\psi\|_b +
                                       C\rho^{n_0 - 1}\|\cL_s\psi\|_b\}\\
                              & \leq & (1 + |b|^\alpha)
                                       \{2C\widetilde{C}C^{-1}_0\|\psi\|_b +
                                       C\widetilde{C}\rho^{n_0 - 1}\|\psi\|_b\}
                                       \leq 
                                       \dfrac{2}{3} (1 + |b|^\alpha) \|\psi\|_b.
  \end{array}
  $$ Hence $\|\cL^{n_0}\|_b \leq 2/3$ and now the claim follows.

  Write $n = 4mn_0 + 1 + r$, with $0 \leq r < 4n_0 - 1$ in the case
  that $n$ is not multiple of $4n_0$ and $r = 4n_0 - 1$, otherwise.
  If $m \geq A \log|b|$, using the claim above and Corollary
  \ref{corollary:bounded-normalized-transfer-operator} we get
  $\|\cL^n_s\|_b \leq \|\cL^r_s\|_b\|\cL^{4mn_0 + 1}_s\|_b \lesssim
  \gamma^m_1 \lesssim (\gamma^{1 / 4n_0}_1)^n.$

  By definition we have
  $P_s \psi = \lambda_\sigma f_\sigma \cL_s(f^{-1}_\sigma\psi)$, so
  using the fact that $\|f_\sigma\|_\alpha$ and
  $\|f^{-1}_\sigma\|_\alpha$ are bounded for $|\sigma| < \varepsilon$
  we obtain
  $\|P^n_s\|_b \lesssim \lambda^n_\sigma\|\cL^n_s\|_b \lesssim
  (\gamma^{1 / 4n_0}_1 \lambda_\sigma)^n$. Taking a smaller $\delta$
  in Remark \ref{remark:pertubation-eigenthings} we can assume that
  $\gamma:=\gamma^{1 / 4n_0}\lambda_\sigma < 1$ for all
  $|\sigma| < \varepsilon$.  Thus, making $A$ bigger if necessary,
  there exists a constant $\widetilde{C} > 0$ such that
  $\|P^n_s\|_b \leq \widetilde{C}\gamma^n$, for all
  $n \geq A \log|b|$. Finally, we can increase $A$ and modify $\gamma$
  to absorb the constant $\widetilde{C}$, finishing the theorem.
\end{proof}

\subsubsection{Proof of Theorem
  \ref{theorem:decay-onedimensional-weaker}}
\label{sub:proof-theorem-decay-onedimensional-weaker}

Here we denote the \emph{correlation function} of $F_t$ with respect
to observables $\varphi \in L^\infty(\Delta^r)$ and
$\psi \in L^1(\Delta^r)$ by
\begin{align*}
\rho_{\varphi, \psi}(t) = \int_{\Delta^r} (\varphi \circ F) \psi \, d\mu^r_F -
\int_{\Delta^r} \varphi \, d\mu^r_F \int_{\Delta^r} \psi\, d\mu^r_F.
\end{align*}
We denote by
$\widehat{\rho}_{\varphi, \psi}(s) = \int^\infty_0 e^{-st}
\rho_{\varphi, \psi}(t)\, dt$ the Laplace Transform of the correlation
function. Note that if $\varphi, \psi \in L^\infty(\Delta^r)$ then
$\widehat{\rho}_{\varphi, \psi}$ is well-defined and analytic on
$\{s \in \CC:\ \Re(s) > 0\}$ since
$|\rho_{\varphi, \psi}(t)| \leq 2 |\varphi|_\infty |\psi|_\infty$ for
all $t > 0$. The key estimate is 

\begin{lemma}\cite[Lemma 2.17
  adapted]{ArMel16}\label{lemma:key-estimate-laplace-transform}
  There exists $\varepsilon > 0$ such that $\widehat{\rho}_{\varphi,
    \psi}$ is analytic on $\{s \in \CC:\ \Re(s) > - \varepsilon\}$ for
  all $\varphi \in L^\infty(\Delta^r)$ and all $\psi \in
  C^\alpha_{\loc}(\Delta^r)$. Moreover, there is a constant $C > 0$
  such that $|\widehat{\rho}_{\varphi, \psi}(s)| \leq C(1 +
  |b|^{1/2})|\varphi|_\infty \|\psi\|_{\alpha, 2}$ for all $s = \sigma + i
  b$ with $\sigma \in [-\varepsilon / 2, 0]$.
\end{lemma}

Assuming that the Lemma \ref{lemma:key-estimate-laplace-transform}
holds we explain how to prove Theorem
\ref{theorem:decay-onedimensional-weaker} following the steps of
\cite{ArMel16}.

\begin{proof}[Proof of Theorem \ref{theorem:decay-onedimensional-weaker}]
  By Lemma \ref{lemma:key-estimate-laplace-transform},
  $\widehat{\rho}_{\varphi,\psi}$ is analytic on
  $\{s \in \CC:\ \Re(s) > -\varepsilon\}$. The inversion formula gives
  \begin{eqnarray}\label{eq:inversion-formula-laplace}
    \rho_{\varphi, \psi}(t) = \int_{\Gamma} e^{s t}
    \widehat{\rho}_{\varphi, \psi}(s)\, ds,
  \end{eqnarray} where we can take
  $\Gamma=\{s \in \CC:\ \Re(s) = -\varepsilon / 2\}$.  Taylor's Theorem on
  $\rho_{\varphi, \psi}$ provides
  \begin{align*}
  \rho_{\varphi, \psi}(t) = \rho_{\varphi, \psi}(0) + \rho_{\varphi,
    \partial_t \psi}(0)t + \int^t_0 \dfrac{(t - v)^2}{2}\rho_{\varphi,
    \partial^2_t \psi}(v)\, dv.
  \end{align*}
  Applying the Laplace Transform on the expression above we
  get
  \begin{align*}
  \widehat{\rho}_{\varphi, \psi}(s) = \rho_{\varphi, \psi}(0)s^{-1}
  + \rho_{\varphi, \partial_t \psi}(0)s^{-2} +
  \widehat{\rho}_{\varphi, \partial^2_t \psi}(t)s^{-3},
  \end{align*}
  where in the last term we have used that
  $t \mapsto \int^t_0 \frac{(t - v)^2}{2}\rho_{\varphi, \partial^2_t
    \psi}(v)\, dv$ is the convolution function between
  $v \mapsto \rho_{\varphi, \partial^2_t \psi}(v)$ and
  $v \mapsto v^2 / 2$. Thus, from~\eqref{eq:inversion-formula-laplace}
  and using the estimate given by Lemma
  \ref{lemma:key-estimate-laplace-transform}, we get that
  $|\rho_{\varphi, \psi}(t)|$ is bounded above by
  \begin{align*}
    \Big|\rho_{\varphi, \psi}(0)
    &\int_\Gamma \dfrac{e^{s t}}{s}\, ds + \rho_{\varphi,
      \partial_t\psi}(0) \int_\Gamma\dfrac{e^{s t}}{s^2}\, ds +
    \int_\Gamma \dfrac{e^{s t}\rho_{\varphi, \partial^2_t
        \psi}(s)}{s^3}\, ds\Big|
  \\
  &\leq C_1 e^{-\frac{\varepsilon
      t}{2}}|\varphi|_\infty|\psi|_\infty + C_2e^{-\frac{\varepsilon
      t}{2}}|\varphi|_\infty|\partial_t\psi|_\infty + C_3C(1 +
  |b|^{1/2}) e^{-\frac{\varepsilon
      t}{2}}|\varphi|_\infty\|\partial^2_t \psi\|_\alpha,
  \end{align*}
  where
  $C_j = 2 |\int_\RR \frac{e^{ibt}}{(-\varepsilon/2 + ib)^j}|\, db$,
  for $j = 1, 2, 3$. Thus, there exists a constant $C > 0$ such that
  $|\rho_{\varphi, \psi}(t)| \leq C e^{-\frac{\varepsilon
      t}{2}}|\varphi|_\infty\|\psi\|_{\alpha, 2}$ as required.
\end{proof}
In what follows we prove Lemma
\ref{lemma:key-estimate-laplace-transform}.  Given
$\varphi, \psi \in L^\infty(\Delta^r)$ and $s \in \CC$ define:
$ \varphi_s(x) = \int^{r(x)}_0 e^{-su} \varphi(x,u))\, du \text{ and }
\psi_s(x) = \int^{r(x)}_0 e^{su}\psi(x,u)\, du, $ for all
$x \in \Delta$.  From \cite[Appendix A]{ArMel16} we write the
correlation function as
$\rho_{\varphi, \psi}(t) = \sum^{\infty}_{n = 0}J_n(t)$ (see
\cite[Proposition A.1]{ArMel16}) and we obtain the next properties
for the Laplace Transform of $J_n$:
\begin{description}
\item[Proposition A.2 of~\cite{ArMel16}]
  $|\widehat{J}_0(s)| \lesssim |\varphi|_\infty|\psi|_\infty$;
  \item[Proposition A.3 of~\cite{ArMel16}]
    $\widehat{J}_n(s) = (\int_\Delta r\, d\mu_F)^{-1}\int_\Delta
    (\varphi_s \circ F^n)e^{-sr_n}\psi_s\, d\mu_F, \forall n \geq 1$.
  \end{description}
  Because
  $\widehat{\rho}_{\varphi, \psi}(s) = \sum^\infty_{n =
    0}\widehat{J}_n(s)$ it follows that we are left to prove the
  estimate of Lemma \ref{lemma:key-estimate-laplace-transform} for
  $\Psi(s):=\sum^\infty_{n = 1}\widehat{J}_n(s)$.

Let $A$ and $D'$ be as in Theorem
\ref{theorem:contraction-non-normalized-transfer-operator}. We split
the proof into three ranges of $n$ and $b$: (i) $|b| \leq D'$, (ii)
$n \leq A \log|b|$, $|b| \geq 2$ and (iii)
$|b| \geq D', n \geq A \log|b|$. Lemma
\ref{lemma:key-estimate-laplace-transform} follows from Lemmas
\ref{lemma:lemma-2-19-ArMel16}, \ref{lemma:lemma-2-22-ArMel16} and
\ref{lemma:lemma-2-23-ArMel16}.  The next two results follow in the
same way as in \cite{ArMel16} since they do not depend on the Hölder
norm.

\begin{proposition}\cite[Proposition 2.18]{ArMel16} \label{prop:proposition-2-18-ArMel16}
  $r e^{\varepsilon r / 2} \leq 2 \varepsilon^{-1}e^{\varepsilon r}$
  and $\int_{\Delta}e^{\varepsilon r} d\mu_F < \infty.$
\end{proposition}

\begin{lemma}[The range $n \leq A \log|b|$, $|b| \geq 2$]\cite[Lemma 2.19]{ArMel16}
   \label{lemma:lemma-2-19-ArMel16}
  There exist $\varepsilon > 0$ and $C > 0$ so that
  $\sum_{1 \leq n \leq A \log|b|}|\widehat{J}_n(s)| \leq C \varepsilon^{-2}(1 +
  |b|^{1 / 2})|\varphi|_\infty |\psi|_\infty,$ for all $\varphi, \psi \in
  L^\infty(\Delta^r)$ and for all $s = \sigma + i b$ with $\sigma \in
  [-\varepsilon / 2, 0]$ and $|b| \geq 2$.
\end{lemma}


\begin{proposition}\label{prop:reduction-observable-one-dimension}
  If $\psi \in C^\alpha_\loc(\Delta^r)$ and $s = \sigma + i b$ with
  $\sigma \leq 0$, then there exists $C
  > 0$ such that $\|\psi_s \circ h\|_\alpha \leq (|r \circ h|_\infty +
  C)\|\psi\|_\alpha$, for all $h \in \cH$.
\end{proposition}

\begin{proof}
  Given any $h \in \cH$ and $x, y \in \Delta$ without loss fo generality we may
  assume that $r(hy) \leq r(hx)$. We have
  $$
  \begin{array}{rcl}
    |\psi_s(hx) - \psi_s(hy)| & = & \left| \displaystyle \int^{r(hx)}_0 e^{su}\psi(hx, u)\,
                                    du - \int^{r(hy)}_0 e^{su} \psi(hy,u)\, du\right|\\
                              & \leq & \displaystyle\left| \int^{r(hy)}_0 e^{su}(\psi(hx,u) - \psi(hy,
                                       u))\, du \right| + \left| \int^{r(hx)}_{r(hy)}e^{su} \psi(hx, u)\, du \right|\\
                              & \leq & |r \circ h|_\infty |\psi|_{\alpha, \loc}|x - y|^\alpha + |r(hx) - r(hy)||\psi|_\infty\\
                              & \leq & |r \circ h|_\infty |\psi|_{\alpha, \loc}|x - y|^\alpha + C|\psi|_\infty|x -
                                       y|,
  \end{array}
  $$ where the constant $C$ above comes from \eqref{item:derivative-r-h} in
  Subsection \ref{subsection:uniformly-expanding-maps}. Hence
  $|\psi_s \circ h|_{\alpha} \leq |r \circ h|_\infty |\psi|_{\alpha,
    \loc} + C|\psi|_\infty.$ Also, from the definition of $\psi_s$ we
  get that
  $|\psi_s \circ h|_\infty \leq |r \circ h|_\infty
  |\psi|_\infty$. Thus, $\|\psi_s \circ h\|_\alpha = |\psi_s \circ h|_\infty + |\psi_s \circ
  h|_{\alpha,\loc}$ is bounded above by
  \begin{align*}
   (|r \circ h|_\infty +
  C) |\psi|_\infty + |r \circ h|_\infty|\psi|_{\alpha, \loc} \leq (|r \circ h|_\infty + C)
  \|\psi\|_\alpha
  \end{align*}
  as claimed.
\end{proof}

\begin{proposition}\label{prop:proposition-2-21-ArMel16}
  There exists a constant $C > 0$ such that
  $$|\varphi_s|_1 \leq C \varepsilon^{-1}|\varphi|_\infty, \ \
  |P_s(f_0\psi_s)|_\infty \leq C \varepsilon^{-1}|\psi|_\infty,\ \ \|P_s(f_0
  \psi_s)\|_\alpha \leq C \varepsilon^{-1}\|\psi\|_\alpha,$$ for all $\varphi \in
  L^\infty(\Delta^r)$, $\psi \in C^\alpha_\loc(\Delta^r)$ and all $s = \sigma +
  i b$ with $\sigma \in [-\varepsilon / 2, 0]$.
\end{proposition}

\begin{proof}
  By Proposition \ref{prop:proposition-2-18-ArMel16},
  $|\varphi_s|_1 =\int_\Delta |\varphi_s|\, d\mu_F$ is bounded above
  by
  \begin{align*}
    \int_\Delta \int^{r(x)}_0 e^{\varepsilon r(x) / 2}|\varphi(x, u)|\, du\, d\mu_F(x)
    \leq  |\varphi|_\infty \int_\Delta re^{\varepsilon
    r / 2}\, d\mu_F \lesssim \varepsilon^{-1}|\varphi|_\infty.
  \end{align*}
  It follows from the definition of $\psi_s$ that $|\psi_s \circ h|_\infty \leq
  |r \circ h|_\infty |\psi|_\infty,$ for all $h \in \cH$ and all $s = \sigma + i
  b$ with $\sigma \leq 0$. The Proposition \ref{prop:proposition-2-18-ArMel16}
  ensures that $|r \circ h|_\infty e^{\varepsilon |r \circ h|_\infty / 2} \leq 2
  \varepsilon^{-1}e^{\varepsilon |r \circ h|_\infty}$. Using these inequalities
  we get
  $$
  \begin{array}{rcl}
    |A_{s,h}(f_0 \psi_s)|_\infty & \leq & e^{\varepsilon |r \circ
                                          h|_\infty / 2}|h'|_\infty |f_0 \circ
                                          h|_\infty |\psi_s \circ h|_\infty\\
                                 & \leq & |f_0|_\infty|\psi|_\infty e^{\varepsilon |r \circ h|_\infty / 2}|r
                                          \circ h|_\infty |h'|_\infty\\
                                 & \leq & 2\varepsilon^{-1}|f_0||\psi|_\infty e^{|r \circ h|_\infty}|h'|_\infty.
  \end{array}
  $$ By condition \eqref{item:expontial-tail} of Subsection \ref{subsection:uniformly-expanding-maps}
  $|P_s(f_0 \psi_s)|_\infty \leq 2\varepsilon^{-1}|f_0|_\infty
  |\psi|_\infty \sum_{h \in \cH} e^{\varepsilon |r \circ
    h|_\infty}|h'|_\infty \lesssim \varepsilon^{-1}|\psi|_\infty.$
  Finally, it follows from the proof of Proposition
  \ref{prop:twisted-operator} that, for all $h \in \cH$
  $$|A_{s,h}(f_0\psi_s)|_\alpha \lesssim e^{\varepsilon |r \circ h|_\infty / 2}
  |h'|_\infty \|\psi_s \circ h\|_\alpha.$$ By Propositions \ref{prop:proposition-2-18-ArMel16} and
  \ref{prop:reduction-observable-one-dimension}
  $$
  \begin{array}{rcl}
    |A_{s,h}(f_0 \psi_s)|_\alpha & \leq & e^{\varepsilon |r \circ h|_\infty / 2}|h'|_\infty(|r \circ h|_\infty + C)\|\psi\|_\alpha \\
     & \leq & e^{\varepsilon |r \circ h|_\infty / 2}|r \circ h|_\infty
              |h'|_\infty(1 + C / \inf r)\|\psi\|_\alpha \\
                                 & \lesssim &  \varepsilon^{-1}e^{\varepsilon |r \circ h|_\infty} |h'|_\infty\|\psi\|_\alpha.
  \end{array}
  $$ Again by condition \eqref{item:expontial-tail} it follows that
  $|P_s(f_0\psi_s)|_\alpha \lesssim \varepsilon^{-1}\|\psi\|_\alpha$ and also
  $\|P_s(f_0\psi_s)\|_\alpha \lesssim \varepsilon^{-1}\|\psi\|_\alpha$.
\end{proof}

For the next lemma we define the family of normalized twisted transfer operators $\cQ_s: C^\alpha(\Delta)
\to C^\alpha(\Delta)$ by $\cQ_s\psi = f^{-1}_0P_s(f_0\psi)$. Note that
$\int_\Delta \varphi \cQ_s \psi\, d\mu_F = \int_\Delta (\varphi \circ F)e^{-s
  r}\psi\, d\mu_F$, for all $\varphi \in L^\infty(\Delta)$ and $\psi \in C^\alpha(\Delta)$.

\begin{lemma}[The range $|b| \leq D'$]\label{lemma:lemma-2-22-ArMel16}
  There exist constants $\varepsilon > 0$ and $C > 0$ such that
  $|\Psi(s)| \leq C \varepsilon^{-1} |\varphi|_\infty\|\psi\|_\alpha$
  for all $\varphi \in L^\infty(\Delta^r)$,
  $\psi \in C^\alpha_\loc(\Delta^r)$ and for all $s = \sigma + ib$
  with $\sigma \in [-\varepsilon / 2, 0]$ and $|b| \leq D'$.
\end{lemma}

\begin{proof}
  Replacing $\psi$ by $\psi - \int_{\Delta^r} \psi\, d\mu^r_F$, we can suppose
  without loss of generality that $\psi$ lies on $\cB = \{\psi \in
  C^\alpha_\loc(\Delta^r):\ \int_{\Delta^r}\psi\, d\mu^r_F = 0\}$. We can
  write
  $$
  \begin{array}{rcl}
    \Psi(s) & = & \displaystyle\sum^\infty_{n =
                  1}\int_\Delta \varphi_s \cQ^{n}_s\psi_s\, d\mu_F
                  = \int_\Delta \varphi_s(1 -
                  \cQ_s)^{-1}Q_s\psi_s\, d\mu_F
                  =  \int_\Delta \varphi_sZ_s\psi\, d\mu_F,
  \end{array}
  $$ where $Z_s\psi := (1 - \cQ_s)^{-1}\cQ_s \psi_s$. Using the expression of
  $\cQ_s$ and Remark \ref{remark:decay-local-holder} it follows that
  $\cQ_s\psi_s \in C^\alpha(\Delta)$ for all
  $\psi \in C^\alpha_\loc(\Delta)$.  In particular
  $Z_s \psi \in C^\alpha(\Delta)$, for all
  $\psi \in C^\alpha_\loc(\Delta^r)$. Following \cite[Lemma
  2.22]{ArMel16} it is possible to prove that the family of operators
  $Z_s: \cB \to C^\alpha(\Delta)$ is analytic on
  $\{s \in \CC:\ \Re(s) > 0\}$ and admits an extension beyond the
  imaginary axis. In particular, there exists $\varepsilon > 0$ such
  that $Z_s$ is analytic on the region
  $[-\varepsilon, 0] \times [-D', D']$ and hence there is a constant
  $C > 0$ such that $\|Z_s\psi\|_\alpha \leq C \|\psi\|_\alpha$ for
  all $\psi \in \cB$ and all
  $s \in [-\varepsilon, 0] \times [-D', D]$. Thus it follows from
  Proposition \ref{prop:proposition-2-21-ArMel16} that
  $|\Psi(s)| \leq |\varphi_s|_1|Z_s\psi|_\infty \lesssim
  \varepsilon^{-1}|\varphi|_\infty\|\psi\|_\alpha.$
\end{proof}

\begin{lemma}[The range $|b| \leq D'$, $n \geq A \log|b|$]\label{lemma:lemma-2-23-ArMel16}
  There exist constants $\varepsilon, C > 0$ such that
  $\sum_{n \geq A \log|b|}|\widehat{J}(s)| \leq C
  \varepsilon^{-2}|\varphi|_\infty \|\psi\|_\alpha,$ for all
  $\varphi \in L^\infty(\Delta^r)$, $\psi \in C^\alpha_\loc(\Delta^r)$
  and for all $s = \sigma + i b$ with
  $\sigma \in [-\varepsilon / 2, 0]$, $|b| \leq D'$.
\end{lemma}

\begin{proof}
  For short, in what follows we denote
  $\overline{r} = \int_\Delta r\, d\mu_F$.  Let $\cQ_s$ be as in Lemma
  \ref{lemma:lemma-2-22-ArMel16}. Note that
  $\widehat{J}_n(s) = \overline{r}^{-1}\int_Y \varphi_s
  \cQ^n_s\psi_s\, d\mu_F = \overline{r}^{-1}\int_\Delta \varphi_s
  f^{-1}_0 P^n_s(f_0 \psi_s)\, d\mu_F.$ Hence, using Theorem
  \ref{theorem:contraction-non-normalized-transfer-operator} and
  Proposition \ref{prop:proposition-2-21-ArMel16}

  $$
  \begin{array}{rcl} \displaystyle \sum_{n \geq A \log|b|}|\widehat{J}_n(s)| &
                                                                               \lesssim & \displaystyle \varepsilon^{-1}|\varphi|_\infty \sum_{n \geq A \log|b|}|P^{n - 1}_s
                                                                                          (P_s(f_0 \psi_s))|_\infty \\
                                                                             & \lesssim & \displaystyle \varepsilon^{-1}|\varphi|_\infty \sum_{n \geq A \log|b|}\|P^{n - 1}_s\|_b \|P_s(f_0\psi_s)\|_b\\
                                                                             & \lesssim &  \varepsilon^{-1}|\varphi|_\infty \|P_s(f_0\psi_s)\|_b \lesssim
                                                                                          \varepsilon^{-2}|\varphi|_\infty \|\psi\|_\alpha,
  \end{array}
  $$ as required.
\end{proof}

For the next result we introduce the Lebesgue measure $\Leb^r_2$ on $\Delta^r$ by
setting $\Leb^r_2 = (\Leb_\Delta \times \Leb_\RR) / \int r\, d\Leb_\Delta$, where
$\Leb_\Delta$ is the Lebesgue measure restricted to the Borelean sets of $\Delta$.

\begin{corollary}[Convergence to equilibrium]
  \label{corollary:convergence-equilibrium-onedimensional}
  In the setting of Theorem
  \ref{theorem:decay-onedimensional-weaker}, there are constants
  $c, C > 0$ so that for all
  $\varphi \in L^\infty(\Delta^r),$
  $\psi \in C^{\alpha, 2}_\loc(\Delta^r)$, $t > 0$ we have
  $\left| \int (\varphi \circ F_t) \psi\, d\Leb^r_2 - \int \varphi\,
    d\mu^r_F \int \psi\, d\Leb^r_2 \right| \leq
  Ce^{-ct}|\varphi|_\infty\|\psi\|_{\alpha, 2}.$ 
\end{corollary}

\begin{proof}
  Since $\mu_F$ is absolutely continuous with respect to $\Leb_\Delta$, it
  follows that $\mu^r_{F}$ is also absolutely continuous with respect to
  $\Leb^r_2$. Moreover, the density satisfies
  $$\frac{d\mu^r_F}{d\Leb^r_2}(x,u) = \frac{d\mu_F}{d\Leb_\Delta}(x),$$ for all
  $(x,u) \in \Delta^r$. Because $d\mu_F / d\Leb_\Delta$ is
  $\alpha$-Hölder and bounded from above and below, it follows that
  $\xi := d\Leb_\Delta / d\mu_F$ is also $\alpha$-Hölder and bounded
  from above and below. Hence we have that
  $\psi \xi \in C^{\alpha, 2}_{\loc}(\Delta^r)$. Finally, using
  Theorem \ref{theorem:decay-onedimensional-weaker} we get
  \begin{align*}
    \Big| \int (\varphi \circ F_t)\psi\, d\Leb^r_2
    &-\!\!
     \int \varphi\, d\mu^r_F\int \psi\, d\Leb^r_2 \Big|
      = \Big|
      \int (\varphi \circ
      F_t)\psi\xi\, d\mu^r_F
      -\!\!
      \int\varphi\, d\mu^r_F\int\psi\xi d\mu^r_F
      \Big|
    \\
    & \leq  C |\varphi|_\infty
      \|\psi \xi\|_{\alpha,
      2}e^{-c t} \leq C|\xi|_\infty
      |\varphi|_\infty
      \|\psi\|_{\alpha,
      2}e^{-c t}
  \end{align*}
  as stated.
\end{proof}

\subsection[Decay of correlations for $C^{1 + \alpha}$ hyperbolic skew
product semiflows]{Decay of correlations for {\boldmath
    $C^{1 + \alpha}$} hyperbolic skew product semiflows}
\label{section:decay-hyperbolic-skew-product-semiflows}

In this section we prove Theorem
\ref{prop:decay-correlations-skew-product-weaker}. Let $\widehat{F}_t:
\widehat{\Delta}^r \to \widehat{\Delta}^r$ be a $C^{1 + \alpha}$ hyperbolic skew
product semiflow with roof function $r$ satisfying the UNI condition as in Subsection \ref{subsec:decay-hyperbolic-skew-product-flow-statement}.

\subsubsection{Disintegration of the measure $\mu_{\widehat{F}}$}
Let $\mathcal{L}: L^1(\Delta) \to L^1(\Delta)$ be the \emph{transfer
  operator} for the map $F$, that is, $\mathcal{L}$ satisfies
$\int (\varphi \circ F) \cdot \psi\, d\mu_F = \int \varphi \cdot
\mathcal{L}\psi\, d\mu_F$ for all $\varphi \in L^\infty(\Delta)$ and
$\psi \in L^1(\Delta)$. It is clear that $\cL := \cL_0$ is a transfer
operator.

In order to deduce Theorem
\ref{prop:decay-correlations-skew-product-weaker} we need the
following properties.  Recall that we denote
$\widehat{F}^n(x,y) = (F^n(x), G_n(x,y))$, for all
$(x,y) \in \widehat{\Delta}$, and for a function
$\psi \in C^0(\widehat{\Delta})$ we denote
$\psi_n(x) = \psi(\widehat{F}^n(x, 0))$, for all $x \in \Delta$, where
$0$ is some element of $\Omega$ that we fix previously.

\begin{proposition}\label{prop:3-6-butterley-melbourne}
  The following properties hold
  \begin{enumerate}[(a)]
  \item For each $\psi \in C^0(\widehat{\Delta})$ the limit
    $\eta_x(\psi) = \lim_{n \to +\infty}(\mathcal{L}^n\psi_n)(x)$
    exists for $\mu_F$-almost every $x \in \Delta$ and defines a
    probability measure supported on $\pi^{-1}(x)$. Moreover, the
    function $\overline{\psi}: \Delta \to \RR$, given by
    $\overline{\psi}(x) = \eta_x(\psi) := \int_{\pi^{-1}(x)} \psi \,
    d\eta_x$, is $\mu_F$-integrable and
    $\int_{\widehat{\Delta}} \psi\, d\mu_{\widehat{F}} = \int_{x \in
      \Delta}\int_{\pi^{-1}(x)}\psi\, d\eta_x\, d\mu_F(x).$
    \item For any $\psi \in C^\alpha_\loc(\widehat{\Delta})$, the function
$\overline{\psi}: \Delta \to \RR$, defined as item (a), lies on $C^\alpha_\loc(\Delta)$ and there exists a
constant $C > 0$ such that $\|\overline{\psi}\|_\alpha \leq C\|\psi\|_\alpha$.
  \end{enumerate}
\end{proposition}

The item (a) of this proposition is \cite[Proposition 3]{BM2015} and
item (b) is a generalization of \cite[Proposition 6]{BM2015} to the
function space $C^\alpha_\loc(\widehat{\Delta})$. Checking that paper
we see that the only result that need to be generalized there is
\cite[Corollary 8]{BM2015} and we do this in Lemma
\ref{lemma:gen-corollary-8-butterley-melbourne}. Let
$\xi \in C^\alpha(\Delta)$ be the density of $\mu_F$ with respect to
the Lebesgue measure in $\Delta$. Because $F$ is uniformly expanding
(with full branches), we have that
$\cL^n\psi_n = \xi^{-1}\sum_{h \in \cH_n}h'(\xi \psi_n) \circ
h$. Recall that from item (a) of Proposition
\ref{prop:3-6-butterley-melbourne} we have that
$\overline{\psi} = \lim_{n \to \infty} \cL^n\psi_n$.

\begin{lemma}\cite[Lemma 7]{BM2015}\label{lemma:7-butterley-melbourne}
  There exists $C > 0$ such that $|G_n(hx_1, y) - G_n(hx_2, y)| \leq C|x_1 -
  x_2|$, for all $h \in \cH_n$, $n \geq 1$, $(x_1,y),
  (x_2,y) \in \widehat{\Delta}$. 
\end{lemma}

Now we generalize Corollary 8 in \cite{BM2015}.

\begin{lemma}\label{lemma:gen-corollary-8-butterley-melbourne} For each $h \in
  \cH_n$ and $\psi \in C^\alpha_\loc(\widehat{\Delta})$ we have $(\xi \psi_n)
  \circ h \in C^\alpha_\loc(\Delta)$. More precisely, there exists $C > 0$ such that $\|(\xi \psi_n) \circ h\|_{\alpha} \leq
  C\|\psi\|_{\alpha}$, for all $\psi \in
  C^\alpha_{\loc}(\widehat{\Delta}),\ h \in \cH_n$ and $n \geq 1.$
\end{lemma}

\begin{proof}
  Let $x_1, x_2 \in \Delta$ and $g \in \cH$. Since $\widehat{F}^n(hx, 0) = (x,
  G_n(hx, 0))$ we have that $\psi_n \circ h(x) = \psi(x, G_n(hx, 0))$, for
  all $x \in \Delta$. Hence, using Lemma \ref{lemma:7-butterley-melbourne} and
  \eqref{item:backward-contraction} from Subsection
  \ref{subsection:uniformly-expanding-maps} we get that $|\psi_n \circ h(gx_1) - \psi_n \circ
    h(gx_2)|$ is bounded above by
    \begin{align*}
      |\psi|_{\alpha, \loc}(|gx_1 - gx_2|^\alpha + |G_n(h(gx_1), 0) -
      G_n(h(gx_2), 0)|)
      \lesssim  |\psi|_{\alpha, \loc}|x_1 - x_2|^\alpha
    \end{align*}
    and we get that
    $|\psi_n \circ h \circ g|_\alpha \lesssim |\psi|_{\alpha,
      \loc}$. Taking the supremum over $g \in \cH$ we get that
    $|\psi_n \circ h|_{\alpha, \loc} \lesssim|\psi|_{\alpha,
      \loc}$. Thus, we obtain
    $\|\psi_n \circ h\|_\alpha = |\psi_n \circ h|_\infty + |\psi_n
    \circ h|_{\alpha, loc}\lesssim \|\psi\|_\alpha$.

    For the density $\xi$ we have that
    $|\xi \circ h|_\infty \leq |\xi|_\infty$ and
  \begin{align*}
  |\xi \circ h(x_1) - \xi \circ h(x_2)|
  & \leq
  |\xi|_\alpha|hx_1 -  hx_2|^\alpha
  \leq |\xi|_\alpha |h'|^\alpha_\infty|x_1 - x_2|^\alpha
  \lesssim
  |\xi|_\alpha  \rho^{\alpha n}|x_1 - x_2|^\alpha.
  \end{align*}
  Therefore, there exists a constant $C > 1$ such that $|\xi \circ h|_\alpha
  \leq C$ and we get that
  $\|\xi \circ h\|_\alpha = |\xi \circ h|_\infty + |\xi \circ
  h|_\alpha \leq 2C$.  Finally,
  $\|(\xi\psi_n) \circ h\|_\alpha \leq \|\xi \circ h\|_\alpha
  \|\psi_n \circ h\|_\alpha \leq 2C \|\psi\|_\alpha$ as needed.
\end{proof}

With Lemma \ref{lemma:gen-corollary-8-butterley-melbourne} the proof
of the next result follows as in \cite{BM2015} even for the class of
observables $C^\alpha_{\loc}(\widehat{\Delta})$.


\begin{lemma}\cite[Lemma 9]{BM2015}\label{lemma:9-butterley-melbourne}
  There exists $C > 0$ such that $\|\cL^n\psi_n\|_\alpha \leq
  C\|\psi\|_\alpha$, for all $\psi \in C^\alpha_{\loc}(\widehat{\Delta})$
  and $n \geq 1$.
\end{lemma}


\begin{proof}[Proof of Proposition \ref{prop:3-6-butterley-melbourne}
  (b)] Because $\cL^n\psi_n$ converges pointwise to $\overline{\psi}$
  and $\sup_n\|\cL^n\psi_n\|_\alpha < \infty$ it follows that
  $\overline{\psi} \in C^\alpha_\loc(\Delta)$. Now using Lemma
  \ref{lemma:9-butterley-melbourne} we get the desired result.
\end{proof}

Define $w_t:\Delta^r \to \RR$ to be the number of visits to $\Delta$
by time $t$, that is,
$$w_t(x,u) = \max\{n \geq 0:\ u + t \geq r_n(x)\}.$$ For the next
proposition, we recall that we are denoting by $\gamma$ the
contraction rate along $\Omega$ for the skew product $\widehat{F}$
(check Subsection \ref{subsection:hyperbolic-skew-products}).

\begin{proposition}\cite[Proposition 3.5]{ArMel16}\label{prop:3.5-ArMel16}
  There exist $\delta, C > 0$ such that $\int_{\Delta^r}\gamma^{\alpha
    w_t}d\mu^r_F \leq C e^{-\delta t}$, for all $t > 0$.
\end{proposition}
  
Now we are ready to prove Theorem
\ref{prop:decay-correlations-skew-product-weaker}. We follow the
approach of \cite[Theorem 3.3]{ArMel16}. Here we denote the
correlation function of $\widehat{F}_t$ with respect to observables
$\varphi \in L^\infty(\widehat{\Delta}^r)$ and
$\psi \in L^1(\widehat{\Delta}^r)$ by
$\PP_{\varphi, \psi}(t) = \int_{\widehat{\Delta}^r}(\varphi \circ
\widehat{F}_t)\cdot \psi\, d\mu^r_{\widehat{F}} -
\int_{\widehat{\Delta}^r} \varphi\,
d\mu^r_{\widehat{F}}\int_{\widehat{\Delta}^r} \psi\,
d\mu^r_{\widehat{F}}.$

\begin{proof}[Proof of Theorem \ref{prop:decay-correlations-skew-product-weaker}]
  Let $\varphi \in C_\loc^\alpha(\widehat{\Delta}^r)$ and
  $\psi \in C_\loc^{\alpha, 2}(\widehat{\Delta}^r)$. Without loss of
  generality, we may suppose that
  $\int_{\widehat{\Delta}^r}\psi \, d\mu^r_{\widehat{F}} = 0$. Define
  $\pi^r: \widehat{\Delta}^r \to \Delta^r$ as
  $\pi^r(x,u) = (\pi x, u)$. This defines a semiconjugacy between
  $F_t$ and $\widehat{F}_t$, that is,
  $\pi^r \circ \widehat{F}_t = F_t \circ \pi^r$ and
  $\pi^r_*\mu^r_{\widehat{F}} = \mu^r_F$.  Define
  $\varphi_t: \Delta^r \to \RR$ by setting
  $\varphi_t(y,u) = \int_{x \in \pi^{-1}(y)}\varphi \circ
  \widehat{F}_t(x,u)\, d\eta_y(x).$ Then
  $\PP_{\varphi, \psi}(2t) = I_1(t) + I_2(t)$, where
  $I_2(t) = \int_{\widehat{\Delta}^r} (\varphi_t \circ F_t \circ
  \pi^r) \psi\, d\mu^r_{\widehat{F}}$ and
  \begin{align*}
I_1(t) = \int_{\widehat{\Delta}^r}(\varphi \circ \widehat{F}_{2t})\psi\, d\mu^r_{\widehat{F}} - \int_{\widehat{\Delta}^r}
  (\varphi_t \circ F_t \circ \pi^r)\psi\, d\mu^r_{\widehat{F}}.
  \end{align*}
   Note that $I_1(t) = \int_{\widehat{\Delta}^r} [(\varphi \circ \widehat{F}_t
  - \varphi_t \circ \pi^r) \circ \widehat{F}_t]\psi\, d\mu^r_{\widehat{F}}$ and so
  $$|I_1(t)| \leq |\psi|_\infty\int_{\widehat{\Delta}^r}|\varphi \circ \widehat{F}_t - \varphi_t \circ
  \pi^r|\, d\mu^r_{\widehat{F}}.$$ Using definitions of $\pi^r$ and
  $\varphi_t$ we get that
  \begin{align*}
    \varphi \circ \widehat{F}_t(x,u) - \varphi_t \circ \pi^r(x,u)
    = 
    \int_{x' \in \pi^{-1}(\pi x)}(\varphi \circ \widehat{F}_t(x,u) -
    \varphi \circ \widehat{F}_t(x',u))\, d\eta_{\pi(x)}(x').
  \end{align*}
  Using the contraction of $\widehat{F}$ along $\Omega$ we have that
  there exists a constat $C > 0$ such that
  $|\varphi \circ \widehat{F}_t(x,u) - \varphi_t \circ \pi^r(x,u)|$ is
  bounded above by
  \begin{align*}
  C \int_{x' \in \pi^{-1}(\pi x)}|\varphi|_{\alpha,
    \loc}\gamma^{\alpha w_t(\pi x, u)}\, d\eta_{\pi(x)}(x')
  =
  C|\varphi|_{\alpha, \loc}\gamma^{\alpha w_t(\pi x, u)} =
  C|\varphi|_{\alpha, \loc}\gamma^{\alpha w_t} \circ \pi^r(x, u).
  \end{align*}
  Hence
  $|I_1(t)| \leq C|\psi|_\infty |\varphi|_{\alpha, \loc}\int \gamma^{\alpha w_t}
  \circ \pi^r\, d\mu^r_{\widehat{F}} = C|\psi|_\infty |\varphi|_{\alpha,
    \loc}\int \gamma^{\alpha w_t}\, d\mu^r_F.$ Now, using Proposition
  \ref{prop:3.5-ArMel16}, $|I_1(t)| \leq C|\psi|_\infty |\varphi|_{\alpha,
    \loc}e^{-\delta t}$, for some $\delta > 0$ and for all $t > 0$.

  Now define $\overline{\psi}: \Delta^r \to \RR$ by setting
  $\overline{\psi}(x,u) = \int_{z \in \pi^{-1}(x)}\psi(z,u)\, d\eta_x(z)$.
  Since $\int_{\widehat{\Delta}^r} \psi\, d\mu^r_{\widehat{F}} = 0$, it follows from item (a) of
  Proposition \ref{prop:3-6-butterley-melbourne} that $\int_{\Delta^r} \overline{\psi}\,
  d\mu^r_F = 0$.

  We also have that
  $I_2(t) = \int_{\Delta^r} (\varphi_t \circ F_t)\overline{\psi}\,
  d\mu^r_F = \rho_{\varphi_t, \overline{\psi}}(t)$, where $\rho$
  denotes the correlation function for $F_t$.
  
  By Proposition \ref{prop:3-6-butterley-melbourne} (b), we
  have that $\overline{\psi} \in C^{\alpha, 2}_\loc(\Delta)$ and
  $\|\overline{\psi}\|_{\alpha, 2} \leq C \|\psi\|_{\alpha,2}$, for some
  constant $C > 0$. Hence, using Theorem \ref{theorem:decay-onedimensional-weaker} there exist constants $c,
  C>0$ such that $|I_2(t)| \leq
  Ce^{-ct}\|\overline{\psi}\|_{\alpha,2}|\varphi_t|_\infty \leq
  Ce^{-ct}\|\psi\|_{\alpha, 2}|\varphi|_\infty$ completing the proof.
\end{proof}

\subsubsection{Convergence to equilibrium for hyperbolic skew product semiflows}
\label{sub:convergence-equilibrium-semiflow}

In this subsection we prove the convergence to equilibrium for the hyperbolic
skew product semiflow $\widehat{F}_t$, that is, we prove Corollary \ref{corollary:convergence-to-equilibrium-skew-product-flow}.

\begin{proof}[Proof of Corollary \ref{corollary:convergence-to-equilibrium-skew-product-flow}]
  Recall from the proof of Corollary
  \ref{corollary:convergence-equilibrium-onedimensional} that
  $\xi = d\Leb^r_2 / d\mu^r_F$ is well-defined and bounded from above
  and below. In particular, it follows from Proposition
  \ref{prop:3.5-ArMel16} that
  $\int_{\Delta^r} \gamma^{\alpha w_t}d\Leb^r_2 \leq
  |\xi|_\infty\int_{\Delta^r} \gamma^{\alpha w_t}\, d\mu^r_F \lesssim
  e^{-\delta t}, $ for some $\delta > 0$ and for all $t > 0$.

  If, in the proof of Theorem
  \ref{prop:decay-correlations-skew-product-weaker}, we define
  $I_1(t)$ and $I_2(t)$ with the measure $\Leb^r_3$ instead of
  $\mu^r_{\widehat{F}}$, then (i) $\pi^r_*\Leb^r_3 = \Leb^r_2$ and we
  can use (ii) Corollary
  \ref{corollary:convergence-equilibrium-onedimensional} instead of
  Theorem \ref{theorem:decay-onedimensional-weaker}; and (iii) the
  above inequality for the integral of the number of visits to
  $\Delta$ instead of Proposition \ref{prop:3.5-ArMel16}, to redo all
  the calculations in the same way to get the desired result.
\end{proof}



 \def\cprime{$'$}

\bibliographystyle{abbrv}
 

\end{document}